\documentclass[a4paper,11pt]{article}
\usepackage{amsmath,amsthm,amssymb,dsfont,graphicx,xspace,epsfig,xcolor}
\usepackage[T1]{fontenc}
\usepackage[plain]{fullpage}
\usepackage{thm-restate}
\usepackage[hidelinks]{hyperref}
\usepackage{bm}
\usepackage{color}
\usepackage{comment}
\usepackage{mathrsfs}
\usepackage[shortlabels]{enumitem}
\usepackage{float}
\usepackage{dutchcal}
\usepackage{pict2e}
\usepackage{mathtools}
\usepackage{tikz}
\usepackage{subfigure}
\usepackage{xstring,refcount}
\usepackage{authblk}

\usetikzlibrary{fit, backgrounds, matrix, arrows.meta}
\usetikzlibrary{calc,arrows,positioning}
\usetikzlibrary{decorations}
\usetikzlibrary{decorations.pathmorphing}
\usetikzlibrary{decorations.pathreplacing}
\usetikzlibrary{decorations.shapes}
\usetikzlibrary{decorations.text}
\usetikzlibrary{decorations.markings}
\usetikzlibrary{decorations.fractals}
\usetikzlibrary{decorations.footprints}

\tikzset{vertex/.style = {circle,fill=black,minimum size=4pt, inner sep=0pt}}
\tikzset{edge/.style = {thick,dashed,->,> = latex,dash pattern={on 5pt off 4pt}}}
\tikzset{arc/.style = {thick,->,> = latex}}
\tikzset{digon/.style = {thick,<->,> = latex}}
\tikzset{bigvertex/.style = {shape=circle,draw}}

\definecolor{g-blue}{rgb}{0.0, 0.5, 1.0}

\newtheorem{theorem}{Theorem}
\newtheorem{lemma}[theorem]{Lemma}
\newtheorem{corollary}[theorem]{Corollary}

\newtheorem{conjecture}[theorem]{Conjecture}

\newtheorem{claim}{Claim}[theorem]
\newcounter{oldthm}{}
\def\@extractthmnum#1.#2{#2}
\newenvironment{claimlabeled}[1][\@nil]{%
    \setcounter{oldthm}{\value{theorem}}%
    \def\tmp{#1}%
    \ifx\tmp\@nnil
    \else
        \refused{#1}
        \setcounter{theorem}{\getrefnumber{#1}}
    \fi
    \begin{claim}%
}{%
    \ifx\tmp\@nnil
    \else
        \setcounter{theorem}{\value{oldthm}}
    \fi
    \end{claim}
}

\newtheorem{subclaim}{Subclaim}[claim]
\newtheorem{remark}[theorem]{Remark}

\theoremstyle{definition}

\newtheorem{problem}[theorem]{Problem}

\newenvironment{proofclaim}[1][]{\par\noindent {\it Proof of claim}. }{ \hfill$\lozenge$\par\addvspace{6pt plus 6pt}}

\newcommand{\opentriangle}{
  \raisebox{0.2pt}{\makebox[0.77778em]{
    \setlength{\unitlength}{0.6em}
    \linethickness{0.4pt}
    \begin{picture}(1,1)
    \polygon(0,0)(1,0)(1,1)
    \end{picture}
  }}
}
\newenvironment{proofsubclaim}[1][]
	{\par\noindent {\it Proof of subclaim}. }{ \hfill$\opentriangle$\par\addvspace{6pt plus 6pt}}

\newcommand{\dic}{\vec{\chi}}
\makeatletter
\newcommand{\overleftrightsmallarrow}{\mathpalette{\overarrowsmall@\leftrightarrowfill@}}
\newcommand{\overrightsmallarrow}{\mathpalette{\overarrowsmall@\rightarrowfill@}}
\newcommand{\overleftsmallarrow}{\mathpalette{\overarrowsmall@\leftarrowfill@}}
\newcommand{\overarrowsmall@}[3]{%
  \vbox{%
    \ialign{%
      ##\crcr
      #1{\smaller@style{#2}}\crcr
      \noalign{\nointerlineskip}%
      $\m@th\hfil#2#3\hfil$\crcr
    }%
  }%
}
\def\smaller@style#1{%
  \ifx#1\displaystyle\scriptstyle\else
    \ifx#1\textstyle\scriptstyle\else
      \scriptscriptstyle
    \fi
  \fi
}
\makeatother
\newcommand{\bid}[1]{\overleftrightsmallarrow{#1}}
\newcommand{\bic}{\overset{\text{\tiny$\leftrightarrow$}}{\omega}}
\newcommand{\delmax}{\Delta_{\max}}
\newcommand{\delmin}{\Delta_{\min}}
\newcommand{\deltil}{\Tilde{\Delta}}
\newcommand{\delplus}{\Delta^+}
\newcommand{\dig}[1]{\langle #1 \rangle}
\newcommand{\inboundary}[1]{\partial^-(#1)}
\newcommand{\outboundary}[1]{\partial^+(#1)}
\newcommand{\inboundarysize}[1]{|\partial^-(#1)|}
\newcommand{\outboundarysize}[1]{|\partial^+(#1)|}
\newcommand{\EE}{\mathbb{E}}
\newcommand{\PP}{\mathbb{P}}
\newcommand{\cSCR}{\mathscr{C}}
\newcommand{\lSCR}{\mathscr{L}}
\newcommand{\zSCR}{\mathscr{Z}}
\newcommand{\sSCR}{\mathscr{S}}
\newcommand{\mSCR}{\mathscr{M}}
\newcommand{\rSCR}{\mathscr{R}}
\newcommand{\tSCR}{\mathscr{T}}
\newcommand{\hX}{X^{\star}}
\newcommand{\hS}{\sSCR^{\star}}
\newcommand{\DeltaDecomp}{\Delta_{\ref{lemma:dense_decomposition}}}
\newcommand{\DeltaDirectedBK}{\Delta_{\ref{thm:directed_bk}}}
\newcommand{\DeltaDirectedBKPlus}{\Delta_{\ref{thm:directed_bk_deltaplus}}}

\DeclareMathOperator{\UG}{UG}
\DeclareMathOperator{\ch}{ch}
\renewcommand{\epsilon}{\varepsilon}
\renewcommand{\emptyset}{\varnothing}
\renewcommand{\hat}{\widehat}
\renewcommand{\phi}{\varphi}
\let\le\leqslant

\let\leq\leqslant
\let\geq\geqslant

\hypersetup{
    colorlinks,
    linkcolor={red!50!black},
    citecolor={blue!50!black},
    urlcolor={blue!80!black}
}

\title{$(\Delta-1)$-dicolouring of digraphs}

\author[1]{Ararat Harutyunyan\thanks{Research supported by ANR-21-CE48-0012.}$^,$}
\author[2]{Ken-ichi Kawarabayashi\thanks{Research supported by JSPS KAKENHI JP20A402 and 22H05001, and by JST ASPIRE JPMJAP2302.}$^,$}
\author[2]{Lucas Picasarri-Arrieta\protect\footnotemark[2]$^,$}
\author[1]{\\Gil Puig i Surroca\protect\footnotemark[1]$^,$}

\affil[1]{LAMSADE, Université Paris Dauphine - PSL, Paris, France}
\affil[2]{National Institute of Informatics, The University of Tokyo, Tokyo, Japan.}

\date{}

\begin{document}

\maketitle

\begin{abstract}

In 1977, Borodin and Kostochka conjectured that every graph with maximum degree $\Delta \geq 9$ is $(\Delta-1)$-colourable, unless it contains a clique of size $\Delta$. In 1999, Reed confirmed the conjecture when $\Delta\geq 10^{14}$.
We propose different generalisations of this conjecture for digraphs, and prove the analogue of Reed's result for each of them.

The {\it dichromatic number} of a digraph is the minimum number of colours needed to partition it into acyclic induced subdigraphs.
A {\it biclique} is a set of vertices that are pairwise linked with two arcs in opposite directions. For a digraph $D$, we let $\deltil(D) = \max_{v\in V(D)} \sqrt{d^+(v) \cdot d^-(v)}$.

Let $\Delta\geq 9$ be any integer and $D$ be any digraph with $\deltil(D) \leq \Delta$. We conjecture that $D$ has dichromatic number at most $\Delta-1$, unless either (i) $D$ contains a biclique of size $\Delta$, or (ii) $D$ contains a biclique $K$ of size $\Delta-2$, a  directed $3$-cycle $\vec{C_3}$ disjoint from $K$, and all possible arcs in both directions between $\vec{C_3}$ and $K$. If true, this implies the conjecture of Borodin and Kostochka. We prove our conjecture when $\Delta$ is large enough, thereby generalising the result of Reed.

We also conjecture that every digraph $D$ with maximum out-degree $\Delta^+\geq 9$ has dichromatic number at most $\Delta^+-1$, unless $D$ contains one of the two obstructions described above. 
Again, we show the conjecture when $\Delta^+$ is large enough, yielding a second generalisation of Reed's result.
We finally give a sufficient condition for a digraph $D$ to have dichromatic number at most $\delmin(D)-1$, assuming that $\delmin(D) = \max_{v\in V(D)} \min (d^+(v),d^-(v))$ is large enough. 
In particular, this holds when $D$ has no biclique of size $\lceil \frac{1}{2}(\delmin(D)-1)\rceil$, and when the underlying graph of $D$ has no clique of size $\delmin(D)$, thus yielding a third independent generalisation of Reed's result. 
We further give a hardness result witnessing that our sufficient condition is best possible.

To obtain these new upper bounds on the dichromatic number, we prove a dense decomposition lemma for digraphs having large maximum degree, which generalises to the directed setting the so-called dense decomposition of graphs due to Molloy and Reed. We believe this may be
of independent interest, especially as a tool in various applications.

\end{abstract}

\section{Introduction}

The relationships between the chromatic number $\chi(G)$, the clique number $\omega(G)$, and the maximum degree $\Delta(G)$ of a graph $G$ have attracted a lot of attention during the last decades. Perhaps the first non-trivial result here is Brooks' Theorem~\cite{brooksMPCPS37}, which states that, for a connected graph $G$, $\chi(G)$ is at most $\Delta(G)$, unless $G$ is an odd cycle or a complete graph. Following this, in 1977, Borodin and Kostochka~\cite{borodinJCTB23} made the following well-known conjecture, which, if true, would give a nice extension of Brooks' theorem for graphs with maximum degree at least $9$.

\begin{conjecture}[{\sc Borodin--Kostochka}~\cite{borodinJCTB23}]
    \label{conj:borodin_kostochka}
    Every graph $G$ with maximum degree $\Delta(G) \geq 9$ and clique number $\omega(G)\leq \Delta(G)-1$ satisfies $\displaystyle \chi(G) \leq \Delta(G)-1$.
\end{conjecture}

The graph obtained from a $5$-cycle by blowing up each vertex into a triangle witnesses that the lower bound on the maximum degree is tight.
Twenty years later, Reed published two fundamental papers~\cite{reedJGT27,reedJCTB76} related to Borodin--Kostochka's conjecture. In the first one, Reed posed the following celebrated conjecture.

\begin{conjecture}[{\sc Reed}~\cite{reedJGT27}]
    \label{conj:reed}
    Every graph $G$ satisfies 
    $\chi(G) \leq \left\lceil \frac{1}{2}\big(\Delta(G)+1+\omega(G)\big) \right\rceil$.
\end{conjecture}

If true, it implies in particular that every graph $G$ with chromatic number at least $\Delta(G)$ contains a clique of size $\Delta(G)-2$, thereby approaching Conjecture~\ref{conj:borodin_kostochka}.
As a partial result, Reed~\cite{reedJGT27} shows that $\chi(G)$ is indeed bounded by a convex combination of $\Delta(G)+1$ and $\omega(G)$ (up to ceiling).

\begin{theorem}[\cite{reedJGT27}]
    \label{thm:reed_epsilon}
    There exists $\epsilon >0$ such that every graph $G$ satisfies 
    \[\chi(G) \leq \left\lceil (1-\epsilon)(\Delta(G)+1)+\epsilon\omega(G) \right\rceil.\]
\end{theorem}
 
Theorem~\ref{thm:reed_epsilon} implies the weaker version of Borodin--Kostochka's Conjecture stating that every graph $G$ with chromatic number at least $\Delta(G)$ contains a clique of size $\Delta(G)-C$, for some absolute constant $C$. One year later, Reed~\cite{reedJCTB76} strengthened this by proving that Conjecture~\ref{conj:borodin_kostochka} holds for large values of $\Delta$.

\begin{theorem}[\cite{reedJCTB76}]
    \label{thm:reed_bk}
    There exists $\Delta_{\ref{thm:reed_bk}}\in \mathbb{N}$ such that the following holds.
    Every graph $G$ with maximum degree $\Delta(G) \geq \Delta_{\ref{thm:reed_bk}}$ and clique number $\omega(G) \leq \Delta(G)-1$ satisfies $\chi(G) \leq \Delta(G)-1$.
\end{theorem}

\subsection*{$\bic$, $\deltil$, and $\dic$}

The purpose of this work is to provide extensions of Conjecture~\ref{conj:borodin_kostochka} and Theorem~\ref{thm:reed_bk} to the digraph setting. We thus introduce three digraph parameters extending $\Delta$, $\omega$, and $\chi$. By an extension, we mean that if $D$ is \textbf{symmetric}, that is, $D$ is obtained from a graph $G$ by replacing each edge by two opposite arcs, then these digraph parameters on $D$ coincide with $\Delta$, $\omega$, and $\chi$ on $G$. 
It is straightforward to verify that this holds for $\deltil$, $\bic$, and $\dic$, defined as follows.

Let $D$ be a digraph. We let $\deltil(D)$ denote the maximum geometric mean of the in- and out-degrees of the vertices of $D$, that is $\deltil(D) = \max \{ \sqrt{d^-(v)\cdot d^+(v)} : v \in V(D) \}$.
A \textbf{biclique} of $D$ is a set of vertices inducing a complete digraph, which is a digraph containing all possible arcs ({\it i.e.}, between every pair of vertices there are two oppositely oriented arcs).
The \textbf{biclique number} $\bic(D)$ of $D$ is the size of the largest biclique of $D$.

A \textbf{$k$-dicolouring} of $D$ is a colouring of its vertex-set $\phi\colon V(D) \to [k]$ such that each colour class $\phi^{-1}(i)$ induces an acyclic subdigraph on $D$. The \textbf{dichromatic number} of $D$, denoted by $\dic(D)$, is the smallest $k\in \mathbb{N}$ for which $D$ admits a $k$-dicolouring.
These colouring notions for digraphs were introduced in the late 1970s by Erd\H{o}s and Neumann-Lara~\cite{neumannlaraJCT33,erdosPNCN1979}, rediscovered by Mohar in the 2000s~\cite{moharJGT43,bokalJGT46}, and have received a lot of attention since then. 
It has emerged that many classical results on graph colouring can be extended to digraphs under these definitions, reinforcing the idea that dicolouring constitutes an appropriate generalisation of proper colouring to digraphs.
For instance, the Strong Perfect Graph Theorem~\cite{andresJGT79}, Gy\'arf\'as--Sumner's conjecture~\cite{aboulkerEJC28,aboulkerEJC31b,cookEJC30}, Hadwiger's Conjecture~\cite{axenovichCPC31, meszarosJGT101}, and Kostochka--Yancey's type of bounds~\cite{kostochkaGC36,aboulkerSIAM36,aboulkerEJC31a,havetJGT107,picasarriDM347} have been generalised to digraphs via the notion of dicolouring.

\medskip

Similar to the undirected case, relationships linking $\dic$, $\deltil$, and $\bic$ have recently gained interest. 
This line of research began with a result of Mohar~\cite{moharLAA432}, which states that $\dic(D)\leq \lceil\deltil(D)\rceil$ for every connected digraph $D$ which is not a directed cycle, a symmetric odd cycle, nor a complete digraph, thereby generalising Brooks' Theorem (see also~\cite{harutyunyanSIDMA25, aboulker2023digraph, gonccalves2024brooks} for stronger versions of this result).
Recently, Kawarabayashi and Picasarri-Arrieta~\cite{kawarabayashiSODA25} posed that Reed's conjecture generalises to digraphs with these parameters.

\begin{conjecture}[\cite{kawarabayashiSODA25}]
    \label{conj:kawarabayashi}
    Every digraph $D$ satisfies 
    $\dic(D) \leq \left\lceil \frac{1}{2}\big(\deltil(D) + 1 + \bic(D)\big)\right\rceil$.
\end{conjecture}

If true, this would imply both Reed's conjecture in the case of symmetric digraphs and an independent conjecture posed by Harutyunyan and Mohar~\cite{harutyunyanEJC18} when applied to digraphs with biclique number~$1$. As an intermediate result, they generalised Theorem~\ref{thm:reed_epsilon} to digraphs (up to the value of~$\epsilon$).

\begin{theorem}[\cite{kawarabayashiSODA25}]
    \label{thm:kawarabayashi}
    There exists $\epsilon >0$ such that every digraph $D$ satisfies 
    \[\dic(D) \leq \lceil (1-\epsilon)(\deltil(D) +1) + \epsilon\bic(D)\rceil.\]
\end{theorem}

One might expect that, like Reed’s conjecture, Borodin--Kostochka's conjecture directly extends to digraphs with these parameters; that is, for every integer $\Delta\geq 9$, digraphs $D$ with maximum degree $\deltil(D) \leq \Delta$ and biclique number $\bic(D) \leq \Delta-1$ have dichromatic number at most $\Delta -1$. It turns out that this is not the case: for every integer $k\geq 1$, let $H_k = \vec{C_3} \boxplus \bid{K}_{k-2}$ be the digraph made of a copy $C$ of $\vec{C}_3$ (the directed cycle of length $3$), a disjoint copy $K$ of $\bid{K}_{k-2}$ (the complete digraph of order $k-2$), and all possible arcs in both directions between $K$ and $C$. Then $H_k$ satisfies $\deltil(H_{k}) = k$, $\bic(H_k)=k-1$, and $\dic(D) = k$. We believe that this is the only obstruction arising from the directed setting and conjecture the following.

\begin{conjecture}
    \label{conj:directed_bk}
    Let $\Delta\geq 9$ be an integer and $D$ a digraph satisfying $\deltil(D)\leq \Delta$ and $\bic(D)\leq \Delta-1$. Then $\dic(D) \leq \Delta-1$ unless $D$ contains $\vec{C_3} \boxplus \bid{K}_{\Delta-2}$.
\end{conjecture}

If true, Conjecture~\ref{conj:directed_bk} would imply Conjecture~\ref{conj:borodin_kostochka}, since graphs in Conjecture~\ref{conj:borodin_kostochka} can be replaced by corresponding symmetric digraphs in Conjecture~\ref{conj:directed_bk}. In particular, condition $\Delta\geq 9$ cannot be relaxed, as it is already best possible is in the undirected case.
 Similarly to the undirected case, Theorem~\ref{thm:kawarabayashi} implies the weaker version stating that digraphs $D$ with $\dic(D) \geq \deltil(D)$ have biclique number at least $\deltil(D)-C$ for some absolute constant $C$. Our first contribution is that Conjecture~\ref{conj:directed_bk} holds for digraphs with large maximum degree, thereby generalising Theorem~\ref{thm:reed_bk} to digraphs (up to the value of $\DeltaDirectedBK$).

\begin{restatable}{theorem}{thmdirectedbk}
    \label{thm:directed_bk}
    There exists $\DeltaDirectedBK\in \mathbb{N}$ such that the following holds.
    Let $\Delta\geq \DeltaDirectedBK$ be an integer and $D$ a digraph with $\deltil(D)\leq \Delta$ and $\bic(D)\leq \Delta-1$. Then $\dic(D) \leq \Delta-1$ unless $D$ contains $\vec{C_3} \boxplus \bid{K}_{\Delta-2}$.
\end{restatable}

Our proof of Theorem~\ref{thm:directed_bk} combines structural and probabilistic arguments. In particular, we prove the following dense decomposition lemma, which generalises to the directed setting the so-called dense decomposition of graphs appearing in Molloy and Reed's series of papers~\cite{reedJGT27, molloyCPC7,molloyComb18, reedJCTB76, farzadJCTB93, molloyJCTB109} (see also~\cite[Chapter~15]{molloyreed}).
We believe it is of independent interest, and that it can be used to extend to the directed setting many other results obtained by Molloy and Reed, as well as have applications in the algorithmic theory of digraphs.

Before we state it, we first need a few additional definitions.
Let $D$ be a digraph. We denote by $\delmax(D)$ the maximum, overall vertices $v$ of $D$, of $\max(d^+(v),d^-(v))$. Suppose that $\delmax(D) = \Delta$, and let $0\leq d \leq \Delta-1$ be any real number. We say that a vertex $v$ is {\bf $d$-sparse} if the digraph induced by its out-neighbourhood contains at most $\Delta(\Delta-1) - d\Delta$ arcs. A vertex that is not $d$-sparse is {\bf $d$-dense}. 
Given a subset of vertices $X$, we denote by $\outboundary{X}$ (respectively $\inboundary{X}$) the set of arcs of $D$ going out from $X$ (respectively entering $X$).
In what follows, by $\omega(1) \leq  d \leq o(\Delta)$ we mean that $d$ is any fixed sublinear function of $\Delta$ going to infinity. 

\begin{restatable}[{\sc Dense Decomposition Lemma}]{lemma}{densedecompositionlemma}
    \label{lemma:dense_decomposition}
    For every $0 < \epsilon < \frac{1}{2}$ and $\omega(1) \leq  d \leq o(\delmax)$, there exists $\DeltaDecomp$ such that the following holds.
    Every digraph $D$ with $\delmax(D) = \delmax \geq \DeltaDecomp$ admits a partition $X_1\sqcup \ldots \sqcup X_t \sqcup S$ of its vertex-set such that:
    \begin{enumerate}[label=\rm(\roman*)]
        \item for every $i\in [t]$, $\delmax- \frac{3}{\epsilon} d < |X_i| < \delmax+1 + 4d$;
        \label{enum:dense_decomposition:1}
        \item for every $i\in [t]$, $\outboundarysize{X_i} \leq \frac{8}{\epsilon} \cdot d \cdot \delmax$ and $\inboundarysize{X_i} \leq \frac{8}{\epsilon} \cdot d \cdot \delmax$;
        \label{enum:dense_decomposition:2}
        \item for every $i\in [t]$ and every $u\in V(D)$,
        $u \in X_i$ if and only if $|N^+(u) \cap X_i| \geq (1-\epsilon) \delmax$; and
        \label{enum:dense_decomposition:3}
        \item vertices in $S$ are $d$-sparse.
        \label{enum:dense_decomposition:4}
    \end{enumerate}
\end{restatable}

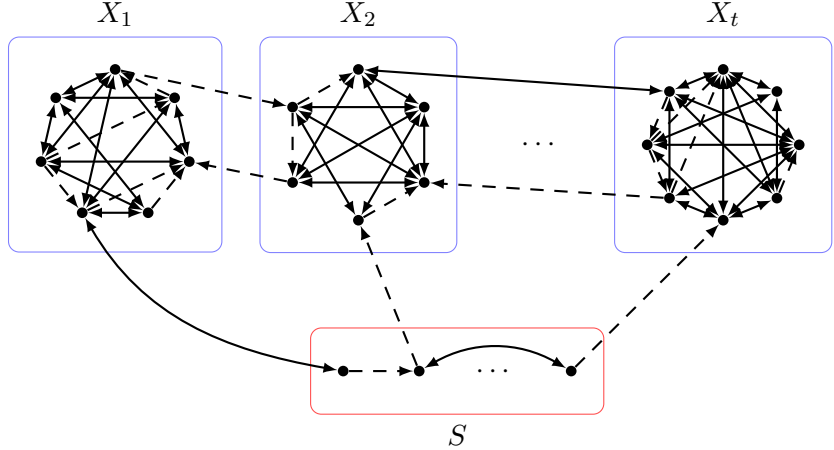
\begin{figure}
\centering
\begin{tikzpicture}[
  Xset/.style={rounded corners, draw=blue!50, inner sep=10pt},
  Sset/.style={rounded corners, draw=red!70, inner sep=10pt},
]

\foreach \i in {1,...,7}{
    \node[vertex] (x1\i) at (360*\i/7+90:1) {};
}
\begin{scope}[on background layer]
  \node[Xset, fit=(x11)(x12)(x13)(x14)(x15)(x16)(x17)(-90:1.07), label=above:$X_1$] {};
\end{scope}

\begin{scope}[xshift=3.2cm]
    \foreach \i in {1,...,6}{
        \node[vertex] (x2\i) at (360*\i/6+90:1) {};
    }
\end{scope}
\begin{scope}[on background layer]
  \node[Xset, fit=(x21)(x22)(x23)(x24)(x25)(x26), label=above:$X_2$] {};
\end{scope}
\node[] (dots) at (5.6,0) {$\cdots$};
\begin{scope}[xshift=8cm]
    \foreach \i in {1,...,8}{
        \node[vertex] (x3\i) at (360*\i/8+90:1) {};
    }
\end{scope}
\begin{scope}[on background layer]
  \node[Xset, fit=(x31)(x32)(x33)(x34)(x35)(x36)(x37)(x38), label=above:$X_t$] {};
\end{scope}
\begin{scope}[xshift=4.5cm]
\node[vertex] (s1) at (-1.5,-3) {};
\node[vertex] (s2) at (-0.5,-3) {};
\node[] (s3) at (0.5,-3) {$\cdots$};
\node[vertex] (s4) at (1.5,-3) {};
\end{scope}

\begin{scope}[on background layer]
  \node[Sset, fit=(s1)(s2)(s3)(s4), label=below:$S$] {};
\end{scope}

\draw[digon] (x11) -- (x12);
\draw[digon] (x11) -- (x14);
\draw[digon] (x11) -- (x16);
\draw[digon] (x11) -- (x17);
\draw[edge] (x12) -- (x13);
\draw[digon] (x12) -- (x14);
\draw[digon] (x12) -- (x15);
\draw[edge] (x12) -- (x16);
\draw[digon] (x12) -- (x17);
\draw[digon] (x13) -- (x14);
\draw[edge] (x13) -- (x15);
\draw[digon] (x13) -- (x16);
\draw[digon] (x13) -- (x17);
\draw[edge] (x14) -- (x15);
\draw[digon] (x15) -- (x16);
\draw[digon] (x15) -- (x17);
\draw[edge] (x16) -- (x17);

\draw[edge] (x21) -- (x22);
\draw[digon] (x21) -- (x23);
\draw[digon] (x21) -- (x24);
\draw[digon] (x21) -- (x25);
\draw[edge] (x21) -- (x26);
\draw[digon] (x22) -- (x24);
\draw[digon] (x22) -- (x25);
\draw[digon] (x22) -- (x26);
\draw[edge] (x23) -- (x24);
\draw[digon] (x23) -- (x25);
\draw[digon] (x24) -- (x25);
\draw[digon] (x24) -- (x26);
\draw[digon] (x25) -- (x26);

\draw[edge] (x31) -- (x32);
\draw[digon] (x31) -- (x33);
\draw[digon] (x31) -- (x35);
\draw[digon] (x31) -- (x36);
\draw[digon] (x31) -- (x38);
\draw[edge] (x32) -- (x33);
\draw[digon] (x32) -- (x34);
\draw[digon] (x32) -- (x36);
\draw[digon] (x32) -- (x37);
\draw[edge] (x32) -- (x38);
\draw[digon] (x33) -- (x34);
\draw[digon] (x33) -- (x36);
\draw[edge] (x33) -- (x38);
\draw[digon] (x34) -- (x35);
\draw[digon] (x34) -- (x36);
\draw[digon] (x34) -- (x38);
\draw[edge] (x35) -- (x36);
\draw[digon] (x35) -- (x37);
\draw[digon] (x35) -- (x38);
\draw[digon] (x36) -- (x38);
\draw[digon] (x37) -- (x38);

\draw[edge] (x17) -- (x21);
\draw[edge] (x22) -- (x15);
\draw[digon] (x13) to[out=-60,in=170] (s1);

\draw[digon] (x31) -- (x26);
\draw[edge] (x33) -- (x24);

\draw[edge] (s1) -- (s2);
\draw[digon] (s2) to[out=30,in=150] (s4);
\draw[edge] (s2) -- (x23);
\draw[edge] (s4) -- (x34);

\end{tikzpicture}
\caption{An illustration of the decomposition $X_1 \sqcup \ldots \sqcup X_t \sqcup S$ provided by Lemma~\ref{lemma:dense_decomposition}. Each set $X_i$ has size close to $\Delta$, and each vertex in $X_i$ has a large fraction of its out-neighbours inside $X_i$. Intuitively, each $X_i$ is thus close to being a biclique on $\Delta$ vertices. Vertices in $S$ are $d$-sparse.}
\label{fig:DDL}
\end{figure}

See Figure~\ref{fig:DDL} for an illustration of a dense decomposition.

We now outline the proof of Theorem~\ref{thm:directed_bk}. Consider a minimal counterexample $D$ along with one of its dense decompositions $(X_1,\ldots,X_t,S)$ generated by Lemma~\ref{lemma:dense_decomposition}, where $\epsilon$ is sufficiently small and $d$ is approximately $\log^3 \Delta$. In a minimal counterexample, it is easy to show that $\deltil$ and $\delmax$ are close, so for clarity, let us assume here that $\Delta = \delmax(D) = \deltil(D)$.

Consider the random colouring process which consists of (i) assigning each vertex of $D$ a colour uniformly at random from $\{1,\dots,\Delta-1\}$, and (ii) simultaneously uncolouring all vertices belonging to a monochromatic directed cycle. The core idea is to show that the resulting partial dicolouring can be extended to $D$ with positive probability.
This method is now standard in the field, and has been used for instance in~\cite{harutyunyanEJC18,kawarabayashiSODA25}.
For a $d$-sparse vertex $v$, applying this process, it is highly likely that many colours are repeated in its out-neighbourhood, thereby allowing us to extend the obtained dicolouring to $v$. This argument is central to the proof of~\cite[Theorem~1.7]{kawarabayashiSODA25}, which establishes that a sparse digraph $H$ can be dicoloured using only a fraction of $\deltil(H)$ colours.

Thus, the main challenge in proving Theorem~\ref{thm:directed_bk} lies in handling $d$-dense vertices. Furthermore, we know that all such vertices belong to some $X_i$, and that each $X_i$ can be coloured almost independently from the rest of the digraph. This is mainly due to property~\ref{enum:dense_decomposition:3} of Lemma~\ref{lemma:dense_decomposition} and the fact that $\epsilon$ is sufficiently small.
Building on this observation, we progressively constrain the possible structure of each $X_i$ until reaching a point where it becomes likely that the random colouring process described above yields a partial dicolouring that can be extended to $X_i$.
At that stage, using the Lovász Local Lemma, we conclude that the obtained partial dicolouring can be extended simultaneously to all sparse vertices and all $X_i$s with positive probability, leading to the desired contradiction.

\subsection*{Analogues for other definitions of \texorpdfstring{$\Delta$}{Delta}}

A natural question is whether $\deltil$ could be replaced by other, potentially more natural, definitions of the maximum degree of a digraph in Conjecture~\ref{conj:directed_bk} and Theorem~\ref{thm:directed_bk}.
Observe that, by definition, all digraphs $D$ satisfy $\delmax(D) \geq \deltil(D)$. Hence, replacing $\deltil$ with $\delmax$ yields weaker statements.
A natural candidate is then the {\bf maximum out-degree} $\delplus$ of a digraph, defined as 
\[ \delplus(D) = \max \{d^+(v) : v\in V(D)\}.\]

Indeed, $\delplus$ and $\deltil$ are independent parameters, in the sense that one cannot be bounded by any function of the other. 
Furthermore, when only the out-degree of a digraph is bounded, probabilistic methods are unlikely to work. Intuitively, if the maximum in-degree of a digraph is arbitrarily large compared to its maximum out-degree, one cannot expect to apply Lov\'asz Local Lemma, as it seems unlikely to bound the degree of ``bad events''. That is why finding upper bounds on the dichromatic number of a digraph $D$ involving $\delplus(D)$ is usually significantly harder than ones involving $\deltil(D)$. The most emblematic open problem of this kind is probably the following one, proposed by Kawarabayashi and Picasarri-Arrieta~\cite{kawarabayashiSODA25}, whose version with $\deltil$ has been proved first by Harutyunyan and Mohar~\cite{harutyunyanEJC18} (see also~\cite{golowichDM339}), and appears to be a particular case of Theorem~\ref{thm:kawarabayashi}.

\begin{problem}[\cite{kawarabayashiSODA25}]
    Show the existence of $\epsilon > 0$ such that every digraph $D$ with $\bic(D)=1$ satisfies $\dic(D) \leq (1-\epsilon)\delplus(D) + O(1)$.
\end{problem}

Nevertheless, we believe that Conjecture~\ref{conj:directed_bk} remains valid when $\deltil$ is replaced by $\delplus$. Furthermore, we are able to prove the analogue of Theorem~\ref{thm:directed_bk} to support it.

\begin{conjecture}
    \label{conj:directed_bk_delplus}
    Every digraph $D$ with maximum out-degree $\delplus(D) =\Delta \geq 9$ and biclique number $\bic(D) \leq \Delta-1$ satisfies $\dic(D) \leq \Delta-1$, unless $D$ contains $\vec{C_3} \boxplus \bid{K}_{\Delta-2}$.
\end{conjecture}

\begin{restatable}{theorem}{thmdirectedbkplus}
    \label{thm:directed_bk_deltaplus}
    There exists $\DeltaDirectedBKPlus\in \mathbb{N}$ such that the following holds. 
    Let $\Delta\geq \DeltaDirectedBKPlus$ be an integer and $D$ a digraph with $\delplus(D)\leq \Delta$ and $\bic(D)\leq \Delta-1$. Then $\dic(D) \leq \Delta-1$ unless $D$ contains $\vec{C_3} \boxplus \bid{K}_{\Delta-2}$.
\end{restatable}

Theorem~\ref{thm:directed_bk_deltaplus} yields a second independent generalisation of Theorem~\ref{thm:reed_bk}.
The proof combines structural and discharging arguments on a well-chosen counterexample, which in particular does not satisfy the requirements of Theorem~\ref{thm:directed_bk}.

\medskip

Finally, leveraging a technique introduced first in~\cite{picasarriJGT106} and formally established in~\cite{kawarabayashiArxiv24},
we derive  from Theorem~\ref{thm:directed_bk_deltaplus} the following two sufficient conditions for a digraph $D$ to satisfy $\dic(D)\leq \delmin(D)-1$, where 
\[\delmin(D) = \max \{ \min(d^-(v), d^+(v)) : v\in V(D)\}.\]

Note that, for digraphs satisfying one of these conditions, this is a strengthening of Theorems~\ref{thm:directed_bk} and~\ref{thm:directed_bk_deltaplus}, as $\delmin$ is always bounded above by both $\deltil$ and $\delplus$. 

\begin{restatable}{corollary}{cordirectedbkmin}
    \label{cor:directed_bk_deltamin}
    For every integer $\Delta\geq \DeltaDirectedBKPlus$ and digraph $D$ with $\delmin(D) \leq \Delta$ and $\bic(D) < \frac{1}{2}(\Delta-1)$, $\dic(D) \leq \Delta-1$.
\end{restatable}

\begin{restatable}{corollary}{cordirectedbkminbis}
    \label{cor:directed_bk_deltamin_bis}
    The following holds for every integer $\Delta\geq \DeltaDirectedBKPlus$. Let $D$ be a digraph with underlying graph $G$. If $\delmin(D) \leq \Delta$ and $\omega(G) \leq \Delta-1$, then $\dic(D) \leq \Delta-1$.
\end{restatable}

Note that Corollary~\ref{cor:directed_bk_deltamin_bis} gives a third independent generalisation of Theorem~\ref{thm:reed_bk}. We are actually able to prove a slightly more general result on digraphs $D$ with $\dic(D) \geq \delmin(D)$, the interested reader is referred to Section~\ref{sec:directed_bk:min} (see Corollary~\ref{thm:directed_bk_deltamin}) for the exact statement.
We finally justify that the condition $\bic(D) < \frac{1}{2}(\Delta-1)$ of Corollary~\ref{cor:directed_bk_deltamin} is tight. In fact, when restricted to digraphs $D$ with $\bic(D) = \lceil \frac{1}{2}(\delmin(D)-1)\rceil$, we prove that determining whether $\dic(D) \leq \delmin (D)-1$ is even NP-complete. 

\begin{restatable}{proposition}{propnpcompleteness}
    \label{prop:npcompleteness}
    For every integer $k \geq 2$, deciding $\dic(D)\leq k$ remains NP-complete when restricted to digraphs $D$ with $\delmin(D) = k+1$ and $\bic(D) = \lceil \frac{k}{2}\rceil$.
\end{restatable}

\paragraph{Organisation of the paper.}
All notation, definitions, and probabilistic tools used in the paper are provided in Section~\ref{sec:notations}.
Section~\ref{sec:dense_decomposition} is dedicated to the proof of Lemma~\ref{lemma:dense_decomposition}.
Section~\ref{sec:directed_bk} is devoted to the proof of Theorem~\ref{thm:directed_bk}. It is mainly divided into two parts: 
in Subsection~\ref{sec:directed_bk:structure}, we use Lemma~\ref{lemma:dense_decomposition} to show a collection of structural properties on a minimum counterexample of Theorem~\ref{thm:directed_bk}. In Subsection~\ref{sec:directed_bk:probabilistic}, we use the probabilistic method to show that a digraph $D$ satisfying these specific properties can actually be coloured with $\Delta - 1$ colours at most, thus showing the result. We prove Theorem~\ref{thm:directed_bk_deltaplus} in Section~\ref{sec:directed_bk:plus}. We finally derive Corollaries~\ref{cor:directed_bk_deltamin} 
 and~\ref{cor:directed_bk_deltamin_bis}, and prove Proposition~\ref{prop:npcompleteness} in Section~\ref{sec:directed_bk:min}.

\section{Preliminaries}
\label{sec:notations}

\subsection{General notation on digraphs}

We only consider digraphs without loops and parallel arcs. Our notation follows~\cite{bang2009}.
A \textbf{digon} is a pair of arcs in opposite directions between the same vertices $u,v$, and is denoted $\dig{u,v}$. A \textbf{simple arc} is an arc that is not in a digon.
A digraph is \textbf{symmetric} if it contains no simple arc.
Given an undirected graph $G$, the digraph $\bid{G}$ is the symmetric digraph obtained from $G$ by replacing every edge with a digon. For an integer $n$, $\bid{K}_n$ is called the \textbf{complete digraph} on $n$ vertices, and a \textbf{biclique} of size $k$ is thus a set of vertices inducing a copy of $\bid{K}_k$. 

Let $D$ be a digraph.
The \textbf{underlying graph} of $D$, denoted by $\UG(D)$, is the undirected graph with vertex-set $V(D)$ in which $uv$ is an edge if and only if $uv$ or $vu$ is an arc of $D$.
We say that $D$ is an {\bf orientation} of $G$ if $D$ does not contain any digon and $G$ is the underlying graph of $D$. A {\bf tournament} is an orientation of a complete graph.

The \textbf{complement} of $D$, denoted by $\overline{D}$, is the digraph with vertex-set $V(D)$ and arc set $(V(D)\times V(D))\setminus A(D)\setminus\{(v,v) : v\in V(D)\}$.
A \textbf{matching} of $D$ is a set of pairwise disjoint arcs of $D$. The maximum size of a matching of $D$ is denoted $\nu(D)$.

Let $xy$ be an arc of $D$. We call $y$ an {\bf out-neighbour} of $x$, and $x$ and {\bf in-neighbour} of $y$. If further $xy$ is simple, we call $y$ a {\bf simple out-neighbour} of $x$, and $x$ a {\bf simple in-neighbour} of $y$. Finally, a {\bf neighbour} of $x$ is either an in- or out-neighbour of $x$.
Let $v$ be any vertex of $D$. The {\bf out-neighbourhood} $N^+(v)$ of $v$ is the set of out-neighbours of $v$, and its {\bf in-neighbourhood} $N^-(v)$ is the set of its in-neighbours. We denote by $N^+[v]$ and $N^-[v]$ the {\bf closed out-neighbourhood} $N^+(v)\cup \{v\}$ and the {\bf closed in-neighbourhood} $N^-(v) \cup \{v\}$  of $v$ respectively.
Furthermore, we denote by $N(v)$ the union of $N^+(v)$ and $N^-(v)$, and by $N^{\pm}(v)$ their intersection.
The {\bf out-} and {\bf in-degrees} of $v$, denoted by $d^+(v)$ and $d^-(v)$, are respectively the number of out- and in-neighbours of $v$.
From this, we can derive the four distinct definitions of the maximum degree of $D$ used throughout the paper:
\begin{flalign*}
    &- \delmax(D) = \max \Big\{ \max(d^-(v),d^+(v)) : v\in V(D)\Big\}, &\\
    &- \deltil(D) = \max \left\{ \sqrt{d^-(v)\cdot d^+(v)} : v\in V(D)\right\}, &\\
    &- \delplus(D) = \max \Big\{ d^+(v) : v\in V(D)\Big\}, \text{~~and} &\\
    &- \delmin(D) = \max \Big\{ \min(d^-(v),d^+(v)) : v\in V(D)\Big\}.
\end{flalign*}

Given two subsets of vertices $X,Y \subseteq V(D)$, we say that $X$ {\bf dominates} $Y$ if $D$ contains every possible arc $uv$ with $u\in X$ and $v\in Y$. 
For an integer $k$, $D$ is called {\bf $k$-critical} if $\dic(D) = k$ and $\dic(H)\leq k-1$ holds for all proper subdigraphs $H$ of $D$.
Given two digraphs $H_1$ and $H_2$, the digraph $H_1 \boxplus H_2$ is the one obtained from disjoint copies of $H_1$ and $H_2$ by adding all possible digons between them.

\subsection{Probabilistic tools}

Our proof of Theorem~\ref{thm:directed_bk} contains some probabilistic arguments. We briefly recall here three classical results used later. 
The reader unfamiliar with the probabilistic method is referred to~\cite{alon2008,molloyreed}.
We first need the symmetric version of the Lov\'asz Local Lemma, due to Erd\H{o}s and Lov\'asz~\cite{erdosIFS10}.

\begin{lemma}[{\sc Lov\'asz Local Lemma}]
\label{lemma:lll}
    Let $A_1,A_2,\dots,A_n$ be events in an arbitrary probability space. Suppose that each event $A_i$ is mutually independent of a set of all the other events but at most $d$, and that $\PP(A_i) \leq p$ for all $1\leq i \leq n$. If 
    \[ ep(d+1) \leq 1,\]
    then $\PP\left(\bigcap_{i=1}^n \overline{A_i}\right)>0$.
\end{lemma}

The second result we need is the celebrated concentration bound due to Talagrand~\cite{talagrand1995}. We use it in the following form, which is similar to that of~\cite{molloyJCTB109}.

\begin{restatable}[{\sc Talagrand's Inequality}]{lemma}{talagrand}
    \label{lemma:talagrand}
    Let $X$ be a random variable valued in $\mathbb{N}$, determined by $n$ independent trials and satisfying the following for some integers $c,r \geq 1$:
    \begin{enumerate}
        \item changing the outcome of any one trial can affect $X$ by at most $c$, and
        \item for every $s\in \mathbb{N}$, if $X\geq s$ then there is a set of at most $rs$ trials whose outcomes certify that $X\geq s$.
    \end{enumerate}
    Then
    \[ 
    \PP\left(|X- \EE(X)| > t \right) \leq 4\exp\left(\frac{-t^2}{32c^2 r(\EE(X) + t)} \right)
    \]
    for any real number $t>126c\sqrt{r\EE(X)} + 344 c^2 r$.
\end{restatable}

The difference between the version above and that of~\cite{molloyJCTB109} is purely technical. In~\cite{molloyJCTB109}, the second condition is asked not only for integers $s$ but for all real numbers $s>0$. To prove that an integer-valued random variable $X$ satisfies such a condition, one needs, in particular, to justify that $X$ is always positive  (as $s$ might be smaller than $\frac{1}{r}$), and to use some rounding arguments. For the sake of clarity in the proofs of our main theorems, we decided to include these technical details directly in the proof of Lemma~\ref{lemma:talagrand}, which we give in Appendix~\ref{sec:appendix}.

The last probabilistic tool we will use is the following consequence~\cite{molloyreed} of the celebrated concentration bound due to Azuma~\cite{azumaTMJ19}.

\begin{lemma}[{\sc Azuma's Inequality}]
    \label{lemma:azuma}
    Let $X$ be a random variable determined by $n$ trials $T_1,\ldots,T_n$, such that, for each $1\leq i\leq n$ and any two possible sequences of outcomes $t_1,\ldots,t_i$ and $t_1,\ldots,t_{i-1},t'_i$ differing only in the last one,
    \[\Big|\,\EE(X\,|\, T_1=t_1,\ldots,T_i=t_i)-\EE(X\,|\, T_1=t_1,\ldots,T_{i-1}=t_{i-1},T_i=t'_i)\,\Big|\leq c_i.\]
    Then, for any $\lambda\geq 0$,
    \[\PP(|X-\EE(X)|>\lambda)\leq 2\exp\left(\frac{-\lambda^2}{2\sum_{i=1}^n c_i^2}\right).\]
\end{lemma}

\subsection{The Directed Brooks' Theorem}

In order to simplify the proof of Theorem~\ref{thm:directed_bk}, we will make use of the directed Brooks' Theorem, due to Mohar~\cite{moharLAA432} (see also~\cite{aboulkerDM113193}).

\begin{theorem}[{\sc Directed Brooks' Theorem}]
    \label{thm:directed_brooks}
    Let $D$ be a digraph with $\dic(D) = \delmax(D) +1$, then one of the following holds for some connected component $H$ of $D$:
    \begin{itemize}
        \item $\delmax(D) = 1$ and $H$ is a directed cycle,
        \item $\delmax(D) = 2$ and $H$ is a symmetric odd cycle, or
        \item $H$ is a complete digraph on $\delmax(D)+1$ vertices.
    \end{itemize}
\end{theorem}

\section{Directed dense decompositions}
\label{sec:dense_decomposition}

This section is devoted to the proof of Lemma~\ref{lemma:dense_decomposition}. Let us first recall it for convenience.

\densedecompositionlemma*
\begin{proof}
    Let us fix $0<\epsilon < \frac{1}{2}$ and $\omega(1) \leq d \leq o(\delmax)$. We do not explicitly give the value of $\DeltaDecomp$; we simply assume that it is large enough so that all inequalities hold in what follows.
    For the sake of conciseness, let us denote $\delmax$ by $\Delta$ along this proof.
    For every $d$-dense vertex $v$, we let $X_v$ be the set of vertices obtained via the following process:
    \begin{enumerate}
        \item at the beginning, $X_v = N^+[v]$;
        \label{enum:dense_decomposition:process:1}
        \item as long as there exists $u\in X_v$ with $|N^+(u) \cap X_v| < (1-\epsilon)\Delta$, remove $u$ from $X_v$;
        \label{enum:dense_decomposition:process:2}
        \item as long as there exists $u\in V(D) - X_v$ with $|N^+(u) \cap X_v| \geq (1-\epsilon)\Delta$, add $u$ to $X_v$.
        \label{enum:dense_decomposition:process:3}
    \end{enumerate}

    \begin{claimlabeled}[lemma:dense_decomposition]
        \label{claim:dense_decomposition:A}
        Every vertex $u\in V(D)$ belongs to $X_v$ if and only if $|N^+(u) \cap X_v| \geq (1-\epsilon)\Delta$.
    \end{claimlabeled}
    \begin{proofclaim}
        This follows directly from the construction of $X_v$.
    \end{proofclaim}

    \begin{claimlabeled}[lemma:dense_decomposition]
        \label{claim:dense_decomposition:C}
        Every $d$-dense vertex $v$ satisfies $d^+(v) \geq \Delta -d$.
    \end{claimlabeled}
    \begin{proofclaim}
        If $d^+(v) < \Delta -d$, then the number of arcs in $D[N^+(v)]$ is less than 
        \[(\Delta-d)(\Delta-d-1) < \Delta(\Delta-1) - d\Delta,\] and $v$ is not a $d$-dense vertex.
    \end{proofclaim}

    \begin{claimlabeled}[lemma:dense_decomposition]
        \label{claim:dense_decomposition:D}
        For every $d$-dense vertex $v$, the number of vertices removed from $X_v$ at step~\ref{enum:dense_decomposition:process:2} is less than $ \frac{2d}{\epsilon}$. In particular, $|N^+(v) \setminus X_v| <  \frac{2d}{\epsilon}$.
    \end{claimlabeled}
    \begin{proofclaim}
        Let $c_\epsilon = \frac{2}{\epsilon}$. Assume for a contradiction that at least $c_\epsilon \cdot d$ vertices are removed at step~\ref{enum:dense_decomposition:process:2} of the process described above. Let $\mathcal{U}$ be the $\lceil c_\epsilon \cdot d\rceil$ first deleted vertices, and let $\mathcal{V}$ be $N^+(v) \setminus \mathcal{U}$. 
        By Claim~\ref{claim:dense_decomposition:C}, $|\mathcal{V}| \geq \Delta-(c_\epsilon +1)d -1$. By construction of $\mathcal{U}$, vertices in $\mathcal{U}$ have less than $(1-\epsilon)\Delta$ out-neighbours in $\mathcal{V}$. Therefore, for each vertex $u \in \mathcal{U}$, we have
        \begin{equation}
            \label{eq:size_V_setminus_Npu}
            |\mathcal{V} \setminus N^+(u)| > \epsilon\Delta - (c_\epsilon +1)d -1.
        \end{equation}
        Observe that, in $D[N^+(v)]$, the number of non-arcs is at least $\sum_{u\in \mathcal{U}} |\mathcal{V} \setminus N^+(u)|$. By~\eqref{eq:size_V_setminus_Npu}, the number of arcs in $D[N^+(v)]$ is thus less than
        \[
            \Delta(\Delta-1) - c_\epsilon d\cdot(\epsilon \Delta - (c_\epsilon +1)d -1) < \Delta(\Delta-1) -d\Delta.
        \]
        This is a contradiction to $v$ being $d$-dense.
    \end{proofclaim}

    \begin{claimlabeled}[lemma:dense_decomposition]
        \label{claim:dense_decomposition:E}
        For every $d$-dense vertex $v$, $\inboundarysize{N^+[v]} \leq (d+2)\Delta$ and $\outboundarysize{N^+[v]} \leq (d+2)\Delta$.
    \end{claimlabeled}
    \begin{proofclaim}
        The vertex $v$ being $d$-dense, there are at least $\Delta^2 - (d+1)\Delta$ arcs in $D[N^+[v]]$. 
        Observe that $\sum_{u\in N^+[v]}d^-(u)$ is exactly $\inboundarysize{N^+[v]} + |A(D[N^+[v]])|$. We thus have 
        \[\Delta(\Delta+1) \geq \sum_{u \in N^+[v]}d^-(u) \geq \inboundarysize{N^+[v]} + \Delta^2 -(d+1)\Delta,\]
        which implies $\inboundarysize{N^+[v]} \leq (d+2)\Delta$. Symmetrically, we have 
        \[\Delta(\Delta+1) \geq \sum_{u \in N^+[v]}d^+(u) \geq \outboundarysize{N^+[v]} + \Delta^2 -(d+1)\Delta,\]
        which implies $\outboundarysize{N^+[v]} \leq (d+2)\Delta$.
    \end{proofclaim}
    
    \begin{claimlabeled}[lemma:dense_decomposition]
        \label{claim:dense_decomposition:F}
        For every $d$-dense vertex $v$, the number of vertices added to $X_v$ at step~\ref{enum:dense_decomposition:process:3} is less than $4d$. In particular, $|X_v \setminus N^+[v]| < 4d$.
    \end{claimlabeled}
    \begin{proofclaim}
        Assume for a contradiction that at least $4d$ vertices are added at step~\ref{enum:dense_decomposition:process:3} of the process described above. Note that the $i^{\text{th}}$ added vertex has at least $(1-\epsilon)\Delta+1-i$ out-neighbours in $N^+[v]$. This implies 
        \[\inboundarysize{N^+[v]} \geq \sum_{i=1}^{\lceil 4d\rceil}\Big( (1-\epsilon)\Delta+1-i \Big) > (d+2)\Delta,\]
        a contradiction to Claim~\ref{claim:dense_decomposition:E}.
    \end{proofclaim}

    \begin{claimlabeled}[lemma:dense_decomposition]
        \label{claim:dense_decomposition:G}
        For every $d$-dense vertex $v$, $\displaystyle \inboundarysize{X_v} < \frac{8}{\epsilon}d\Delta$ and $\displaystyle \outboundarysize{X_v} < \frac{8}{\epsilon}d\Delta$.
    \end{claimlabeled}
    \begin{proofclaim}
        After step~\ref{enum:dense_decomposition:process:1} of the process, $\inboundarysize{X_v} \leq (d+2)\Delta$ and $\outboundarysize{X_v} \leq (d+2)\Delta$ by Claim~\ref{claim:dense_decomposition:E}. During step~\ref{enum:dense_decomposition:process:2}, each time we remove a vertex $u$ from $X_v$, both $\inboundarysize{X_v}$ and  $\outboundarysize{X_v}$ increase by at most $\Delta$. Similarly, every time we add a vertex during step~\ref{enum:dense_decomposition:process:3}, both increase by at most $\Delta$. Since less than $\frac{2d}{\epsilon}$ vertices are removed and less than $4d$ are added by Claims~\ref{claim:dense_decomposition:D} and~\ref{claim:dense_decomposition:F}, the claim follows.
    \end{proofclaim}

    It remains to show the existence of a set of $d$-dense vertices $\mathcal{W}$ such that $(X_v)_{v\in \mathcal{W}}$ are pairwise disjoint, and such that every $d$-dense vertex belongs to $\bigcup_{v\in \mathcal{W}}X_v$. The result then follows from Claims~\ref{claim:dense_decomposition:A},~\ref{claim:dense_decomposition:C},~\ref{claim:dense_decomposition:D},~\ref{claim:dense_decomposition:F}, and~\ref{claim:dense_decomposition:G}. First observe that every $d$-dense vertex $v$ indeed belongs to $X_v$ by Claims~\ref{claim:dense_decomposition:A},~\ref{claim:dense_decomposition:C}, and~\ref{claim:dense_decomposition:D}.
    The following thus implies that such a set $\mathcal{W}$ can be constructed greedily by adding to $\mathcal{W}$ a $d$-dense vertex $u \notin \bigcup_{v\in \mathcal{V}}X_v$.

    \begin{claimlabeled}[lemma:dense_decomposition]
        Let $u$ and $v$ be two $d$-dense vertices. If $X_v \cap X_u \neq \emptyset$, then $v\in X_u$ and $u \in X_v$.
    \end{claimlabeled}
    \begin{proofclaim}
        By symmetry, we only have to show that $u\in X_v$.
        For this, it is sufficient to show that $|N^+(u) \cap X_v| \geq (1-\epsilon)\Delta$, the claim then follows from Claim~\ref{claim:dense_decomposition:A}.
        Let $x$ be any vertex in $X_u \cap X_v$. By Claim~\ref{claim:dense_decomposition:A} applied to both $X_u$ and $X_v$, we have:
        \[ \Delta \geq d^+(x) \geq 2(1-\epsilon)\Delta - |X_u\cap X_v|+1, \]
        which implies $|X_u\cap X_v| \geq (1-2\epsilon)\Delta+1$.
        Together with Claim~\ref{claim:dense_decomposition:F} applied to $u$, we thus have
        \begin{equation}
            \label{eq:size_Npu_cap_Xv}
            |N^+(u)\cap X_v| > (1-2\epsilon)\Delta -4d.
        \end{equation}
        Let us now count the number $p$ of arcs going from $N^+(u)\cap X_v$ to $N^+(u)\setminus X_v$. First, by Claim~\ref{claim:dense_decomposition:G}, observe that 
        \begin{equation}
            \label{eq:upperbound_p}
            p\leq \outboundarysize{X_v} < \frac{8}{\epsilon}d\Delta.
        \end{equation}
        On the other hand, by denseness of $u$, we know that less than $d\Delta$ arcs are missing in the out-neighbourhood of $u$. In particular, we have
        \begin{equation}
            \label{eq:lowerbound_p}
            p \geq |N^+(u)\cap X_v| \cdot |N^+(u)\setminus X_v| - d\Delta.
        \end{equation}
        Combining~\eqref{eq:size_Npu_cap_Xv},~\eqref{eq:upperbound_p}, and~\eqref{eq:lowerbound_p}, straightforward calculations imply
        $ |N^+(u) \setminus X_v| < \frac{9}{\epsilon(1-2\epsilon)}d$. 
        Since $d^+(u) \geq \Delta -d$ by Claim~\ref{claim:dense_decomposition:C}, we finally obtain 
        \[|N^+(u)\cap X_v| > \Delta -\left(1+\frac{9}{\epsilon(1-2\epsilon)}\right)\cdot d > (1-\epsilon)\Delta,\]
        concluding the proof of the Lemma.
    \end{proofclaim}
\end{proof}

\section{Proof of Theorem~\ref{thm:directed_bk}}
\label{sec:directed_bk}

As explained in the introduction, this section is devoted to the proof of Theorem~\ref{thm:directed_bk}. 
Let us first recall it for the reader's convenience.

\thmdirectedbk*

We do not explicitly give the value of $\DeltaDirectedBK$, we simply assume that it is large enough so that all inequalities claimed along the proof actually hold.

We assume the statement does not hold, which means that there exists an integer $\Delta\geq \DeltaDirectedBK$ and a digraph $D$ satisfying $\dic(D) \geq  \Delta \geq \deltil(D)$ that does not contain any copy of $\bid{K}_\Delta$ nor $\vec{C_3}\boxplus \bid{K}_{\Delta-2}$. Among all such digraphs, we choose $D$ for which $|V(D)| + |A(D)|$ is minimum.
We further set $d=\log^3(\Delta)$.

In Subsection~\ref{sec:directed_bk:structure}, with Lemma~\ref{lemma:dense_decomposition} in hand, we show a collection of structural properties on $D$. In Subsection~\ref{sec:directed_bk:probabilistic}, we use the probabilistic method to show that a digraph satisfying these specific properties can be partially coloured with $\Delta - 1$ colours at most, in such a way that the obtained dicolouring can be greedily extended to the whole digraph, thus yielding a contradiction.

\subsection{Structure of a minimum counterexample}
\label{sec:directed_bk:structure}

Remark first that $D$ is $\Delta$-critical, for otherwise we can remove an arc or a vertex without decreasing the dichromatic number, thus yielding a smaller counterexample.
We start with two well-known basic properties of critical digraphs.

\begin{lemma}
    \label{lemma:mindeg}
    Every vertex $v\in V(D)$ satisfies $\Delta-1 \leq d^-(v), d^+(v) \leq \Delta+1$. Moreover, we have $\Delta \leq \delmax(D) \leq \Delta+1$.
\end{lemma}
\begin{proof}
    If some vertex $v\in V(D)$ satisfies $\min (d^-(v), d^+(v)) \leq \Delta-2$, we can extend a $(\Delta-1)$-dicolouring of $D-v$ (which exists as $D$ is $\Delta$-critical) to $D$ by choosing for $v$ a colour that is not appearing in the in- or out-neighbourhood of $v$, thus contradicting $\dic(D)=\Delta$. This shows the lower bound.
    As for the upper bound, assume for a contradiction that $v$ satisfies $\max (d^-(v), d^+(v)) \geq \Delta+2$. Using the lower bound we just proved, we thus have
    \[
    d^+(v) \cdot d^-(v) \geq (\Delta-1)(\Delta+2) > \Delta^2 \geq \deltil(D)^2.
    \]
    This is a contradiction to the definition of~$\deltil(D)$.
    This shows the first part of the statement.

    For the second part, it remains to justify that $\delmax(D) \neq \Delta-1$. Suppose that $\delmax(D) = \Delta-1$, then by the directed Brooks' Theorem (Theorem~\ref{thm:directed_brooks}), and because $\dic(D)\geq \Delta$, some connected component of $D$ must be isomorphic to $\bid{K}_\Delta$, a contradiction. 
\end{proof}

\begin{lemma}
    \label{lemma:property_simplearc}
    If $uv$ is a simple arc of $D$, then $v$ has a simple out-neighbour, and $u$ has a simple in-neighbour.
\end{lemma}
\begin{proof}
    By criticality of $D$, let $\phi$ be a $(\Delta-1)$-dicolouring of $D\setminus \{uv\}$. Since $\phi$ is not a dicolouring of $D$, $D$ coloured with $\phi$ contains a monochromatic directed cycle, and among all of them we let $\cSCR$ be an induced one. Clearly, $\cSCR$ contains $uv$, for otherwise $\cSCR$ is a monochromatic directed cycle of $D\setminus \{uv\}$ coloured with $\phi$.
    Since $\cSCR$ is induced and $uv$ is simple, the successor of $v$ in $\cSCR$ is indeed one of its simple out-neighbours. Similarly, the predecessor of $u$ in $\cSCR$ is one of its simple in-neighbours.
\end{proof}

By definition, $D$ does not contain any copy of $\bid{K}_{\Delta}$ nor of $\vec C_3\boxplus\bid{K}_{\Delta-2}$. The following shows that it remains true even when we add one or two arcs between two vertices of $D$.

\begin{lemma}
    \label{lemma:no_KD_minus_uv}
    For every pair of vertices $u,v\in V(D)$, $D\cup \{uv,vu\}$ does not contain any copy of $\bid{K}_{\Delta}$ nor of $\vec C_3\boxplus\bid{K}_{\Delta-2}$. 
\end{lemma}
\begin{proof}
    Throughout the proof, we will often rely on the fact that $\min(d^-(w),d^+(w))\leq\Delta$ for every vertex $w$ of $D$. Let us see first that $D\cup\{uv,vu\}$ does not contain any induced copy of $\vec C_3\boxplus\bid{K}_{\Delta-2}$. (If $D\cup\{uv,vu\}$ contains a non-induced copy of $\vec C_3\boxplus\bid{K}_{\Delta-2}$, then it contains a copy of $\bid{K}_{\Delta}$ and we will deal with this later.) Assume for a contradiction that this is not true. Let $T$ be a set of vertices inducing a directed $3$-cycle in $D\cup\{uv,vu\}$ and $W$ a biclique of $D\cup\{uv,vu\}$ of size $\Delta-2$, such that all arcs between $T$ and $W$ are in $A(D)\cup\{uv,vu\}$. Note that $u,v\in T\cup W$; otherwise, $D$ itself would contain a copy of $\vec C_3\boxplus\bid{K}_{\Delta-2}$, a contradiction.
    Let $\varphi$ be a $(\Delta-1)$-dicolouring of $D-T-W$, which exists due to the $\Delta$-criticality. We are going to extend $\varphi$ to a $(\Delta-1)$-dicolouring of $D$, yielding a contradiction. 
    
    Since $T$ induces a directed $3$-cycle in $D\cup\{uv,vu\}$, at least one of $\{u,v\}$, say $v$, belongs to $W$. Let $x$ and $y$ be two distinct vertices in $T\setminus\{u\}$. Since none of $\{\dig{x,y},\dig{u,v}\}$ is a digon of $D$, and each of these four vertices has at most three out-neighbours outside $T\cup W$, we can extend $\varphi$ to $V(D)\setminus W\cup\{u,v\}$ in a way that $\varphi(x)=\varphi(y)\neq\varphi(u)=\varphi(v)$.
    Then, each $w\in W\setminus\{u,v\}$ has at most $\Delta-2$ forbidden colours, so we can extend $\varphi$ to $V(D)$ greedily.

    Let us now see that $D\cup\{uv,vu\}$ does not contain any copy of $\bid{K}_{\Delta}$. Suppose for a contradiction that we can find a set $W$ of $\Delta-2$ vertices such that $X = \{u,v\}\cup W$ is a biclique of $D\cup\{uv,vu\}$ of size $\Delta$.
    
    \begin{claim}\label{claim:no_KD_minus_uv:A} 
        For every $w\in W$, $d^-(w)=d^+(w)=\Delta$.
    \end{claim}
    \begin{proofclaim} 
    Since $\deltil(D) \leq \Delta$, it is enough to show that $d^-(w),d^+(w)\geq\Delta$. Since $D$ is $\Delta$-critical, there is a $(\Delta-1)$-dicolouring $\varphi$ of $D-X$. Note that $\dig{u,v}$ is not a digon, otherwise $X$ would be a biclique of $D$, a contradiction. Moreover, by Lemma~\ref{lemma:mindeg}, both $u$ and $v$ have at most three out-neighbours outside $W$. Therefore, $\varphi$ can be extended to $V(D)\setminus W$ by assigning the same colour to $u$ and $v$. Now, if there was a vertex $w\in W$ with $\min(d^-(w),d^+(w))\leq\Delta-1$, then, due to $u$ and $v$ having the same colour, we could extend $\varphi$ to a $(\Delta-1)$-dicolouring of $D$ greedily, by taking care that $w$ is coloured in the last place, a contradiction.
    \end{proofclaim}
    
    Let us denote by $w_1,\ldots,w_{\Delta-2}$ the vertices in $W$, and, for $1\leq i\leq\Delta-2$, let $w_i^-$ and $w_i^+$ be the unique in-neighbour and out-neighbour of $w_i$ in $V(D)\setminus X$, given by Claim~\ref{claim:no_KD_minus_uv:A}. Note that  $w_i^-$ and $w_i^+$ are not necessarily distinct, and that some (non-trivial) identities of the kind $w_i^-=w_j^-$, $w_i^+=w_j^+$, or $w_i^-=w_j^+$ may hold. Generalising the argument used in the proof of Claim~\ref{claim:no_KD_minus_uv:A}, we actually obtain the following.
    
    \begin{claim}
        \label{claim:no_KD_minus_uv:B} 
        Let $\varphi$ be a $(\Delta-1)$-dicolouring of $D-X$. Then,
    \begin{enumerate}[label=$(\roman*)$]
        \item\label{enum:no_KD_minus_uv:monochromatic_path} for every $1\leq i\leq\Delta-2$ there is a monochromatic directed path from $w_i^+$ to $w^-_i$ in $D-X$ (possibly of length $0$), and
        \item\label{enum:no_KD_minus_uv:monochromatic_neighbours} $\varphi(w^-_1)=\varphi(w^+_1)=\ldots=\varphi(w^-_{\Delta-2})=\varphi(w^+_{\Delta-2})$.
    \end{enumerate}
    \end{claim}
    \begin{proofclaim}
        Assume that for some $1\leq i\leq\Delta-2$, property~\ref{enum:no_KD_minus_uv:monochromatic_path} does not hold. Then we can extend $\varphi$ to a $(\Delta-1)$-dicolouring of $D$ as follows. First, we assign the same colour to $u$ and $v$ (we have already justified in the proof of Claim~\ref{claim:no_KD_minus_uv:A} that this is possible), and then, we colour the vertices in $W\setminus\{w_i\}$ greedily. Finally, we assign to $w_i$ a colour not assigned to any vertex in $X\setminus\{w_i\}$. This last step does not create any monochromatic directed cycle in $D[X]$, nor, by our assumption, any involving $w_i$ and the vertices outside $X$. So $\varphi$ can be extended to a $(\Delta-1)$-dicolouring of $D$, a contradiction.

        Now let us  see~\ref{enum:no_KD_minus_uv:monochromatic_neighbours}. By~\ref{enum:no_KD_minus_uv:monochromatic_path} we already know that $\varphi(w_i^-)=\varphi(w_i^+)$ for every $1\leq i\leq\Delta-2$. Assume that $\varphi(w_i^+)\neq\varphi(w_j^+)$ for some $i\neq j$. Then we can extend $\varphi$ to a $(\Delta-1)$-dicolouring of $D$ as follows. First, we assign the same colour to $u$ and $v$, with the additional requirement that this colour is different from $\varphi(w_j^+)$. The justification that this can be done without creating monochromatic directed cycles is the same as before. Then, we assign to $w_i$ the colour $\varphi(w_j^+)$. By assumption, this does not create any monochromatic directed cycle. We can then colour the vertices in $W\setminus\{w_j\}$ greedily, because they have two out-neighbours sharing a colour and one uncoloured out-neighbour. Finally, $w_j$ has two pairs of out-neighbours sharing a colour, namely, $u,v$ and $w_i,w_j^+$, so $w_j$ can also be coloured greedily using one of the original $\Delta-1$ colours, a contradiction.
    \end{proofclaim}
    
    The point is that now some structure is forced.

    \begin{claim}
    \label{claim:no_KD_minus_uv:C} 
    Let $1\leq i,j\leq\Delta-2$.
    \begin{enumerate}[label=$(\roman*)$]
        \item\label{enum:no_KD_minus_uv:obstruction_1} If $w^-_i\neq w_j^+$ and $w^-_j\neq w^+_i$, then the digraph $D'= D-X\cup\{w^-_i w^+_j,w^-_j w^+_i\}$ contains a copy of $\bid{K}_{\Delta}$ or of $\vec C_3\boxplus\bid{K}_{\Delta-2}$.
        \item\label{enum:no_KD_minus_uv:obstruction_2} If $w^-_i\neq w_j^+$ and $w^-_j=w^+_i$, then the digraph $D-X\cup\{w^-_i w^+_j\}$ contains a copy of $\bid{K}_{\Delta}$ or of $\vec C_3\boxplus\bid{K}_{\Delta-2}$.
    \end{enumerate}
    \end{claim}
    \begin{proofclaim}
        We prove~\ref{enum:no_KD_minus_uv:obstruction_1};~\ref{enum:no_KD_minus_uv:obstruction_2} can be proved analogously. Let us begin by showing that $\vec\chi(D')\geq\Delta$. Assume for a contradiction that there is a $(\Delta-1)$-dicolouring $\varphi$ of $D'$. Since $\varphi$ is also a dicolouring of $D-X$, by Claim~\ref{claim:no_KD_minus_uv:B} we can find directed paths from $w_i^+$ to $w_i^-$ and from $w_j^+$ to $w_j^-$, in which all the vertices have been assigned the colour $\varphi(w_i^+)$. Together with the arcs $w_i^- w_j^+$ and $w_j^- w_i^+$, this yields a monochromatic directed cycle in $D'$. Hence, $\vec\chi(D')\geq\Delta$.

        By construction, observe that both inequalities $d^-_{D'}(x) \leq d^-_{D}(x)$ and $d^+_{D'}(x) \leq d^+_{D}(x)$ hold for all vertices $x\in V(D')$. Therefore, $D'$ satisfies  $\dic(D')\geq \Delta \geq \deltil(D')$. The choice of $D$ implies that $D'$ contains a copy of $\bid{K}_{\Delta}$ or of $\vec C_3\boxplus\bid{K}_{\Delta-2}$.
    \end{proofclaim}

    As a consequence, we obtain the following.

    \begin{claim}\label{claim:no_KD_minus_uv:Cbis}
    Let $1\leq i,j\leq\Delta-2$ such that $w_i^-\neq w_j^+$. Then, either $N^{\pm}(w_i^-)\cap N^{\pm}(w_j^+)$ contains a biclique of $D-X$ of size $\Delta-4$, or $w^-_j\neq w_i^+$ and $N^{\pm}(w_j^-)\cap N^{\pm}(w_i^+)$ contains such a biclique.
    \end{claim}
    \begin{proofclaim}
        Let $K$ be the copy of $\bid K_{\Delta}$ or of $\vec C_3\boxplus\bid{K}_{\Delta-2}$ found by~\ref{enum:no_KD_minus_uv:obstruction_1} or~\ref{enum:no_KD_minus_uv:obstruction_2} of Claim~\ref{claim:no_KD_minus_uv:C}. Recall that $K$ cannot be a subdigraph of $D$, so it uses at least one arc of $\{w_i^- w_j^+,w_j^- w_i^+\}$.
        Assume that $w_i^- w_j^+ \in A(K)$, the other case being symmetric.

        Note that, in $K$, for each pair of vertices $\{x,y\}$, there exists a biclique of size $\Delta-3$ in $N^\pm_{D'}(x)\cap N^\pm_{D'}(y)$. In particular $K$ contains a biclique $Y_{ij}$ of size $\Delta-3$ in $N^\pm_{D'}(w_i^-) \cap N^\pm_{D'}(w_j^+)$. Note that $Y_{ij}$ is not necessarily a biclique of $D$ as, in $K$, it may contain the arc $w_j^- w_i^+$. 
        However, either $Y_{ij}\setminus \{w_j^-\}$ or $Y_{ij}\setminus \{w_i^+\}$ is a biclique of $D$ in $N^\pm_{D}(w_i^-) \cap N^\pm_{D}(w_j^+)$. The result follows.
    \end{proofclaim}

    With the next claim, we make sure that Claim~\ref{claim:no_KD_minus_uv:Cbis} can be used.
    
    \begin{claim}\label{claim:no_KD_minus_uv:D} There exist $1\leq i,j\leq\Delta-2$ such that $w^-_{i}\neq w^+_{j}$.
    \end{claim}
    \begin{proofclaim}
        Suppose that the statement is not true. Then, there is a vertex $w$ such that $w^-_i=w^+_i=w$ for every $1\leq i\leq\Delta-2$. In other words, $w$ forms a digon with every vertex of $W$, just as $u$ and $v$. Observe that the subdigraph of $D$ induced by $\{u,v,w\}$ is acyclic. Indeed, if $D[\{u,v,w\}]$ had a digon, the vertices of the digon together with $W$ would form a biclique of $D$ of size $\Delta$. And if, otherwise, $D[\{u,v,w\}]$ was a copy of $\vec C_3$, then $D[X\cup \{w\}]$ would be isomorphic to $\vec C_3\boxplus\bid{K}_{\Delta-2}$, a contradiction.
        
        Now, by the $\Delta$-criticality of $D$, there is a $(\Delta-1)$-dicolouring of $D-X-w$. Since the set $\{u,v,w\}$ is acyclic (and because, by Lemma~\ref{lemma:mindeg}, $u$, $v$ and $w$ have at most $3$ out-neighbours outside $W$ each), $\varphi$ can be extended to a $(\Delta-1)$-dicolouring of $D-W$ by assigning the same colour to $u$, $v$, and $w$. Then, by Claim~\ref{claim:no_KD_minus_uv:A}, this can be extended to a $(\Delta-1)$-dicolouring of $D$ by colouring the vertices in $W$ greedily, a contradiction. 
    \end{proofclaim}

    The final claim shows the existence of vertices with in-degree and out-degree at least $\Delta+2$. This contradicts Lemma~\ref{lemma:mindeg}, thus completing the proof of Lemma~\ref{lemma:no_KD_minus_uv}.
    
    \begin{claim}\label{claim:no_KD_minus_uv:E} There exists some $1\leq i\leq\Delta-2$ and two sets of vertices $Z^-,Z^+\subseteq V(D)$ such that
    \begin{enumerate}[$(\roman*)$]
        \item\label{enum:no_KD_minus_uv:induction1} $|Z^-|+|Z^+|\geq 11$;
        \item\label{enum:no_KD_minus_uv:induction2} $w_i^-\in Z^-$ and, unless $Z^+$ is empty, $w_i^+\in Z^+$;
        \item\label{enum:no_KD_minus_uv:induction3} every vertex in $Z^-\cup Z^+$ is of the form $w_j^-$ or $w_j^+$ for some $1\leq j\leq\Delta-2$.
        \item\label{enum:no_KD_minus_uv:induction4} for every $z\in Z^-\setminus\{w_i^-\}$, $N^{\pm}(w_i^-)\cap N^{\pm}(z)$ contains a biclique $Y^-_z$ of $D-X$ of size $\Delta-4$;
        \item\label{enum:no_KD_minus_uv:induction5} for every $z\in Z^+\setminus\{w_i^+\}$, $N^{\pm}(w_i^+)\cap N^{\pm}(z)$ contains a biclique $Y^+_z$ of $D-X$ of size $\Delta-4$.
    \end{enumerate}
    \end{claim}
    \begin{proofclaim}
        The proof is by induction on $|Z^-|+|Z^+|$. Let us prove first that the statement holds if~\ref{enum:no_KD_minus_uv:induction1} is replaced by $|Z^-|\geq 2$. By Claims~\ref{claim:no_KD_minus_uv:Cbis} and~\ref{claim:no_KD_minus_uv:D}, it is enough to take $Z^-=\{w_i^-,w_j^+\}$ and $Z^+=\emptyset$ for some adequate $1\leq i,j\leq\Delta-2$. 

        Let us now assume that the statement holds if~\ref{enum:no_KD_minus_uv:induction1} is replaced by $|Z^-|+|Z^+|\geq t$, with $2\leq t\leq 10$, and prove that it is also true with $|Z^-|+|Z^+|\geq t+1$. By Lemma~\ref{lemma:mindeg}, each vertex $z\in Z^-\cup Z^+$ has at most five in-neighbours and five out-neighbours in $W$. Since $\Delta$ is large, there exists some $1\leq k\leq\Delta-2$ with $w_k^-,w_k^+\notin Z^-\cup Z^+$. We apply Claim~\ref{claim:no_KD_minus_uv:Cbis} to $w_i^-$ and $w_{k}^+$. There are two possible outcomes.
        If $N^{\pm}(w_{i}^-)\cap N^{\pm}(w_k^+)$ contains a biclique of $D-X$ of size $\Delta-4$, then we simply add $w_k^+$ to $Z^-$. Else, $N^{\pm}(w_{k}^-)\cap N^{\pm}(w_i^+)$ contains a biclique of $D-X$ of size $\Delta-4$, and we add $w_k^-$ to $Z^+$, together with $w^+_i$ if it was not already in $Z^+$.
    \end{proofclaim}

    Let $1\leq i\leq\Delta-2$ and $Z^-,Z^+\subseteq V(D)$ given by Claim~\ref{claim:no_KD_minus_uv:E}. Let us assume that $|Z^-|\geq 6$; if $|Z^+|\geq 6$, the proof goes the same way. We can further assume that $Z^-$ has exactly $6$ elements; otherwise, we simply truncate it.  Observe that, for any $z\in Z^-\setminus\{w_i^-\}$, the biclique $Y_z^-$ given by Claim~\ref{claim:no_KD_minus_uv:E} is disjoint from $Z^-$: indeed, by Claims~\ref{claim:no_KD_minus_uv:B} and~\ref{claim:no_KD_minus_uv:E}, any $(\Delta-1)$-dicolouring of $D-X$ is constant on $Z^-$, and $z$ forms a digon with every vertex of $Y^-_z$. Now, by Lemma~\ref{lemma:mindeg}, 
    for any $z\in Z^-\setminus\{w_i^-\}$, $|N^\pm(w_i^-) \setminus Y_{z}^-| \leq 4$. Therefore, 
    \[
        \left|\bigcap_{z\in Z^-\setminus\{w_i^-\}}Y_z^-\right| \geq (\Delta-4) - 4(|Z^-|-2) >0.
    \]
    Let thus $y\in\bigcap_{z\in Z^-\setminus\{w_i^-\}}Y_z^-$ and $z\in Z^-\setminus\{w_i^-\}$. Note that $y$ forms a digon with every vertex in $Y_z^-\sqcup Z^-\setminus\{y\}$. Thus $|N^\pm(y)|\geq\Delta+1$, a contradiction to $\deltil(D) \leq \Delta$.
    \end{proof}

    With the same type of arguments, we now forbid copies of $\vec{C_3}\boxplus \bid{K}_{\Delta-3}$ in $D$.

    \begin{lemma}
        \label{lemma:no_C3_K_deltaminus3}
        The digraph $D$ does not contain any copy of $\vec{C_3}\boxplus \bid{K}_{\Delta-3}$.
    \end{lemma}
    \begin{proof}
        Assume for a contradiction that there exist two disjoint sets $\{x,y,z\},W \subseteq V(D)$ such that $D[\{x,y,z\}]$ contains a $3$-cycle, $W$ is a biclique of size $\Delta-3$, and $D$ contains all possible digons between $\{x,y,z\}$ and $W$. We let $X=\{x,y,z\}\cup W$. Note that, by Lemma~\ref{lemma:no_KD_minus_uv}, $D[\{x,y,z\}]$ contains at most one digon. We thus assume without loss of generality that both $xy$ and $yz$ are simple arcs.
        
        \begin{claim}
            For every $w\in W$, $d^-(w)=d^+(w)=\Delta$.
        \end{claim}
        \begin{proofclaim}
            The proof is similar to that of Claim~\ref{claim:no_KD_minus_uv:A}.
            By Lemma~\ref{lemma:mindeg}, it is enough to show that $d^-(w),d^+(w)\geq\Delta$. Let $\phi$ be a $(\Delta-1)$-dicolouring of $D-X$. 
            We extend $\phi$ to $V(D) \setminus W$ by assigning the same colour to $x$ and $y$. 
            Now, if there was a vertex $w\in W$ with $\min(d^-(w),d^+(w))\leq\Delta-1$, then, due to $x$ and $y$ having the same colour, we could extend $\varphi$ to a $(\Delta-1)$-dicolouring of $D$ greedily, by taking care that $w$ is coloured in the last place, a contradiction.
        \end{proofclaim}
        
        We label $w_1,\dots,w_{\Delta-3}$ the vertices of $W$, and for every $i\in [\Delta-3]$, we denote by $w_i^-$ and $w_i^+$ the unique in- and out-neighbour of $w_i$ outside $X$. 
        Following exactly the proof of Claim~\ref{claim:no_KD_minus_uv:B}, we obtain the following. We give it for completeness.
        
        \begin{claim}
            \label{claim:colouring_D-X}
            Let $\phi$ be a $(\Delta-1)$-dicolouring of $D-X$. Then,
            \begin{enumerate}[label=$(\roman*)$]
                \item for every $1\leq i\leq\Delta-3$ there is a monochromatic directed path from $w_i^+$ to $w^-_i$ in $D-X$ (possibly of length $0$), and
                \label{enum:no_KD_minus_uv:monochromatic_path:2}
                \item $\varphi(w^-_1)=\varphi(w^+_1)=\ldots=\varphi(w^-_{\Delta-3})=\varphi(w^+_{\Delta-3})$.
                \label{enum:no_KD_minus_uv:monochromatic_neighbours:2}
            \end{enumerate}
        \end{claim}
        \begin{proofclaim}
            If~\ref{enum:no_KD_minus_uv:monochromatic_path:2} does not hold for some $i\in [\Delta-3]$, then we can extend $\varphi$ to a $(\Delta-1)$-dicolouring of $D$ as follows. First, we assign a same colour to $x$ and $y$, and then, we colour the vertices in $\{z\}\cup W\setminus\{w_i\}$ greedily, starting with $z$. Finally, we assign to $w_i$ a colour not assigned to any vertex in $X\setminus\{w_i\}$. 

            For~\ref{enum:no_KD_minus_uv:monochromatic_neighbours:2}, assume that $\varphi(w_i^+)\neq\varphi(w_j^+)$ for some $i\neq j$. Then we can extend $\varphi$ to a $(\Delta-1)$-dicolouring of $D$ as follows. First, we assign the same colour $c_1$ to $x$ and $y$, and a colour $c_2$ to $z$, with the additional requirement that both $c_1$ and $c_2$ are different from $\varphi(w_j^+)$. Then, we assign to $w_i$ the colour $\varphi(w_j^+)$. We can then colour the vertices in $W\setminus\{w_j\}$ greedily. At the end, $w_j$ has two pairs of out-neighbours sharing a colour, namely, $x,y$ and $w_i,w_j^+$, so $w_j$ can also be coloured greedily, a contradiction.
        \end{proofclaim}
        
        Actually, all vertices $w_i^-,w_i^+$ must be the same one vertex.
        
        \begin{claim}
            $w_1^- = w_1^+ = \ldots = w_{\Delta-3}^-=w_{\Delta-3}^+$
        \end{claim}
        \begin{proofclaim}
            Assume first that, for some $i\in [\Delta-3]$, $w_i^- \neq w_i^+$. Let $D'$ be the digraph obtained from $D-X$ by further adding the arc $w_i^-w_i^+$. Note that $\deltil(D) \geq \deltil(D')$, and that $D'$ does not contain any copy of $\bid{K}_{\Delta}$ nor of $\vec{C_3}\boxplus \bid{K}_{\Delta-2}$ by Lemma~\ref{lemma:no_KD_minus_uv}. Therefore, the minimality of $|V(D)|+|A(D)|$ implies the existence of a $(\Delta-1)$-dicolouring of $D'$, which by construction is a dicolouring of $D-X$ without any monochromatic directed path from $w_i^+$ to $w_i^-$, a contradiction to Claim~\ref{claim:colouring_D-X}.

            Assume then that, for some $i,j$, $w_i^- = w_i^+ \neq w_j^- = w_j^+$. Similarly, let $D'$ be the digraph obtained from $D-X$ by further adding the digon $\dig{w_i^-,w_j^-}$. Since $w_i^- = w_i^+$ and $ w_j^- = w_j^+$, we have $\deltil(D) \geq \deltil(D')$, and $D'$ does not contain any copy of $\bid{K}_{\Delta}$ nor of $\vec{C_3}\boxplus \bid{K}_{\Delta-2}$ by Lemma~\ref{lemma:no_KD_minus_uv}. Therefore, the minimality of $|V(D)|+|A(D)|$ implies the existence of a $(\Delta-1)$-dicolouring $\phi$ of $D'$, which by construction is a dicolouring of $D-X$ satisfying $\phi(w_i^-) \neq \phi(w_j^-)$, a contradiction to Claim~\ref{claim:colouring_D-X}.
        \end{proofclaim}
        
        Let $w\in V(D) \setminus X$ be such that $W\subseteq N^{\pm}(w)$. 
        Note that $D$ contains all possible digons between $W$ and $\{x,y,z,w\}$. By Lemma~\ref{lemma:no_KD_minus_uv}, we get that, in $D[\{x,y,z,w\}]$, no vertex is incident to two distinct digons. It is straightforward to check that, in particular, $\overline{D}[\{x,y,z,w\}]$ admits a matching $\{x_1y_1,x_2y_2\}$.

        Let $\phi$ be a $(\Delta-1)$-dicolouring of $D-(X\cup \{w\})$. We can extend $\phi$ to $D$ as follows. We first give a same colour $c_1$ to both $x_1$ and $y_1$ that is not appearing in $N^+(x_1)\cup N^+(y_1)$, and similarly a same colour $c_2$ to both $x_2$ and $y_2$.
        Now, two colours are repeated in the out- and in-neighbourhoods of every vertex of $W$, so $\phi$ can be greedily extended to $D$, yielding the contradiction.
    \end{proof}

For the next lemma, we first need a few specific definitions.
We say that $X\subseteq V(D)$ is a {\bf quasi-biclique} of $D$ if one of the following holds:
\begin{enumerate}[label=(\arabic*)]
    \item $X$ is a biclique;
    \label{enum:def_quasibiclique:1}
    \item $D[X] \cup \{xy\} = \bid{K}_{|X|}$ for some $x,y\in X$;
    \label{enum:def_quasibiclique:2}
    \item $D[X]=\vec{C_3} \boxplus \bid{K}_{|X|-3}$;
    \label{enum:def_quasibiclique:3}
    \item $D[X] -x= \bid{K}_{|X|-1}$ and $|N^+(x) \cap X|\geq \frac{99}{100}\Delta$ hold for some $x\in X$; or
    \label{enum:def_quasibiclique:4}
    \item $X$ admits a partition $(R,W,K)$ such that $|R|\leq 4$, $|K| > \frac{1}{2}|X|$, $W\cup K$ is a biclique, $R$ dominates $K$, and vertices in $K$ have in-degree exactly $\Delta+1$.
    \label{enum:def_quasibiclique:5}
\end{enumerate}
For $\rm i\in [5]$, we say that $X$ is a {\bf quasi-biclique of type~(i)} if $X$ satisfies property~(i) above and none of the properties~(j) for $\rm j<i$. 

In what follows, by applying the dense-decomposition lemma to $D$, we show that $D$ is made of a collection of pairwise-disjoint quasi-bicliques and a set of sparse vertices. Recall that $d$ was set to $\log^3\Delta$ at the beginning of the section.

\begin{lemma}
    \label{lemma:quasiclique_decomposition}
    There is a partition $X_1\sqcup \ldots \sqcup X_t \sqcup S$ of $V(D)$ such that:
    \begin{enumerate}[label=$(\roman*)$]
        \item each $X_i$ is a quasi-biclique of $D$ of size $|X_i| > \Delta - 300d$; 
        \item for every $i\in [t]$ and $u\in V(D)\setminus X_i$, $|N^+(u) \cap X_i| < \frac{100}{101}\Delta$.
        \item vertices in $S$ are $d$-sparse.
    \end{enumerate}
\end{lemma}
\begin{proof}
    From now on, we denote $\delmax(D)$ by $\delmax$ for conciseness.
    Recall that, by Lemma~\ref{lemma:mindeg}, $\Delta \leq \delmax \leq \Delta+1$.
    We apply Lemma~\ref{lemma:dense_decomposition} with $\epsilon = \frac{1}{100}$ and $d=\log^3(\Delta)$. 
    We thus obtain a partition $X_1\sqcup \ldots \sqcup X_t \sqcup S$ of $V(D)$ such that:
    \begin{itemize}
        \item for every $i\in [t]$, $|X_i| > \Delta - 300d$;
        \item for every $i\in [t]$ and every $u\in V(D)$, $u\in X_i$ if and only if $|N^+(u) \cap X_i| \geq \frac{99}{100}\delmax$; and
        \item vertices in $S$ are $d$-sparse.
    \end{itemize}

    In particular, the second property implies that vertices $u\in V(D)\setminus X_i$ satisfy $|N^+(u) \cap X_i| < \frac{99}{100}(\Delta+1) < \frac{100}{101}\Delta$.
    Clearly, it thus remains to prove that each $X_i$ is a quasi-biclique. We thus fix $i\in [t]$, and denote $D[X_i]$ by $D_i$. To prove that $X_i$ is a quasi-biclique, our strategy is first to show that the complement $\overline{D_i}$ of $D_i$ has no matching of size $3$. We first prove that $\nu\left(\overline{D_i}\right)$ is indeed at most a small fraction of $\Delta$.

    \begin{claim}
        \label{claim:matching_Di_linear}
        $\displaystyle \nu\left( \overline{D_i} \right) <  \frac{1}{50} \delmax$.
    \end{claim}
    \begin{proofclaim}
        Assume for a contradiction that $\nu\left( \overline{D_i} \right) \geq \delmax/50$, and let $M=\{u_1v_1,\dots,u_{\ell}v_{\ell}\}$ be a matching of $\overline{D_i}$ of size exactly $\ell = \lceil \delmax/50 \rceil$. Note that, for every $j\in [\ell]$, there is at most one arc between $u_j$ and $v_j$ in $D$. 

        Let $\phi$ be a $(\Delta-1)$-dicolouring of $D-X_i$, which exists as $D$ is $\Delta$-critical. We will show that $\phi$ can be extended into a $(\Delta-1)$-dicolouring of $D$, thus yielding a contradiction. We first colour the vertices covered by $M$ by sequentially choosing, for every $j\in [\ell]$, a common colour $c_j\in [\Delta-1]$ for both $u_j$ and $v_j$ that is not already appearing in $N^+(u_j) \cup N^+(v_j)$. Note that this does not create any monochromatic directed cycle as there exists at most one arc between $u_j$ and $v_j$. Let us show that such a colour $c_j$ exists. By Lemma~\ref{lemma:dense_decomposition}\ref{enum:dense_decomposition:3}, both $u_j$ and $v_j$ have at most $\delmax/100$ out-neighbours in $V(D)\setminus X_i$. Moreover, at most $j-1$ colours appear in $X_i$ when colouring $\{u_j,v_j\}$ (namely $c_1,\dots, c_{j-1}$). Hence, the number of colours that appear in $N^+(u_j) \cup N^+(v_j)$ is at most 
        \[ \frac{1}{50} \delmax + \ell-1 \leq \frac{1}{25} \delmax < \Delta-1.\]
        Once the vertices covered by $M$ are coloured, we colour the remaining vertices of $X_i$ greedily, by choosing for each of them a colour of $[\Delta-1]$ that is not appearing in its out-neighbourhood. We claim that this is always possible, so assume for a contradiction that,
        for some $x\in X_i\setminus V(M)$, all $\Delta-1$ colours appear in $N^+(x)$ when $x$ is considered.
        Let $\iota(x)$ be the number of indices $j$ such that $\{u_j,v_j\} \subseteq N^+(x)$. Since any two vertices matched in $M$ received the same colour, the number of colours appearing in $N^+(x)$ is at most $d^+(x) - \iota(x)$. We thus have $d^+(x) - \iota(x) \geq \Delta-1$, which together with Lemma~\ref{lemma:mindeg} implies 
        $\iota(x) \leq 2$.
        Therefore, at least $\ell - 2$ vertices covered by $M$ do not belong to $N^+(x)$. Together with Lemma~\ref{lemma:dense_decomposition}\ref{enum:dense_decomposition:1}, this implies 
        \begin{align*}
            |N^+(x) \cap X_i| &\leq |X_i| - \ell + 2
                < \frac{98}{100}\delmax + 4d + 4
                < \frac{99}{100}\delmax.
        \end{align*}
        By Lemma~\ref{lemma:dense_decomposition}\ref{enum:dense_decomposition:3}, $x$ does not belong to $X_i$, a contradiction. The claim follows.
    \end{proofclaim}

    We can now prove that indeed, $\nu(\overline{D_i}) \leq 2$.

    \begin{claim}
        \label{claim:matching_Di_2}
        $\nu(\overline{D_i}) \leq 2$.
    \end{claim}
    \begin{proofclaim}
        Let $M=\{u_1v_1,\dots,u_{|M|}v_{|M|}\}$ be a maximum matching of $\overline{D_i}$, and assume for a contradiction that $|M| \geq 3$. By Claim~\ref{claim:matching_Di_linear}, $|M| < \frac{1}{50} \delmax$.
        We denote by $Y$ the vertices of $X_i$ that are not covered by $M$. Note that, by maximality of $M$, $Y$ is a biclique of $D$. 
        We define $Z$ as the set of vertices in $Y$ that are out-neighbours of all vertices $u_1,u_2,u_3,v_1,v_2,v_3$, so formally:
        \[ Z = Y \cap \left( \bigcap_{j=1}^{3} N^+(u_j)\right) \cap \left( \bigcap_{j=1}^{3} N^+(v_j)\right).\]

        \begin{subclaim}
            \label{subclaim:size_Z}
            $|Z| \geq 3$.
        \end{subclaim}
        \begin{proofsubclaim}
            Let $u$ be any vertex covered by $M$. Recall that, by Lemma~\ref{lemma:dense_decomposition}\ref{enum:dense_decomposition:3}, 
            \[|N^+(u)\cap X_i| \geq \frac{99}{100}\delmax.\]
            Together with Claim~\ref{claim:matching_Di_linear}, this implies $|N^+(u) \cap Y| \geq \frac{95}{100}\delmax$. 
            Since, by Lemma~\ref{lemma:dense_decomposition}\ref{enum:dense_decomposition:1}, $|Y| < |X_i| \leq \delmax + 1 + 4d$, we obtain 
            \[ |Y\setminus N^+(u)| < \frac{5}{100} \delmax + 4d +1< \frac{6}{100} \delmax. \]
            Since this holds for every $u$ saturated by $M$, and in particular for every $u\in \{u_1,u_2, u_{3},v_1,v_2, v_{3}\}$, we have 
            \begin{equation}
                \label{eq:size_YsetminusZ}
                |Y \setminus Z| < \frac{36}{100} \delmax.
            \end{equation}
            Recall that the size of $Y$ is larger than $|X_i| - \frac{4}{100}  \delmax$ by Claim~\ref{claim:matching_Di_linear}, which together with Lemma~\ref{lemma:dense_decomposition}\ref{enum:dense_decomposition:1} implies 
            \begin{equation}
                \label{eq:size_Y_lower}
                |Y| > \frac{96}{100}\delmax -300d.
            \end{equation}
            Combining~\eqref{eq:size_YsetminusZ} and~\eqref{eq:size_Y_lower}, straightforward calculations imply $|Z| > \frac{1}{2}\delmax \geq 3$,
            as desired.
        \end{proofsubclaim}

    We let $\phi$ be a $(\Delta-1)$-dicolouring of $D-Y$, with the extra property that $\phi(u_j) = \phi(v_j)$ holds for every $j\in [|M|]$. We skip the proof that $\phi$ exists, it is shown by taking a $(\Delta-1)$-dicolouring of $D-X_i$ and extending it to $D-Y$ by choosing, for every $j\in [|M|]$, a colour that is not appearing in $N^+(u_j)\cup N^+(v_j)$, exactly as in the proof of Claim~\ref{claim:matching_Di_linear}.
    
    We now prove that $\phi$ can be extended into a $(\Delta-1)$-dicolouring of $D$, thus yielding the contradiction. We first greedily colour the vertices of $Y-Z$ by choosing for each of them a colour that does not appear in its out-neighbourhood. Let us show that this is possible. Since $Y$ is a biclique, and because $|Z|\geq 3$ by Subclaim~\ref{subclaim:size_Z}, note that when we consider $v\in Y-Z$, it has at least three uncoloured neighbours (namely the vertices of $Z$). 
    Since, by Lemma~\ref{lemma:mindeg}, $d^+(v) \leq \Delta+1$, this implies that the number of colours appearing in $N^+(v)$ (when considered) is at most $\Delta-2$.

    We then extend the obtained dicolouring to $D$ by choosing for every vertex $z\in Z$ a colour that does not appear in its in-neighbourhood. Again, let us show that this is possible. Since all vertices $u_1,u_2,u_3,v_1,v_2,v_3$ belong to $N^-(z)$ (by definition of $Z$), and because $\phi(u_j) = \phi(v_j)$ for every $j\in [3]$, the number of colours appearing in the in-neighbourhood of $z$ (when considered) is at most $d^-(z) - 3$, which is at most $\Delta-2$ by Lemma~\ref{lemma:mindeg}. This yields a contradiction, and concludes the proof of Claim~\ref{claim:matching_Di_2}.
    \end{proofclaim}

    We will now prove that $X_i$ is a quasi-biclique. 
    We distinguish two cases, depending on the value of $\nu(\overline{D_i})$.
    \begin{description}
        \item[Case 1:] {\it $\nu(\overline{D_i}) \leq 1$}.
        
        If $\nu(\overline{D_i}) = 0$, then $X_i$ satisfies~\ref{enum:def_quasibiclique:1}. We thus assume that $\nu(\overline{D_i}) = 1$ and we let $u$ and $v$ be two vertices with at most one arc between them. Observe first that if $(X_i \setminus \{u,v\})\subseteq N^{\pm}(u)$, then $X_i$ satisfies~\ref{enum:def_quasibiclique:4}, and $X_i$ is a quasi-biclique. We can thus assume that there exists $x_u\in X_i \setminus \{u,v\}$ such that there is at most one arc between $x_u$ and $u$. Similarly, we can assume that there exists $x_v\in X_i \setminus \{u,v\}$ such that there is at most one arc between $x_v$ and $v$.
        Since $\nu(\overline{D_i}) = 1$, we must have $x_u = x_v$, and $D$ contains every possible digon between $\{u,v,x_u\}$ and $X_i\setminus \{u,v,x_u\}$. Note that $D[\{u,v,x_u\}]$ does not contain any digon by construction. Hence, if it contains a directed cycle, it has size exactly $3$, and $X_i$ satisfies~\ref{enum:def_quasibiclique:3}.

        Assume now that $D[\{u,v,x_u\}]$ is acyclic and let us reach a contradiction. Let $\phi$ be any $(\Delta-1)$-dicolouring of $D-X_i$. We can extend $\phi$ to $\{u,v,x_u\}$ in such a way that $\phi(u) = \phi(v) =\phi(x_u)$. Again, this is possible because $D[\{u,v,x_u\}]$ is acyclic, and because $D$ contains all digons between $\{u,v,x_u\}$ and $X_i-\{u,v,x_u\}$. Now, $\phi$ can be extended to the remaining vertices of $X_i$ because each such vertex has one colour used three times in both its in- and out-neighbourhoods.
         
        \item[Case 2:] {\it $\nu(\overline{D_i}) = 2$}.

        Let $\{u_1v_1,u_2v_2\}$ be a maximum matching of $\overline{D_i}$. Let $R=\{u_1,v_1,u_2,v_2\}$ and $K$ be the set of vertices in $X_i\setminus R$ that are in the out-neighbourhood of all vertices in $R$, so formally
        \[
            K = (X_i\setminus R) \cap \left( \bigcap_{u\in R} N^+(u) \right). 
        \]
        By definition, note that $K$ is disjoint from $R$. Let $W = X_i \setminus (R \cup K)$. We are now going to show that the partition $(R,W,K)$ of $X_i$ satisfies~\ref{enum:def_quasibiclique:5}.
        By construction, we have:
        \begin{itemize}
            \item $|R| = 4$,
            \item $W \cup K$ is a biclique (as $M$ is maximum), and
            \item $R$ dominates $K$.
        \end{itemize}
        It remains to show that $K> \frac{1}{2}|X_i|$ and that every vertex in $K$ has in-degree $\Delta+1$.
         We first show the former. The proof is similar to that of Subclaim~\ref{subclaim:size_Z}. Let $y$ be any vertex in $R$. By Lemma~\ref{lemma:dense_decomposition}\ref{enum:dense_decomposition:3}, $|N^+(y)\cap X_i| \geq \frac{99}{100}\delmax$. 
        Since, by Lemma~\ref{lemma:dense_decomposition}\ref{enum:dense_decomposition:1}, $|X_i| < \delmax +1 + 4d$, we have $
        |(X_i\setminus N^+(y))\cup R| < \frac{1}{100} \delmax + 4 + 4d < \frac{2}{100} \delmax$. 
        Since this holds for every $y\in R$, we have 
        \begin{equation}
            \label{eq:size_XisetminusK2}
            |X_i \setminus K| < \frac{8}{100} \delmax.
        \end{equation}
        Recall that $\delmax - 300d < |X_i| < \delmax+4d+1$ by Lemma~\ref{lemma:dense_decomposition}\ref{enum:dense_decomposition:1}. Hence~\eqref{eq:size_XisetminusK2} implies $|K| > \frac{1}{2}|X_i|$, as desired.

        To prove that every vertex $x\in K$ has in-degree at least $\Delta+1$, assume for a contradiction that $d^-(x) \leq \Delta$ for some $x\in K$.
        Let $\phi$ be a $(\Delta-1)$-dicolouring of $D-(K\cup W)$, with the extra property that $\phi(u_1) = \phi(v_1)$ and $\phi(u_2) =\phi(v_2)$.

        Let $y$ be any vertex in $K$ distinct from $x$. Since two colours (namely $\phi(u_1)$ and $\phi(u_2)$) are repeated in the in-neighbourhood of both $x$ and $y$, and because $d^-(x) \leq \Delta$, one can extend $\phi$ to $D$ by first colouring the vertices of $(K\cup W)\setminus \{x,y\}$ (choosing for each of them a colour that is not appearing either in its in- or out-neighbourhood), and then colouring $y$ and $x$ in this order. This yields the desired contradiction, shows that $X_i$ is a biclique of type~\ref{enum:def_quasibiclique:5}, and concludes the proof of Lemma~\ref{lemma:quasiclique_decomposition}. \qedhere
    \end{description}
\end{proof}

From now on, we fix a decomposition $X_1\sqcup \ldots \sqcup X_t\sqcup S$ of $V(D)$ satisfying the properties of Lemma~\ref{lemma:quasiclique_decomposition}.
For the next lemma, we need a few additional technical definitions. Let $X$ be a quasi-biclique of $D$, and let $r$ be any positive real number. We define $\zSCR_r(X)$ as the set of vertices in $V(D) \setminus X$ with less than $r$ in- and out-neighbours in $X$, that is
\[
    \zSCR_r(X) = \left\{ z \in V(D) \setminus X :  \max\Big(|N^-(z) \cap X|,|N^+(z) \cap X|\Big) < r \right\}. 
\]
The {\bf saving-part $\sSCR(X)$} of $X$ is defined as:
\begin{itemize}
    \item $X$ itself if $X$ is of type~\ref{enum:def_quasibiclique:1};
    \item $X\setminus \{x,y\}$ if $X$ is of type~\ref{enum:def_quasibiclique:2};
    \item $X\setminus V(\vec{C}_3)$ if $X$ is of type~\ref{enum:def_quasibiclique:3};
    \item $X\setminus \{x\}$ if $X$ is of type~\ref{enum:def_quasibiclique:4}; and
    \item $K$ if $X$ is of type~\ref{enum:def_quasibiclique:5}.
\end{itemize}
We say that a vertex $x\in \sSCR(X)$ is an {\bf $r$-saviour} of $X$ if one of the following holds:
\begin{enumerate}[label=(\alph*)]
    \item $d^+(x) = \Delta-1$ and $\zSCR_r(X) \cap N^+(x) \neq \emptyset$;
    \label{enum:def_saviour_a}
    \item $d^-(x) = \Delta-1$ and $\zSCR_r(X) \cap N^-(x) \neq \emptyset$;
    \label{enum:def_saviour_b}
    \item $d^+(x) \leq \Delta$ and $|\zSCR_r(X) \cap N^+(x)| \geq 2$; or
    \label{enum:def_saviour_c}
    \item $d^-(x) \leq \Delta$ and $|\zSCR_r(X) \cap N^-(x)| \geq 2$.
    \label{enum:def_saviour_d}
\end{enumerate}

Intuitively speaking, if $r$ is small enough compared to $\Delta$, when applying the pseudo-random colouring process of the next section, a saviour will have a quite high probability to be uncoloured and to have either two repeated colours in its out- or in-neighbourhood (if (c) or (d) holds) or only one repeated colour but one missing arc (if (a) or (b) holds). If this property holds for at least three saviours of $X$, and if these saviours are in the out-neighbourhood of every vertex of $X$, we will be able to extend a partial dicolouring of $D$ to $X$.

The following two lemmas say that, in a minimum counterexample to Theorem~\ref{thm:directed_bk}, large quasi-bicliques indeed have many saviours.
From now on, we fix $r=\log^4 \Delta$.

\begin{lemma}
    \label{lemma:linear_number_saviours_type5}
    For every $i\in [t]$, if $X_i$ is of type~\ref{enum:def_quasibiclique:5} with partition $(R_i,K_i,W_i)$, then elements of $K_i$ are $r$-saviours.
\end{lemma}
\begin{proof}
    We fix $i\in [t]$ such that $X_i$ is of type~\ref{enum:def_quasibiclique:5}.
    For the sake of conciseness, we denote $X_i,R_i,K_i$, and $W_i$ by $X,R,K$, and $W$ respectively. We further denote $\sSCR(X)$ by $\sSCR$ and denote $\zSCR_{r}(X)$ by $\zSCR_{r}$.
    By definition, recall that $|R|\leq 4$, $R$ dominates $K$, $K\cup W$ is a biclique, $|K| > \frac{1}{2}|X|$, and vertices in $K$ have in-degree $\Delta+1$. We claim that all vertices in $K$ satisfy~\ref{enum:def_saviour_a}, thus implying the statement.
    
    Assume for a contradiction that $s\in K$ does not satisfy~\ref{enum:def_saviour_a}. By definition of $\Delta$ and because $d^-(s) = \Delta+1$, we have $d^+(s) = \Delta-1$. Since $d^-(s) \neq d^+(s)$, $s$ is adjacent to a simple arc, and by Lemma~\ref{lemma:property_simplearc} we know that $s$ has a simple out-neighbour $y$. Since $R$ dominates $K$, and $K\cup W$ is a biclique, we thus necessarily have $y\in V(D)\setminus X$. As $s$ does not satisfy~\ref{enum:def_saviour_a}, $y\notin \zSCR_r$, that is, $y$ has at least $r$ out-neighbours or $r$ in-neighbours in $X$. Let $u_y$ be any vertex in $(K\cup W)\setminus N^+(y)$ distinct from $s$. Note that $u_y$ exists because $y\notin X$, so $|N^+(y)\cap X|<\frac{100}{101}\Delta$, while $|K\cup W|=|X|-|R|>\Delta-300d-4$.
    
    Let $\phi$ be a $(\Delta-1)$-dicolouring of $D-(X\cup \{y\})$. We claim that there exists a colour $c$ that is not appearing either in $N^-(y) \cup N^-(u_y)$ or in $N^+(y) \cup N^+(u_y)$. 
    To see this, assume that $y$ has $r$ out-neighbours in $X$. Then $y$ has $r$ uncoloured out-neighbours, and $u_y$ has at least $|K\cup W|-1 \geq |X|-5$ uncoloured out-neighbours. The total number of colours appearing in $N^+(y) \cup N^+(u_y)$ is thus at most
     \[
        (\Delta+1)-r + (\Delta+1)-(\Delta-301d) = \Delta -\log^4(\Delta) +301\log^3(\Delta) +2,
    \]
    which is less than $\Delta-1$ when $\Delta$ is large enough.
    Similarly, if $y$ has $r$ in-neighbours in $X$, then 
    at most $\Delta+1-r$ colours appear in its in-neighbourhood, while at most $301d$ colours appear in the in-neighbourhood of $u_y$.
    
    We thus extend $\phi$ by first choosing colour $c$ for both $y$ and $u_y$. Then, we greedily colour the vertices of $R$ and $W$ in this order. This is possible because $K$ is large, uncoloured, and in the out-neighbourhood of all vertices of $X$. We then colour the vertices of $K$, keeping $s$ for the end. Since the vertices in $K\setminus \{s\}$ have out-degree $\Delta-1$ and an uncoloured out-neighbour, this is always possible for them. For $s$, this is possible because $d^+(s) = \Delta-1$ and one colour (namely $c$) is repeated in its out-neighbourhood. This yields a contradiction.
\end{proof}

\begin{lemma}
    \label{lemma:linear_number_saviours_type1-4}
    For every $i\in [t]$, if $X_i$ is of type~\ref{enum:def_quasibiclique:1}--\ref{enum:def_quasibiclique:4}, it has at least $\frac{1}{2}|X_i|$ $r$-saviours.
\end{lemma}
\begin{proof}
    We fix $i\in [t]$ and assume that $X_i$ is of type~\ref{enum:def_quasibiclique:1}--\ref{enum:def_quasibiclique:4}. 
    Using ideas analogous to the one used in the proof of Lemma~\ref{lemma:linear_number_saviours_type5}, we first show a collection of structural results around non-$r$-saviour vertices. We then leverage these to show that many vertices must indeed be $r$-saviours. 
    
    For the sake of conciseness, again, we denote $X_i$ by $X$ and $\sSCR(X)$ by $\sSCR$. Similarly, for every $r'$, we denote $\zSCR_{r'}(X)$ by $\zSCR_{r'}$. We further denote by $\rSCR$ the set of $r$-saviours of $X$. 
    Finally, we let $Y$ be $V(D) \setminus X$, unless $X$ is of type~\ref{enum:def_quasibiclique:4} in which case $Y$ is $(V(D)\setminus X) \cup \{x\}$, where $x$ is the special vertex witnessing that $X$ is of type~\ref{enum:def_quasibiclique:4}.

    Along the proof, we will use many times that $|\sSCR |\geq |X|-3$ and that every vertex $s\in \sSCR$ is linked with digons to all vertices in $X\setminus (Y\cup \{s\})$. Note that these two properties hold because $X$ is of type~\ref{enum:def_quasibiclique:1}--\ref{enum:def_quasibiclique:4}. We further have the following.

    \begin{claim}
        \label{claim:size_X}
        Every vertex $s\in \sSCR$ satisfies $|N(s) \cap (X\setminus Y)| \leq \Delta-2$.
    \end{claim}
    \begin{proofclaim}
        If $X$ is of type~\ref{enum:def_quasibiclique:1} or~\ref{enum:def_quasibiclique:4}, this is because $D$ does not contain any biclique of size $\Delta$.
        Else, if $X$ is of type~\ref{enum:def_quasibiclique:2}, this is a consequence of Lemma~\ref{lemma:no_KD_minus_uv}.
        Finally, if $X$ is of type~\ref{enum:def_quasibiclique:3}, this is a consequence of Lemma~\ref{lemma:no_C3_K_deltaminus3}.
    \end{proofclaim}

    Intuitively, and similarly to the proof of Lemma~\ref{lemma:linear_number_saviours_type5}, when some vertex $y\in Y$ has many neighbours in $(X\setminus Y)$, we mainly want to argue that we can extend a $(\Delta-1)$-colouring of $D-(X\cup \{y\})$ by first colouring $y$ and a vertex $u_y \in X\setminus Y$ with a same colour, and then greedily colour the remaining vertices in $X$ in a specific order. The first step consists of showing that actually, we can find such a vertex $u_y$ that is not linked with a digon to $y$.
    
    \begin{claim}
        \label{claim:existence_uy}
        Let $s\in \sSCR$ and $y\in Y$ be a neighbour of $s$. There exists a vertex $u_y \in X\setminus Y$, distinct from $s$, such that $\dig{y,u_y}$ is not a digon of $D$. Moreover, if $y_1,y_2 \in Y$ are two distinct neighbours of $s$, then there exist distinct vertices $u_{1}, u_{2} \in X\setminus Y$, both distinct from $s$, such that none of $\{\dig{y_1,u_{1}},\dig{y_2,u_{2}}\}$ are digons of $D$.
    \end{claim}
    \begin{proofclaim}
        We first prove the first part of the claim.
        Assume that this is not the case, so $y$ is linked with a digon to all vertices in $X\setminus (Y\cup \{s\})$. Since vertices outside $X$ have less than $\frac{100}{101}\Delta$ out-neighbours in $X$, and because $|X| > \Delta - 300d$, we have $y\in X$. Since $y$ is chosen in $Y$, it means that $Y\cap X \neq \emptyset$, hence $X$ is of type~\ref{enum:def_quasibiclique:4}, and $y$ is the special vertex such that $X\setminus \{y\}$ is a biclique of $D$. Since $y$ is a neighbour of $s$, adding at most one arc between $y$ and $s$ to $D[X]$ yields a complete digraph. Hence $X$ is either of type~\ref{enum:def_quasibiclique:1} or~\ref{enum:def_quasibiclique:2}, a contradiction (recall that $X$ being of type ($\rm i$) implies that $X$ is not of type ($\rm j$) for any $\rm j <i$).

        To show the second part, note that at most one of $y_1,y_2$ belongs to $X$, so assume by symmetry $y_2\notin X$. We apply the first part of the statement to $y_1$ and thus find $u_{1}$. Since $y_2\notin X$, we can let $u_{2}$ be any vertex of $X\setminus (N^+(y_2) \cup Y\cup  \{s,u_{1}\})$, which is non-empty as $y_2$ has less than 
        $\frac{100}{101}\Delta$ out-neighbours in $X$, and $|X| > \Delta-300d$.
    \end{proofclaim}

    From now on, let us show some structure around vertices in $\sSCR\setminus \rSCR$, {\it i.e.} non $r$-saviours of $\sSCR$. 
    
    \begin{claim}
        \label{claim:bad_vertices_are_balanced}
        Every vertex $s\in \sSCR \setminus \rSCR$ satisfies $d^-(s)= d^+(s) =\Delta$.
    \end{claim}
    \begin{proofclaim}
        Assume for a contradiction that $\min(d^+(s),d^-(s)) = \Delta-1$. The proof is the same as that of Lemma~\ref{lemma:linear_number_saviours_type5}, but we give it for completeness.
        
        If $d^+(s) = \Delta-1$, we let $y$ be any out-neighbour of $s$ in $Y$, otherwise $d^-(s) = \Delta-1$ and we let $y$ be any in-neighbour of $s$ in $Y$. In both cases, $y$ exists by Claim~\ref{claim:size_X}. 
        Since $s\notin \rSCR$, we know that $y\notin \zSCR_r$, for otherwise $s$ would satisfy either~\ref{enum:def_saviour_a} or~\ref{enum:def_saviour_b}. By definition of $\zSCR_r$, either $y\in X$, $|N^-(y)\cap X|\geq r$, or $|N^+(y) \cap X| \geq r$. If $y\in X$, then $X$ being a quasi-biclique implies $|N^+(y) \cap X| \geq r$, so we can ignore this case. Let us thus assume that $|N^+(y)\cap X|\geq r$ or $|N^-(y)\cap X|\geq r$.

        Let $u_y$ be any vertex in $X\setminus (\{s\}\cup Y)$ such that $\dig{y,u_y}$ is not a digon of $D$, which exists by Claim~\ref{claim:existence_uy}.
        Let $v$ be any vertex in $(N(y)\cap \sSCR) \setminus \{u_y,s\}$, which exists as $|\sSCR|\geq |X|-3$.
        
        Let $\phi$ be a $(\Delta-1)$-dicolouring of $D-(X\cup \{y\})$. We greedily extend $\phi$ to $D$, thus yielding the contradiction. We first give a same colour $c$ to both $y$ and $u_y$ that is not appearing in $N^+(y) \cup N^+(u_y)$ or $N^-(y) \cup N^-(u_y)$. This is possible because $y$ has at least $r$ uncoloured in- or out-neighbours and $u_y$ is linked with digons to at least $|X|-2>\Delta-301d$ uncoloured vertices. We then greedily colour the vertices of $X$, starting with the possible vertex in $Y\cap X$, and keeping $v$ and $s$ for the end (in this order). Since both $v$ and $s$ are linked with digons to all vertices in $X\setminus Y\setminus\{v,s\}$, it is possible to choose for each of them. For $v$, if $d^+(v)= \Delta-1$ or $d^-(v)= \Delta-1$, then we can pick a colour that is not appearing in $N^+(v)$ or in $N^-(v)$ as $\dig{v,s}$ is a digon and $s$ is uncoloured. Else, $d^-(v) = d^+(v) = \Delta$, $s$ is an uncoloured in-neighbour of $v$, and $v$ has a repeated colour (namely $c$) in its in- or out-neighbourhood. We can thus choose a colour that is not appearing in $N^-(v)$ or in $N^+(v)$. We finally consider $s$, for which either $d^+(s) =\Delta-1$ and one colour (namely $c$) is repeated in $N^+(s)$, or $d^-(s) =\Delta-1$ and $c$ repeated in $N^-(s)$. In both cases, we can extend $\phi$ to $D$, a contradiction.
    \end{proofclaim}

    For vertices $s$ with in-degree and out-degree $\Delta$, the argument used above does not work anymore, as we can only save one colour for $s$. However, if two neighbours $y_1,y_2$ outside $X\setminus Y$ indeed have a lot of neighbours in $X$, then we can manage to save two colours for $s$, and the argument works.

    \begin{claim}
        \label{claim:no_two_bad_z}
        For every vertex $s\in \sSCR$ with $\max(d^+(s),d^-(s))\leq \Delta$,
        \[
            |(N^+(s) \cap Y) \setminus \zSCR_{r-3}|\leq 1 \mbox{~~~and~~~} |(N^-(s) \cap Y) \setminus \zSCR_{r-3}| \leq 1.
        \]
    \end{claim}
    \begin{proofclaim}
        Assume for a contradiction that some $s\in \sSCR$ with $d^+(s) \leq \Delta$ has two distinct out-neighbours $y_1,y_2\in Y$, both of them having at least $r-3$ in- or out-neighbours in $X$. The case of $s$ having $d^-(s) \leq \Delta$ and two in-neighbours $y_1,y_2\in Y$, none of them belonging to $\zSCR_{r-3}$, is proved analogously.

        Let $\phi$ be any $(\Delta-1)$-dicolouring of $D-(X\cup \{y_1,y_2\})$. We will show that $\phi$ extends to $D$, thus yielding the contradiction. By Claim~\ref{claim:existence_uy}, we let $u_{1}$ and $u_{2}$ be two distinct vertices of $X\setminus Y$, both distinct from $s$, such that none of $\{\dig{y_1,u_{1}},\dig{y_2,u_{2}}\}$ are digons of $D$.
        We finally let $v$ be any vertex in $N(y_1)\cap \sSCR$ distinct from $u_1$, $u_2$, and $s$. Note that we excluded the case of $X$ being of type~\ref{enum:def_quasibiclique:5}, so $|\sSCR|\geq |X|-3$, and $v$ exists as $|N(y_1)\cap (\sSCR\setminus \{u_1,u_2,s\})| \geq r-9 >0$. 
        By definition, observe that $\dig{s,v}$ is a digon of $D$, and that $D$ contains all possible digons between $\{s,v\}$ and $\{u_1,u_2\}$.
        
        We first extend $\phi$ by choosing for $y_1$ and $u_1$ a common colour that is not appearing in $N^+(y_1)\cup N^+(u_1)$ if $y_1$ has $r-3$ out-neighbours in $X$, or in $N^-(y_1)\cup N^-(u_1)$ otherwise. Similarly, we choose for $y_2$ and $u_2$ a common colour that is not appearing either in $N^+(y_2)\cup N^+(u_2)$ or $N^-(y_2)\cup N^-(u_2)$. To see that this is possible, note that both $y_1$ and $y_2$ have at least $r-4$ uncoloured in- or out-neighbours, and that both $u_1$ and $u_2$ have at least $\Delta-301d$ uncoloured in- and out-neighbours (the ones in $X\setminus\{u_1\}$). The number of forbidden colours is thus at most
        \[
            (\Delta+1)-(r-4) + (\Delta+1)-(\Delta-301d) = \Delta -\log^4(\Delta) +301\log^3(\Delta) +6,       
        \]
        which is smaller than $\Delta-1$ when $\Delta$ is large enough.
        We finally greedily colour the vertices in $X$ that are not already coloured, starting with the possible vertex in $Y$ (if it is not $y_1$ or $y_2$), and finishing with $v$ and $s$ in this order. The possible vertex in $Y$ can be coloured because it has at least $\frac{99}{100}\Delta-2$ uncoloured out-neighbours. Vertices in $X\setminus (Y\cup \{v,s\})$ can be coloured because they are linked with digons to both $v$ and $s$, which are uncoloured. 
        
        At the end, $v$ has a colour repeated in its in- or out-neighbourhood (the colour of $u_1$ and $y_1$) and one uncoloured neighbour $s$ to which $v$ is linked with a digon. Therefore, a colour is missing either in $N^-(v)$ or in $N^+(v)$. Finally, two colours are repeated in the out-neighbourhood of $s$ (namely the colour of $y_1$ and $u_1$, and the colour of $y_2$ and $u_2$). Since we assumed $d^+(s) \leq \Delta$, we can actually colour $s$, yielding the contradiction.
    \end{proofclaim}

    Recall that $|\sSCR| \geq |X| - 3$ and that every vertex $s\in \sSCR$ is linked with digons to all vertices in $X\setminus Y$.
    Therefore, since $|X| > \Delta-300d$,  the statement follows from Claims~\ref{claim:bad_vertices_are_balanced} and~\ref{claim:no_two_bad_z} when every vertex $s$ of $\sSCR$ satisfies $|N^{\pm}(s)\cap (X\setminus Y)|<\Delta-2$ (note that $\zSCR_{r-3} \subseteq \zSCR_r$).
    Since vertices in $\sSCR$ have at most $\Delta-2$ neighbours in $X\setminus Y$ by Claim~\ref{claim:size_X}, we are done unless one of the following holds:
    \begin{enumerate}
        \item $D[X\setminus Y]$ is a biclique of size $\Delta-1$ (if $X$ is of type~\ref{enum:def_quasibiclique:1} or~\ref{enum:def_quasibiclique:4});
        \item $D[X] \cup \{xy\} = \bid{K}_{\Delta-1}$ for some $x,y\in X$ (if $X$ is of type~\ref{enum:def_quasibiclique:2}); or
        \item $D[X] = \vec{C_3} \boxplus \bid{K}_{\Delta-4}$ (if $X$ is of type~\ref{enum:def_quasibiclique:3}).
    \end{enumerate}
    
    From now on, we thus assume that we are in one of the three cases above. Furthermore, combining Claims~\ref{claim:size_X},~\ref{claim:bad_vertices_are_balanced}, and~\ref{claim:no_two_bad_z}, we obtain that every vertex $s \in \sSCR \setminus \rSCR$ has exactly two in-neighbours $y_s^-,z_s^-$ and exactly two out-neighbours $y_s^+,z_s^+$ in $Y$ (we may have $y_s^- = y_s^+$ and/or $z_s^-=z_s^+$). Furthermore, as $s\notin \rSCR$ and by Claim~\ref{claim:no_two_bad_z}, exactly one of $\{y_s^-,z_s^-\}$ belongs to $\zSCR_{r-3}$, and the other one has at least $r$ in- or out-neighbours in $X$. Similarly, exactly one of $\{y_s^+,z_s^+\}$ belongs to $\zSCR_{r-3}$, and the other one has at least $r$ in- or out-neighbours in $X$.
    
    Henceforth, for every vertex $s\in \sSCR \setminus \rSCR$, we thus denote by $z_s^-$ the unique vertex in $N^-(s) \cap \zSCR_{r-3}$ and by $y_s^-$ the other vertex in $N^-(s) \cap Y$. Similarly, we denote by $z_s^+$ the unique vertex in $N^+(s) \cap \zSCR_{r-3}$ and by $y_s^+$ the other vertex in $N^+(s) \cap Y$. It turns out that $z_s^-$ and $z_s^+$ are necessarily the same vertex.
    
    \begin{claim}
        \label{claim:z_s_digon}
        For every vertex $s\in \sSCR \setminus \rSCR$, $z_s^- = z_s^+$.
    \end{claim}
    \begin{proofclaim}
        Assume for a contradiction that, for some $s\in \sSCR \setminus \rSCR$, $z_s^-s$ and $sz_s^+$ are two simple arcs.
        Let $D'$ be the digraph obtained from $D$ by removing $X \cup \{y_s^-,y_s^+\}$ and adding the arc $z_s^-z_s^+$ (note that we may have $y_s^-=y_s^+$). Note that $\deltil(D')\leq \deltil(D)$, and that by Lemma~\ref{lemma:no_KD_minus_uv}, $D'$ does not contain any copy of $\bid{K}_\Delta$ nor $\vec{C_3}\boxplus \bid{K}_{\Delta-2}$. The minimality of $D$ implies the existence of a $(\Delta-1)$-dicolouring $\phi$ of $D'$.

        Hence, $\phi$ is a dicolouring of $D-(X\cup \{y_s^-,y_s^+\})$. Let us show that it extends to $D$, thus yielding the contradiction. Let $u_y$ be any vertex in $\sSCR$ such that $\dig{y_s^-,u_y}$ is not a digon of $D$, the existence of which is guaranteed by Claim~\ref{claim:existence_uy}. We first extend $\phi$ by colouring $y_s^-$ and $u_y$ with a common colour $c$ that is not appearing either in $N^+(y_s^-) \cup N^+(u_y)$ or in $N^-(y_s^-) \cup N^-(u_y)$. If $y_s^- \neq y_s^+$, we further colour $y_s^+$ with a colour that is not appearing in $N^+(y_s^+) \cup \{z_s^-\}$ or in $N^-(y_s^-) \cup N^-(u_y)$.

        Next, we greedily extend the obtained dicolouring to $D-\{s,v\}$ by choosing for every vertex a colour that is not appearing either in its in- or out-neighbourhood, where $v$ is an arbitrary fixed vertex in $N(y_s^-) \cap (\sSCR\setminus \{s,u_y\})$. Since $s$ is uncoloured, $\dig{s,v}$ is a digon, and one colour is repeated in either $N^-(v)$ or $N^+(v)$, we can further extend this dicolouring to $D-s$. Finally, if one colour is not appearing in $N^-(s)$, we choose this colour for $s$, and we are done. If not, every colour of $[\Delta-1]$ is appearing exactly once in $N^-(s)$, except the colour of $u_y$ and $y_s^-$ that is appearing twice exactly. We thus choose colour $\phi(z_s^-)$ for $s$. We claim that the obtained colouring $\psi$ is a dicolouring of $D$, so assume for a contradiction that $D$, coloured with $\psi$, contains a monochromatic directed cycle $\cSCR$. 
        
        Since $X$ is of type~\ref{enum:def_quasibiclique:1} to~\ref{enum:def_quasibiclique:4}, we have 
        $(X\setminus (Y\cup\{s\})) \subseteq N^{\pm}(s)$. Moreover, we explicitly fixed $\psi(y_s^+)$ to a colour distinct from $\psi(z_s^-)$, so $V(\cSCR)\setminus V(D')$ is exactly $\{s\}$. In particular, the successor of $s$ in $\cSCR$ is $z_s^+$. Since $\psi$ is equal to $\phi$ on $V(D')$, this implies that $D'$, coloured with $\phi$, contains a monochromatic path from $z_s^+$ to $z_s^-$. Since $D'$ contains the arc $z_s^-z_s^+$, this is a contradiction.
    \end{proofclaim}

    Henceforth, for better readability, given a vertex $s\in \sSCR \setminus \rSCR$, we thus let  $z_s = z_s^- = z_s^+$.

    \begin{claim}
        \label{claim:out_neighbourhood_zs}
        For every vertex $s\in \sSCR \setminus \rSCR$, $|N^+(z_s) \cap X| \leq 3$ and $|N^-(z_s) \cap X| \leq 3$.
    \end{claim}
    \begin{proofclaim}
        Assume for a contradiction that, for some $s\in \sSCR \setminus \rSCR$, $|N^+(z_s) \cap X| \geq 4$, the other case being proved symmetrically.
        
        Let $u_y$ be any vertex in $X\setminus Y$, distinct from $s$, such that $\dig{y_s^-,u_y}$ is not a digon of $D$, the existence of which is guaranteed by Claim~\ref{claim:existence_uy}. If possible, we further choose $u_y$ such that $u_y\notin N(y_s^-)$.
        We then let $u_z$ be any vertex in $N(y_s^-)\cap \sSCR$ distinct from vertices of $N^+(z_s) \cup \{u_y\}$ (note that it is also distinct from $s$ as $s\in N^+(z_s)$). 
        Let us justify that $u_z$ exists. Recall that $z_s\in \zSCR_{r-3}$, and $|\sSCR \setminus X|\leq 3$, hence
        \[
            \Big| (N(y_s^-) \cap \sSCR) \setminus (N^+(z_s) \cup \{u_y\}) \Big|  > \Big| N(y_s^-) \cap X \Big| - r - \Big| N(y_s^-) \cap \{u_y\} \Big|.
        \]
        If $u_y$ is not a neighbour of $y_s^-$, this is non-negative as $| N(y_s^-) \cap X|\geq r$. If $u_y$ is a neighbour of $y_s^-$ then, by choice of $u_y$, $(X\setminus Y)\subseteq N(y_s^-)$, which implies $|N(y_s^-) \cap X|\geq |X|-1 \geq \Delta -2$ and, again, $|(N(y_s^-) \cap \sSCR) \setminus (N^+(z_s) \cup \{u_y\})|>0$. This shows the existence of $u_z$.
        
        We finally let $v$ be any vertex in $N(y_s^-)\cap \sSCR$ that is distinct from $s$, $u_y$, and $u_z$.
        Let $\phi$ be a $(\Delta-1)$-dicolouring of $D-(X\cup \{y_s^-, z_s\})$. Let us show that it extends to $D$, thus yielding the contradiction.
        
        We first give to $z_s$ and $u_z$ a common colour $c$ that is not appearing in $N^+(z_s)$, and that is neither appearing in one of $N^-(u_z), N^+(u_z)$. Note that this does not create any monochromatic directed cycle because $zu_z\notin A(D)$. To see that such a colour $c$ exists, note that at most $\Delta-3$ colours appear in $N^+(z_s)$, and at most one colour appears either in $N^-(u_z)$ or in $N^+(u_z)$. For $u_z$, this is because $u_z\in \sSCR\cap N(y_s^-)$, so $u_z$ is linked with digons to at least $\Delta-2$ vertices in $X$, and it is adjacent to $y_s^-$, which is uncoloured.

        We then give to $y_s^-$ and $u_y$ a common colour that is not appearing either in $N^-(u_y) \cup N^-(y_s^-)$ or in $N^+(u_y) \cup N^+(y_s^-)$, which is possible because $y_s^-$ has $r$ in- or out-neighbours in $X$.
        We finally greedily colour the vertices in $X$ that are not already coloured, starting with the possible vertex in $Y$, and finishing with $v$ and $s$ in this order. The possible vertex in $Y$ can be coloured because it has at least $\frac{99}{100}\Delta-2$ uncoloured out-neighbours. Vertices in $X\setminus (Y\cup \{v,s\})$ can be coloured because they are linked with digons to both $v$ and $s$, which are uncoloured. At the end, $v$ has a colour repeated in its in- or out-neighbourhood (the colour of $y_s^-$ and $u_y$) and one uncoloured neighbour $s$ to which $v$ is linked with a digon. Hence, $v$ can be coloured. Finally, two colours are repeated in the in-neighbourhood of $s$ (namely the colour of $y_s^-$ and $u_y$, and the one of $z_s$ and $u_z$). Since $d^-(s) = \Delta$ by Claim~\ref{claim:bad_vertices_are_balanced}, we can actually colour $s$, yielding the contradiction.
    \end{proofclaim}

    \begin{claim}
        \label{claim:number_bad_in_Np(ys)}
        For every vertex $s\in \sSCR \setminus \rSCR$, we have
        $|N(y_s^-)\cap (\sSCR \setminus \rSCR)| \leq 10$.
    \end{claim}
    \begin{proofclaim}
        Assume for a contradiction that, for some $s\in \sSCR \setminus \rSCR$, $|N(y_s^-)\cap (\sSCR \setminus \rSCR)| \geq 11$. Therefore, either $N^- = N^-(y_s^-)\cap (\sSCR \setminus \rSCR)$ or $N^+ = N^+(y_s^-)\cap (\sSCR \setminus \rSCR)$ has order at least $6$. We let $N=N^-$ if $|N^-| > |N^+|$, and $N=N^+$ otherwise.
        
        Let $u_y$ be any vertex in $X\setminus Y$, distinct from $s$, such that $\dig{y_s^-,u_y}$ is not a digon of $D$, the existence of which is guaranteed by Claim~\ref{claim:existence_uy}.
        We further let $s'$ be any vertex in 
        \[
            N  \setminus \Big( N^+(z_s) \cup \{u_y,s\} \Big),
        \]
        which is non-empty by Claim~\ref{claim:out_neighbourhood_zs}. We finally let $v$ be any vertex in $N(y_s^-) \cap (X\setminus Y)$ distinct from $s$, $u_y$, and $s'$.

        Note that $s'$ is chosen in $(\sSCR \setminus \rSCR)$, which allows us to consider vertex $z_{s'}$.
        Thus, let $D'$ be the digraph obtained from $D$ by first removing $X \cup \{y_s^-\}$, and further adding the digon $\dig{z_s, z_{s'}}$. Note that $z_s$ and $z_{s'}$ are indeed distinct because $s'\notin N^+(z_s)$, and $z_{s'}$ is an in-neighbour of $s'$ by definition. Moreover, $z_{s'}$ is distinct from $y_s^-$ by Claim~\ref{claim:out_neighbourhood_zs} applied to $s'$. By Claim~\ref{claim:z_s_digon}, both $z_s$ and $z_{s'}$ are linked with a digon to a vertex in $X$ (namely $s$ and $s'$), so in particular $\deltil(D') \leq \deltil(D)$. 
        Furthermore, by Lemma~\ref{lemma:no_KD_minus_uv}, $D'$ does not contain copy of $\bid{K}_{\Delta}$ nor $\vec{C_3} \boxplus \bid{K}_{\Delta-2}$.
        The minimality of $D$ thus implies the existence of a $(\Delta-1)$-dicolouring $\phi$ of $D'$. Let us show that it extends to $D$, thus yielding the contradiction.
        
        We first give to $s'$ colour $\phi(z_s)$. This is possible because $s'\in (\sSCR \setminus \rSCR) \cap N(y_s^-)$, so exactly one colour appears in both the in- and out-neighbourhood of $s'$, namely the one of $z_{s'}$, and $\phi(z_{s'}) \neq \phi(z_s)$ because $\dig{z_s,z_{s'}}$ is a digon of $D'$.

        We then give to $y_s^-$ and $u_y$ a common colour that is not appearing either in $N^+(u_y) \cup N^+(y_s^-)$ or in $N^-(u_y) \cup N^-(y_s^-)$. We finally greedily colour the vertices in $X$ that are not already coloured, starting with the possible vertex in $Y$, and finishing with $v$ and $s$ in this order. Again, this is always possible using the fact that we gave the same colour to both $z_s$ and $s'$, and the same to both $y_s^-$ and $u_y$.
    \end{proofclaim}

    With Claim~\ref{claim:number_bad_in_Np(ys)} in hand, we are now able to prove that almost all vertices in $\sSCR$ are in $\rSCR$. 
    To do this, let us count the following sum in two different ways:
    \[
        S = \sum_{s\in \sSCR\setminus \rSCR} |N(y_s^-) \cap \rSCR|.
    \]
    Let $s$ be any vertex in $\sSCR \setminus \rSCR$. Recall that $|X \setminus \sSCR| \leq 3$. Hence, together with Claim~\ref{claim:number_bad_in_Np(ys)} we thus have
    \[ 
    |N(y_s^-) \cap \rSCR| \geq |N(y_s^-) \cap X|-13 \geq r-13.
    \]
    On the other hand, recall that vertices $s\in \sSCR$ (and, a fortiori, vertices in $\rSCR$) satisfy $|N^{\pm}(s)\cap (X\setminus Y)| \geq \Delta-2$, because we are in one of the following cases:
    \begin{enumerate}
        \item $D[X\setminus Y]$ is a biclique of size $\Delta-1$ (if $X$ is of type~\ref{enum:def_quasibiclique:1} or~\ref{enum:def_quasibiclique:4});
        \item $D[X] \cup \{xy\} = \bid{K}_{\Delta-1}$ for some $x,y\in X$ (if $X$ is of type~\ref{enum:def_quasibiclique:2}); or
        \item $D[X] = \vec{C_3} \boxplus \bid{K}_{\Delta-4}$ (if $X$ is of type~\ref{enum:def_quasibiclique:3}).
    \end{enumerate}
    This implies in particular that vertices $s\in \sSCR$ have at most four neighbours in $Y$, {\it i.e.} $|N(s) \cap Y|\leq 4$.
    Since, for every vertex $s\in \sSCR \setminus \rSCR$, $y_s^-$ belongs to $Y$, we have $S \leq 4|\rSCR|$. Together with the inequality above, we thus have
    \[
        4|\rSCR| \geq S  \geq |\sSCR\setminus \rSCR| \cdot (r-13),
    \]
    which implies 
    \[
        |\rSCR| \cdot \left( \frac{4}{r-13} +1 \right)\geq |\sSCR| \geq |X|-3. 
    \]
    Since $r=\log^4(\Delta)$ and $|X| > \Delta - 300d$, and because $\Delta$ is arbitrarily large, we obtain $|\rSCR| > \frac{1}{2}|X|$, concluding the proof of Lemma~\ref{lemma:linear_number_saviours_type1-4}.
\end{proof}

Now that we know that each $X_i$ has a large number of saviours, we can derive the following. For every $i\in [t]$, we denote by $\hX_i$ the set of vertices of $X_i$ that are in the closed out-neighbourhood of all vertices of $X_i$, so formally
\[
    \hX_i = X_i \cap \left( \bigcap_{x\in X_i} N^+[x] \right).
\]

\begin{lemma}
    \label{lemma:disjoint_uniform_saviours}
    For every $i\in [t]$, there exists a collection $\tSCR_i$ of $\lceil \frac{1}{50r}\Delta \rceil$ $k_i$-tuples of vertices of $D$, such that every vertex appears in at most one tuple and only once within that tuple, and such that one of the following holds.
    \begin{enumerate}[label={\rm (\Alph*)}]
        \item\label{enum:disjoint_uniform_saviours_A} $k_i=2$, and every $(x,u)\in\tSCR_i$ satisfies $x\in \hX_i$, $d^+(x) = \Delta-1$, and $u \in N^+(x) \cap \zSCR_r(X_i)$.
        \item\label{enum:disjoint_uniform_saviours_B} $k_i=2$, and every $(x,u)\in\tSCR_i$ satisfies $x\in \hX_i$, $d^-(x) = \Delta-1$, and $u \in N^-(x) \cap \zSCR_r(X_i)$.
        \item\label{enum:disjoint_uniform_saviours_C} $k_i=3$, and every $(x,u,v)\in\tSCR_i$ satisfies $x\in \hX_i$, $d^+(x) \leq \Delta$, and $u,v \in N^+(x) \cap \zSCR_r(X_i)$.
        \item\label{enum:disjoint_uniform_saviours_D} $k_i=3$, and every $(x,u,v)\in\tSCR_i$ satisfies $x\in \hX_i$, $d^-(x) \leq \Delta$, and $u,v \in N^-(x) \cap \zSCR_r(X_i)$.
    \end{enumerate}
\end{lemma}
\begin{proof}
    Let us fix $i\in [t]$ and let $\hS_i = \hX_i \cap \sSCR(X_i)$. By Lemma~\ref{lemma:quasiclique_decomposition}, $X_i$ is a quasi-biclique. If it is a quasi-biclique of type~\ref{enum:def_quasibiclique:1},~\ref{enum:def_quasibiclique:2} or \ref{enum:def_quasibiclique:3}, then $\hS_i = \sSCR(X_i)$ and $|\hS_i| \geq |X_i|-3$ by definition of $\sSCR(X_i)$. If it is of type~\ref{enum:def_quasibiclique:4}, then $|\hS_i| \geq \frac{99}{100}\Delta$ by definition.

    Therefore, in cases~\ref{enum:def_quasibiclique:1}--\ref{enum:def_quasibiclique:4}, since $\Delta-300d < |X_i| \leq \Delta+1$ (the upper bound coming from the fact that $D$ does not contain a biclique of size $\Delta$), and because $\Delta$ is large enough, we have $|\hS_i| \geq \frac{98}{100}|X_i|$, and in particular $|\sSCR_i\setminus \hS_i| \leq \frac{2}{100}|X_i|$.
    By Lemma~\ref{lemma:linear_number_saviours_type1-4}, $\sSCR_i$ contains at least $\frac{1}{2}|X_i|$ $r$-saviours of $X_i$.  Using the remark above, $\hS_i$ contains at least $\frac{48}{100}|X_i|$ $r$-saviours. 

    If $X_i$ is of type~\ref{enum:def_quasibiclique:5} and $(R_i,K_i,W_i)$ is the given partition of $X_i$, then by definition $K_i = \hS_i$, and $|K_i| > \frac{1}{2}|X_i|$. Therefore, by Lemma~\ref{lemma:linear_number_saviours_type5}, $\hS_i$ contains more than $\frac{48}{100}|X_i|$ $r$-saviours again.

    Since each $r$-saviour satisfies~\ref{enum:def_saviour_a},~\ref{enum:def_saviour_b},~\ref{enum:def_saviour_c}, or~\ref{enum:def_saviour_d}, by the pigeonhole principle, let $\rm \gamma \in \{{ a,b,c,d}\}$ be such that $\hS_i$ contains at least $\frac{12}{100}|X_i|$ $r$-saviours, all satisfying property~$\rm (\gamma)$, and let $Y$ be the set of such $r$-saviours.
    We now greedily construct $\tSCR_i$ by repeating the following process as long as $Y$ is non-empty:
    \begin{enumerate}
        \item choose any vertex $x\in Y$;
        \item choose any vertex $u$ in $N^+(x) \cap \zSCR_r(X_i)$ if $\gamma\in \{\rm a,c\}$ and in $N^-(x) \cap \zSCR_r(X_i)$ otherwise;
        \item choose any vertex $v\neq u$ in $N^+(x) \cap \zSCR_r(X_i)$ if $\gamma = \rm c$ and in $N^-(x) \cap \zSCR_r(X_i)$ if $\gamma = \rm d$;
        \item if $\gamma\in \{\rm a,b\}$, add $(x,u)$ to $\tSCR_i$ and remove $N(u)$ from $Y$;
        \item else if $\gamma\in \{\rm c,d\}$, add $(x,u,v)$ to $\tSCR_i$ and remove $N(u) \cup N(v)$ from $Y$.
    \end{enumerate}
    Note that steps 2 and 3 above are always possible by definition of $x$ being an $r$-saviour.
    By construction, the obtained $\tSCR_i$ thus satisfies one of~\ref{enum:disjoint_uniform_saviours_A},~\ref{enum:disjoint_uniform_saviours_B},~\ref{enum:disjoint_uniform_saviours_C} or~\ref{enum:disjoint_uniform_saviours_D}. Furthermore, steps 4 and 5 guarantee that the elements of $\tSCR_i$ are pairwise disjoint. Finally, at each step, the size of $Y$ decreases by at most $4r$. Hence,
    \[
        |\tSCR_i| \geq \frac{3}{100r}|X_i| \geq \frac{1}{50r}\Delta,
    \]
    where in the second inequality we used $|X_i| \geq \Delta - 300d$ and the fact that $\Delta$ is large enough. The result follows.
\end{proof}

From now on, for every $i\in [t]$, we fix a collection $\tSCR_i$ as in the statement of Lemma~\ref{lemma:disjoint_uniform_saviours}. 

\begin{remark} \label{remark:size_Xstar}
Let us note, for future use, that $|\hX_i|\geq\frac{49}{100}\Delta$ for every $i\in[t]$, as is implied by the proof of Lemma~\ref{lemma:disjoint_uniform_saviours}.
\end{remark}

\subsection{The probabilistic analysis}
\label{sec:directed_bk:probabilistic}

Let $\phi$ be a partial $(\Delta-1)$-dicolouring of $D$. For an integer $i\in [t]$, we say that $\phi$ is {\bf $i$-extendable} if there exist three distinct uncoloured vertices $x\in \hX_i$, each of which satisfying one of the following:
\begin{enumerate}[label=$(\roman*)$]
    \item $d^+(x) = \Delta-1$ and one colour is repeated in $N^+(x)$,
    \item $d^-(x) = \Delta-1$ and one colour is repeated in $N^-(x)$,
    \item $d^+(x) \leq  \Delta$ and two colours are repeated in $N^+(x)$, or
    \item $d^-(x) \leq \Delta$ and two colours are repeated in $N^-(x)$.
\end{enumerate}

We say that $\phi$ is {\bf extendable} if it is $i$-extendable for every $i\in [t]$, and for every $d$-sparse vertex $s\in S$, three colours are repeated in $N^+(s)$.

Let us briefly justify that an extendable partial dicolouring of $D$ can actually be extended to $D$. We colour the uncoloured vertices of $X_1,\dots,X_t,S$ in this order. When colouring the vertices of $X_i$, we keep the three special vertices described above for the end. Since they are uncoloured and belong to $X_i^\star$, for each other vertex of $X_i$, there is a colour not appearing in its out-neighbourhood. When we end up with these three vertices, clearly we can choose a colour that does not appear in either the in- or out-neighbourhood.
At the end, we can colour each vertex $s\in S$ with a colour that does not appear in its out-neighbourhood, as three colours are repeated in it.

\medskip

The following lemma thus implies $\dic(D) \leq \Delta-1$, yielding a contradiction and concluding the proof of Theorem~\ref{thm:directed_bk}.

\begin{lemma}\label{lemma:extendable}
    The digraph $D$ admits an extendable partial $(\Delta-1)$-dicolouring.
\end{lemma}
\begin{proof}
    Let $\psi$ be a random colouring of the vertices of $D$, where each vertex is assigned a colour from $[\Delta-1]$ independently and uniformly at random. Let $V_{\psi}$ be the set of vertices $x$ such that $\psi(x)=\psi(x^-)=\psi(x^+)$ for some $x^-\in N^-(x)$ and some $x^+\in N^+(x)$.
    We consider the random partial dicolouring $\varphi$ obtained from $\psi$ by uncolouring simultaneously all the vertices in $V_{\psi}$. By construction, $\varphi$ is indeed a partial $(\Delta-1)$-dicolouring of $D$. We are going to show, using the Lov\'asz Local Lemma, that $\varphi$ is extendable with positive probability.

    Let us define the events that we want to avoid. For each $d$-sparse vertex $s\in S$, let $A_s$ be the event that fewer than three colours are repeated in $N^+(s)$. And, for each $i\in[t]$, let $B_i$ be the event that $\varphi$ is not $i$-extendable. We note that $A_s$ depends only on the colour choices for the vertices within distance 2 of $s$ in the underlying undirected graph $\UG(D)$ of $D$. Likewise, $B_i$ depends only on the colour choices for the vertices within distance 2 of $X_i$ in $\UG(D)$. Therefore, each of these events is mutually independent of all the others, except at most $(2\delmax)^5$ of them. By the Lov\'asz Local Lemma (Lemma~\ref{lemma:lll}) and Lemma~\ref{lemma:mindeg}, it will be enough to show that each of these events holds with probability at most $\Delta^{-6}$. This is the purpose of Claims~\ref{claim:events_A} and~\ref{claim:events_B}.

    \begin{claim}\label{claim:events_A} For any vertex $s\in S$, $A_s$ holds with probability at most $e^{-\log^2\Delta}$.
    \end{claim}
    \begin{proofclaim}
        We let $M_s$ be the random variable counting the number of colours $c$ such that:
        \begin{itemize}
            \item $c$ is assigned to two distinct out-neighbours $x$ and $y$ of $s$ not inducing a digon, and
            \item for every distinct $x,y\in N^+(s)$ not inducing a digon, if $x$ and $y$ are assigned colour $c$, then both of them retain the colour.
        \end{itemize}
        It suffices to show that $M_s<3$ holds with probability at most $e^{-\log^2\Delta}$. To do so, we first estimate $\EE(M_s)$ and then show that $M_s$ is concentrated around its expectation.  In what follows, we denote by $U_s$ the set of unordered pairs of distinct out-neighbours of $s$ not inducing a digon. 
        \begin{subclaim}
        \label{subclaim:events_A_1}
        $\displaystyle \EE(M_s)\geq\frac{|U_s|}{e^6(\Delta-1)}$.
        \end{subclaim}
        \begin{proofsubclaim}
            Let $\tilde M_s$ be the number of colours $c$ such that:
            \begin{itemize}
                \item $c$ is assigned to exactly two out-neighbours $x$ and $y$ of $s$,
                \item $\dig{x,y}$ is not a digon, and
                \item both $x$ and $y$ retain colour $c$.
            \end{itemize}
            Observe that $M_s\geq\tilde M_s$. Given a colour $c\in[\Delta-1]$ and two distinct out-neighbours $x$ and $y$ of $s$ such that $\dig{x,y}$ is not a digon, let us consider the event $E_{\{x,y\},c}$ that $x$ and $y$ are assigned colour $c$, no other vertex in $N^+(s)$ is assigned colour $c$, and both $x$ and $y$ retain the colour. 
            Let $\mathscr{W}_{x,y} = N^+(s)\cup N(x)\cup N(y)\setminus\{x,y\}$.
            Since $\dig{x,y}$ is not a digon, we have
            \begin{align*}
                 \PP(E_{\{x,y\},c}) &\geq\PP\Big(\text{for every $w\in \mathscr{W}_{x,y}$, }  \psi(w)\neq c =  \psi(x) = \psi(y)\Big)\\
                &\geq\left(\frac{1}{\Delta-1}\right)^2\left(\frac{\Delta-2}{\Delta-1}\right)^{5\delmax}\geq\frac{1}{e^6(\Delta-1)^2},
            \end{align*}
            where in the last inequality we used that $\Delta$ is large enough, ${\displaystyle \lim_{\Delta \to +\infty}} \left(\frac{\Delta-2}{\Delta-1}\right)^{5\Delta} = \frac{1}{e^5} > \frac{1}{e^6}$, and that $\Delta \leq \delmax \le \Delta+1$.
            By linearity of the expectation, we conclude that
            \[
            \EE(M_s)\geq\EE(\tilde M_s)\geq (\Delta-1)\cdot|U_s|\cdot\frac{1}{e^6(\Delta-1)^2}=\frac{|U_s|}{e^6(\Delta-1)}.
            \]
        \end{proofsubclaim}
    \begin{subclaim}
    \label{subclaim:events_A_2} 
    For any $\lambda>504\sqrt{\frac{|U_s|}{\Delta-1}} + 2752$, we have
    \[
    \PP(|M_s-\EE(M_s)|>\lambda)\leq 8\exp\left(\frac{-\lambda^2}{512\left(\frac{|U_s|}{\Delta-1}+\lambda\right)}\right).\]
    \end{subclaim}
    \begin{proofsubclaim}
        We consider two auxiliary random variables: $M^\star_s$, which counts the number of colours $c$ assigned to both $x$ and $y$ for some $\{x,y\}\in U_s$, and $ M^{\mathrm{del}}_s$, which counts the number of colours $c$ such that, for some $\{x,y\}\in U_s$, $c$ is assigned to $x$ and $y$, but not retained by both of them. Observe that $M_s=M^\star_s-M^{\mathrm{del}}_s$. Let us show that these auxiliary random variables are concentrated, hence implying that $M_s$ is concentrated as well.

        Given $\{x,y\}\in U_s$ and a colour $c\in[\Delta-1]$, let $E^\star_{\{x,y\},c}$ be the event that $x$ and $y$ are assigned colour $c$. We have that
        \[\PP(E^\star_{\{x,y\},c})=\frac{1}{(\Delta-1)^2}.\]
        Thus, by the union bound, 
        
        \begin{equation}
        \label{eq:EM^*}
        \EE(M^\star_s)\leq(\Delta-1)\cdot|U_s|\cdot\frac{1}{(\Delta-1)^2}=\frac{|U_s|}{\Delta-1}.
        \end{equation}

        We note that changing the outcome of any colour assignment affects $M^\star_s$ by at most $1$. Moreover, whenever $M^\star_s\geq k$ for some integer $k$, this can be certified by revealing $2k$ colour assignments: for each colour $c$ counted by $M^\star_s$, it suffices to exhibit two vertices $x$ and $y$ with $\{x,y\}\in U_s$ that have been assigned colour $c$. By Talagrand's inequality (Lemma~\ref{lemma:talagrand}) and~\eqref{eq:EM^*},
        \begin{equation}\label{eq:M^*}
            \PP(|M^\star_s-\EE(M^\star_s)|>\lambda)\leq 4\exp\left(\frac{-\lambda^2}{32\cdot 2\cdot(\EE(M^\star_s)+\lambda)}\right)\leq 4\exp\left(\frac{-\lambda^2}{64\left(\frac{|U_s|}{\Delta-1}+\lambda\right)}\right)
        \end{equation}
        for any $\lambda > 126\sqrt{2\frac{|U_s|}{\Delta-1}} + 688 \geq 126\sqrt{2\EE(M^\star_s)}+344\cdot 2$.

        Similarly, $M^{\mathrm{del}}_s$ is affected by at most $1$ when the outcome of any colour assignment is changed.
        Moreover, whenever $M^{\mathrm{del}}_s\geq k$, this can be certified by revealing at most $4k$ colour assignments: for each one of the $k$ colours, it suffices to exhibit two vertices $x$ and $y$ with $\{x,y\}\in U_s$, and one in-neighbour $x^-$ and one out-neighbour $x^+$ of $x$, such that $x$, $y$, $x^-$ and $x^+$ have all been assigned that colour. By Talagrand's inequality (Lemma~\ref{lemma:talagrand}),~\eqref{eq:EM^*}, and the fact that $M^{\mathrm{del}}_s\leq M^\star_s$, we have
        \begin{equation}\label{eq:M^del}
            \PP(|M^{\mathrm{del}}_s-\EE(M^{\mathrm{del}}_s)|>\lambda)
            \leq4\exp\left(\frac{-\lambda^2}{32\cdot 4\cdot(\EE(M^{\mathrm{del}}_s)+\lambda)}\right)
            \leq 4\exp\left(\frac{-\lambda^2}{128\left(\frac{|U_s|}{\Delta-1}+\lambda\right)}\right)
        \end{equation}
        for any $\lambda> 252\sqrt{\frac{|U_s|}{\Delta-1}} + 1376 \geq  126\sqrt{4\EE(M^{\mathrm{del}}_s)} + 344\cdot 4$.

        \medskip
        
        Recall that $M_s = M_s^\star-M^{\mathrm{del}}_s$. Hence, by linearity of the expectation and the triangle inequality, we have $|M_s-\EE(M_s)| \leq |M^\star_s-\EE(M^\star_s)|+|M^{\mathrm{del}}_s-\EE(M^{\mathrm{del}}_s)|$.
        Combining~\eqref{eq:M^*} and~\eqref{eq:M^del}, we thus obtain
        \begin{align*}
        \PP(|M_s-\EE(M_s)|>\lambda)
        &\leq
        \PP(|M^\star_s-\EE(M^\star_s)|+|M^{\mathrm{del}}_s-\EE(M^{\mathrm{del}}_s)|>\lambda) \\
        &\leq\PP(|M^\star_s-\EE(M^\star_s)|>\tfrac{\lambda}{2})+\PP(|M^{\mathrm{del}}_s-\EE(M^{\mathrm{del}}_s)|>\tfrac{\lambda}{2}) \\
        &\leq 4\exp\left(\frac{-\lambda^2}{256\left(\frac{|U_s|}{\Delta-1}+\lambda\right)}\right)+4\exp\left(\frac{-\lambda^2}{512\left(\frac{|U_s|}{\Delta-1}+\lambda\right)}\right)\\
        &\leq 8\exp\left(\frac{-\lambda^2}{512\left(\frac{|U_s|}{\Delta-1}+\lambda\right)}\right)
        \end{align*}
        for any $\lambda>504\sqrt{\frac{|U_s|}{\Delta-1}} + 2752$.
    \end{proofsubclaim}

    To end the proof of the claim, we note that, since $s$ is $d$-sparse, out of the $d^+(s)(d^+(s)-1)$ possible arcs between the vertices of $N^+(s)$, at least
    \[d^+(s)(d^+(s)-1)-\left(\delmax(\delmax-1)-d\delmax\right)
    \geq(\Delta-1)(\Delta-2)-(\Delta+1)\Delta+d(\Delta-1)
    \geq\frac{d\Delta}{2}\]
    of them are not present in $D$. This implies that
    $|U_s|\geq\frac{d\Delta}{4}$.
    In particular, if we take $\lambda=\frac{|U_s|}{e^6(\Delta-1)}-3$, 
    the condition of Subclaim~\ref{subclaim:events_A_2} holds. By Subclaims~\ref{subclaim:events_A_1} and~\ref{subclaim:events_A_2}, we obtain
    \[
    \PP(A_s)\leq\PP(M_s<3)\leq\PP(\EE(M_s)-M_s>\lambda)\leq 8\exp\left(\frac{-\lambda^2}{512\left(\frac{|U_s|}{\Delta-1}+\lambda\right)}\right)\leq e^{-\log^2\Delta},
    \]
    as desired.
    \end{proofclaim}
    
    \begin{claim}\label{claim:events_B} For any $i\in [t]$, $B_i$ holds with probability at most $e^{\frac{-\Delta}{r^3}}$.
    \end{claim}
    \begin{proofclaim} 
    Let $\tSCR_i$ be a set given by Lemma~\ref{lemma:disjoint_uniform_saviours} (of size exactly $\lceil\frac{\Delta}{50r}\rceil$). For the sake of conciseness, we write the proof assuming that $\tSCR_i$ satisfies~\ref{enum:disjoint_uniform_saviours_C}; the other cases can be argued analogously. Let $T_i$ be the set of vertices appearing in some triple of $\tSCR_i$. Let $M_i$ be the random variable counting the number of triples $(x,u,v)\in\tSCR_i$ such that:
    \begin{enumerate}[label=$(\arabic*)$]
        \item\label{enum:good_triple_1} vertex $x$ is uncoloured;
        \item\label{enum:good_triple_2} both $u$ and $v$ retain their colours, which both appear on some vertex of $\hX_i$; and
        \item\label{enum:good_triple_3} for every vertex of the triple, its colour is not assigned to any other vertex of $T_i$.
    \end{enumerate}
    We note that, by~\ref{enum:good_triple_1} and~\ref{enum:good_triple_2}, $\varphi$ is $i$-extendable with probability at least $\PP(M_i\geq 3)$. (In cases~\ref{enum:disjoint_uniform_saviours_B} and~\ref{enum:disjoint_uniform_saviours_D}, we take into account that $\hX_i$ is a biclique.) Subclaim~\ref{subclaim:events_B_1} shows that $M_i$ is large on average, and Subclaim~\ref{subclaim:events_B_2} that it is concentrated around its expected value. The claim is then obtained by combining these two facts. The purpose of condition~\ref{enum:good_triple_3} is to prepare the ground for using Azuma's inequality.
    \begin{subclaim}\label{subclaim:events_B_1} $\displaystyle \EE(M_i)\geq\frac{\Delta}{e^{37}r}$.
    \end{subclaim}
    \begin{proofsubclaim}
        Given a triple $(x,u,v)\in\tSCR_i$, three distinct vertices $x'\in\hX_i\setminus T_i$, $u'\in\hX_i\setminus (T_i\cup N(u))$, and $v'\in\hX_i\setminus (T_i\cup N(v))$, and three distinct colours $c_1,c_2,c_3\in [\Delta-1]$, let $E_{x,u,v,x',u',v',c_1,c_2,c_3}$ be the event that
        \begin{itemize}
            \item $c_1$ is assigned to $x$ and $x'$, but not to any vertex in $\hX_i\cup T_i\setminus\{x,x'\}$; 
            \item $c_2$ is assigned to $u$ and $u'$, but not to any vertex in $\hX_i\cup T_i\cup N(u)\cup N(u')\setminus\{u,u'\}$;
            \item $c_3$ is assigned to $v$ and $v'$, but not to any vertex in $\hX_i\cup T_i\cup N(v)\cup N(v')\setminus\{v,v'\}$.
        \end{itemize}
        Since $\dig{x,x'}$ is a digon, these conditions ensure that the triple $(x,u,v)$ is counted by $M_i$. Since $|T_i|= 3\lceil\frac{\Delta}{50r}\rceil$, we have that
        \[\PP(E_{x,u,v,x',u',v',c_1,c_2,c_3})\geq\left(\frac{1}{\Delta-1}\right)^6\left(\frac{\Delta-4}{\Delta-1}\right)^{|\hX_i|+|T_i|+8\Delta}\geq\left(\frac{1}{\Delta-1}\right)^6 e^{-30}. \]
        By Remark~\ref{remark:size_Xstar} and the fact that $u,v\in\zSCR_r(X_i)$, there are at least $\left(\frac{48}{100}\Delta\right)^3$ ways of choosing vertices $x'$, $u'$, and $v'$; moreover, there are $(\Delta-1)(\Delta-2)(\Delta-3)$ possible ways of choosing colours $c_1$, $c_2$, and $c_3$. Since for every fixed $(x,u,v)\in\tSCR_i$, the events of the form $E_{x,u,v,x',u',v',c_1,c_2,c_3}$ are pairwise disjoint, by linearity of the expectation, we have 
        \[
        \EE(M_i)\geq|\tSCR_i|\left(\frac{48}{100}\Delta\right)^3(\Delta-1)(\Delta-2)(\Delta-3)\left(\frac{1}{\Delta-1}\right)^6 e^{-30}\geq\frac{\Delta}{e^{37}r}.
        \]
    \end{proofsubclaim}
    
    \begin{subclaim}
    \label{subclaim:events_B_2} 
    For any $\lambda\geq 0$, $\PP(|M_i-\EE(M_i)|>\lambda)\leq 2\exp\left(\frac{-\lambda^2}{44\Delta}\right)$.
    \end{subclaim}
    \begin{proofsubclaim}
        We will be interested in how much $M_i$ can vary when the colour assignment of an arbitrary vertex is altered. For that, we first fix some notation. Let $\tau,\tau':V(D)\rightarrow[\Delta-1]$ be two arbitrary colour assignments that differ everywhere. We consider the event $E$ that $\psi=\tau$ and, for every vertex $w$ of $D$, the event $E_w$ that $\psi|_{V(D)\setminus\{w\}}=\tau|_{V(D)\setminus\{w\}}$ and $\psi(w)=\tau'(w)$. We let $\mSCR$ and $\mSCR_w$ be the sets of triples of $\tSCR_i$ satisfying~\ref{enum:good_triple_1},~\ref{enum:good_triple_2} and~\ref{enum:good_triple_3} when $E$ and $E_w$ hold, respectively. 
        In particular, $|\mSCR|$ and $|\mSCR_w|$ are the values $M_i(E)$ and $M_i(E_w)$.

        We are going to consider the random outcomes defining $\psi$ in a specific order, depending on the given $i\in[t]$. We fix a labelling $w_1,\ldots,w_n$ of the vertices of $D$ such that $w_1,\ldots,w_q\in V(D)\setminus(\hX_i\cup T_i)$, and $w_{q+1},\ldots w_n\in\hX_i\cup T_i$ for some $q$. For $1\leq j\leq n$, we denote by $E_j$ and $E'_j$ the events that $\psi(w_1)=\tau(w_1),\ldots,\psi(w_j)=\tau(w_j)$ and that $\psi(w_1)=\tau(w_1),\ldots,\psi(w_{j-1})=\tau(w_{j-1})$ and $\psi(w_j)=\tau'(w_j)$. 
    
        Let us now discuss how much $M_i$ can vary when changing the random assignment of an arbitrary vertex; in other words, let us fix $1\leq j\leq n$ and examine the quantity $\big||\mSCR|-|\mSCR_{w_j}|\big|$. We note that a triple $(x,u,v)\in\tSCR_i$ with $\{\tau(x),\tau(u),\tau(v)\}\cap\{\tau(w_j),\tau'(w_j)\}=\emptyset$ belongs to $\mSCR$ if and only if it belongs to $\mSCR_{w_j}$. On the other hand, by~\ref{enum:good_triple_3}, both $\mSCR$ and $\mSCR_{w_j}$ contain at most two triples $(x,u,v)\in\tSCR_i$ with $\{\tau(x),\tau(u),\tau(v)\}\cap\{\tau(w_j),\tau'(w_j)\}\neq\emptyset$. Therefore, 
        \begin{equation}
            \label{eq:azuma_general}
            \left|M_i(E)-M_i(E_{w_j})\right|\leq 2.
        \end{equation}
        For $1\leq j\leq q$, the previous bound~\eqref{eq:azuma_general} can be improved to $0$ in some situations. Let us denote by $\mathcal N_j$ the set of neighbours of $w_j$ in $\hX_i\cup T_i$, and let $F_j$ be the event that some vertex in $\mathcal N_j$ is assigned colour $\tau(w_j)$ or $\tau'(w_j)$. Let us discuss when do $M_i(E)$ and $M_i(E_{w_j})$ differ. If $\mSCR\neq\mSCR_{w_j}$, then, since $w_j$ is not in $\hX_i\cup T_i$, changing the outcome of $\psi(w_j)$ from $\tau(w_j)$ to $\tau'(w_j)$ has to affect the instruction of uncolouring a vertex $y$ from $\hX_i\cup T_i$. In particular, $y$ must be adjacent to $w_j$, and its assigned colour $\tau(y)$ under $E$ and $E_{w_j}$ must be $\tau(w_j)$ or $\tau'(w_j)$. Hence, $E,E_{w_j}\subseteq F_j$. Therefore, for every $1\leq j\leq q$, changing the outcome of $\psi(w_j)$ does not affect $M_i$ on $\overline{F_j}$. We deduce that
        \begin{equation}\label{eq:azuma_no_change}
            \EE\left(M_i\,|\,E_j\cap \overline{F_j}\right)=\EE\left(M_i\,|\,E'_j\cap\overline{F_j}\right).
        \end{equation}
    
        We also need some control over the probability of $F_j$, for every $1\leq j\leq q$. We note that, by the definition of the ordering $w_1,\ldots,w_n$, $F_j$ is independent of both $E_j$ and $E'_j$. Therefore, by the union bound,
        \begin{equation}
        \label{eq:azuma_probability_bound}
        \PP(F_j\,|\,E_j)=\PP(F_j\,|\,E'_j)=\PP(F_j)\leq\frac{2|\mathcal N_j|}{\Delta-1}
        \end{equation}
        for every $1\leq j\leq q$.
        We now bound
        \[\sum_{j=1}^n\left|\EE(M_i\,|\,E_j)-\EE(M_i\,|\, E'_j)\right|^2,\]
        which will allow us to apply Azuma's inequality. We split the sum into two chunks. For the terms $q+1\leq j\leq n$, \eqref{eq:azuma_general} yields 
        \[\sum_{j=q+1}^n\left|\EE(M_i\,|\,E_j)-\EE(M_i\,|\, E'_j)\right|^2\leq 2^2|\hX_i\cup T_i|\leq 5\Delta.\]
        For $1\leq j\leq q$, \eqref{eq:azuma_general}, \eqref{eq:azuma_no_change}, and~\eqref{eq:azuma_probability_bound} yield
        \begin{align*}
            &\sum_{j=1}^q\left|\EE(M_i\,|\,E_j)-\EE(M_i\,|\, E'_j)\right|^2 \\
            &\leq
            \sum_{j=1}^q 2\left|\EE(M_i\,|\,E_j)-\EE(M_i\,|\, E'_j)\right| \\
            &\begin{multlined}=2\sum_{j=1}^q\left|\EE\left(M_i\,|\,E_j\cap F_j\right)\PP\left(F_j\,|\,E_j\right)+\EE\left(M_i\,|\,E_j\cap\overline{ F_j}\right)\PP\left(\overline{F_j}\,|\,E_j\right)\right.\\ 
            \left.-\,\EE\left(M_i\,|\, E'_j\cap F_j\right)\PP\left(F_j\,|\,E'_j\right)-\EE\left(M_i\,|\,E'_j\cap\overline{ F_j}\right)\PP\left(\overline{F_j}\,|\,E'_j\right)\right|\end{multlined} \\
            &=2\sum_{j=1}^q\left|\EE(M_i\,|\,E_j\cap F_j)\PP(F_j\,|\,E_j)-\EE(M_i\,|\, E'_j\cap  F_j)\PP(F_j\,|\,E'_j)\right| \\
            &\leq 2\sum_{j=1}^q\left|\EE(M_i\,|\,E_j\cap F_j)-\EE(M_i\,|\,E'_j\cap F_j)\right|\frac{2|\mathcal N_j|}{\Delta-1} \\
            &\leq\frac{8}{\Delta-1}\sum_{j=1}^q|\mathcal N_j|\leq\frac{16}{\Delta-1}\delmax|\hX_i\cup T_i|\leq 17\Delta.
        \end{align*}
        Altogether,
        \begin{equation*}
            \sum_{j=1}^n\left|\EE(M_i\,|\,E_j)-\EE(M_i\,|\, E'_j)\right|^2\leq 22\Delta.
        \end{equation*}
        The statement now follows from Azuma's inequality (Lemma~\ref{lemma:azuma}).
    \end{proofsubclaim}
    
    To end the proof of the Claim~\ref{claim:events_B}, we combine Subclaims~\ref{subclaim:events_B_1} and~\ref{subclaim:events_B_2} to obtain 
    \[\PP(B_i)\leq\PP(M_i<3)\leq\PP\left(\EE(M_i)-M_i>\frac{\Delta}{e^{37}r}-3\right)\leq2\exp\left(\frac{-\left(\frac{\Delta}{e^{37}r}-3\right)^2}{44\Delta}\right)\leq e^{\frac{-\Delta}{r^3}},\]
    as desired.
    \end{proofclaim}
    It follows from Claims~\ref{claim:events_A} and~\ref{claim:events_B} that each of the bad events occurs with probability at most $\Delta^{-6}$. As discussed above, this allows us to conclude the proof of Lemma~\ref{lemma:extendable} with the Lov\'asz Local Lemma (Lemma~\ref{lemma:lll}), and thus the proof of Theorem~\ref{thm:directed_bk}.
\end{proof}

\section{Deriving the analogous result for \texorpdfstring{$\Delta^+$}{Dplus}}
\label{sec:directed_bk:plus}

As explained in the introduction, we now derive Theorem~\ref{thm:directed_bk_deltaplus} from Theorem~\ref{thm:directed_bk}. Let us first restate it here for convenience.

\thmdirectedbkplus*

\begin{proof}
    We show the result for $\DeltaDirectedBKPlus = \max(\DeltaDirectedBK,556)$.
    Let us thus fix $\Delta\geq \DeltaDirectedBKPlus$, and assume for a contradiction that there exists a digraph $D$ with $\bic(D)\leq \Delta-1$, $\delplus(D)\leq \Delta$, and $\dic(D)\geq \Delta$, and such that $D$ does not contain $\vec{C_3}\boxplus \bid{K}_{\Delta-2}$.
    Among all such digraphs, we choose $D$ such that:
    \begin{enumerate}[label=(\arabic*)]
        \item $|V(D)|$ is minimum, and
        \label{enum:proof_delplus:1}
        \item with respect to~\ref{enum:proof_delplus:1}, the number of simple arcs of $D$ is minimum.
        \label{enum:proof_delplus:2}
    \end{enumerate}
    
    We start with the following easy property of $D$ that we will need later. 
    \begin{claimlabeled}[thm:directed_bk_deltaplus]
        \label{claim:delplus:no_P3_Kdel}
        $D$ does not contain any induced copy of $\vec{P}_3 \boxplus \bid{K}_{\Delta-2}$.
    \end{claimlabeled}
    \begin{proofclaim}
        Assume for a contradiction that $X\subseteq V(D)$ induces a copy of $\vec{P}_3 \boxplus \bid{K}_{\Delta-2}$. Let $u,v,w \in X$ be the three vertices inducing the copy of $\vec{P_3}$.
        
        By minimality of $D$, let $\phi$ be a $(\Delta-1)$-dicolouring of $D-X$. We extend $\phi$ to $D$ by first giving a common colour to $u,v,w$. This is possible because $u,v,w$ have at least $\Delta-2$ uncoloured out-neighbours, so one colour is not appearing in $N^+(u) \cup N^+(v) \cup N^+(w)$.

        We then greedily extend $\phi$ to $X\setminus \{u,v,w\}$ by choosing for every such vertex $x$ a colour that is not appearing in its out-neighbourhood. This is always possible because $\{u,v,w\} \subseteq N^+(x)$ and $\phi(u) = \phi(v) = \phi(w)$.
    \end{proofclaim}

    Similarly to the proof of Theorem~\ref{thm:directed_bk}, by~\ref{enum:proof_delplus:1}, every vertex $v$ satisfies $d^+(v),d^-(v) \geq \Delta-1$, but this time we do not have any upper bound on $d^-(v)$.
    However, we show that there exists some vertex with in-degree exactly $\Delta-1$ and out-degree at least $\Delta$, for otherwise we are directly done by Theorem~\ref{thm:directed_bk}.

    \begin{claimlabeled}[thm:directed_bk_deltaplus]
        \label{claim:V_deltam1_nonempty}
        There exists $v\in V(D)$ satisfying $d^-(v) = \Delta-1 < d^+(v)$.
    \end{claimlabeled}
    \begin{proofclaim}
        Let us prove that there exists $v$ with $d^-(v) < d^+(v)$. The fact that $d^-(v) = \Delta-1$ then follows as inequality $\Delta-1 \leq d^+(v) \leq \Delta$ holds for every vertex $v$.
        Assume for a contradiction that $d^-(v) \geq d^+(v)$ for every vertex $v\in V(D)$. Since 
        \[ \sum_{v\in V(D)}d^-(v) = |A(D)| = \sum_{v\in V(D)}d^+(v), \]
        we thus have $d^+(v) = d^-(v)$ for every vertex $v\in V(D)$. In particular, $\deltil(D) =\delplus(D) \leq  \Delta$, and Theorem~\ref{thm:directed_bk} implies $\dic(D) \leq \Delta-1$, a contradiction.
    \end{proofclaim}

    Let $\mathcal{V}_{\Delta-1}$ be the set of vertices with in-degree $\Delta-1$ and out-degree at least $\Delta$. The whole proof consists in identifying a vertex $v\in \mathcal{V}_{\Delta-1}$ with a specific structure, and then justifying that a $(\Delta-1)$-dicolouring of $D-N^-[v]$ can be extended to $D$, thus yielding the contradiction. Let us first prove the following two properties of the vertices of $\mathcal{V}_{\Delta-1}$.
    \begin{claimlabeled}[thm:directed_bk_deltaplus]
        \label{claim:simple_in_out}
        Every vertex $v\in \mathcal{V}_{\Delta-1}$ has a simple in-neighbour $u$ and a simple out-neighbour $w$.
    \end{claimlabeled}
    \begin{proofclaim}
        Clearly $v$ has a simple out-neighbour $w$ as $d^+(v) > d^-(v)$. To see that $v$ has a simple in-neighbour, let $D'=D\setminus \{vw\}$. By choice of $D$, $D'$ admits a $(\Delta-1)$-dicolouring $\phi$. Clearly, $D$ coloured with $\phi$ contains a monochromatic induced directed cycle $\cSCR$, which necessarily contains $vw$. Since $\cSCR$ is induced, the predecessor of $v$ is a simple in-neighbour of $v$.
    \end{proofclaim} 

    \begin{claimlabeled}[thm:directed_bk_deltaplus]
        \label{claim:in-neighbourhood_biclique}
        For every $v\in \mathcal{V}_{\Delta-1}$, $N^-(v)$ is a biclique of $D$.
    \end{claimlabeled}
    \begin{proofclaim}
        Let us fix $v\in \mathcal{V}_{\Delta-1}$, and let $D'$ be the digraph obtained from $D$ by first removing all arcs of
        \[
            \{vw : w\in N^+_D(v)\setminus N^-_D(v)\},
        \]
        and then adding all arcs of
        \[
            \{vu : u\in N^-_D(v) \setminus N^+_D(v)\}.
        \]
        Note that the $d_{D'}^+(x) \leq \Delta$ for all vertices $x\neq v$ because we only added arcs from $v$ to $V(D) \setminus \{v\}$. We further have $d^+_{D'}(v) = d^-_{D}(v) = \Delta-1$, so $\Delta^+(D') \leq \Delta$. Claim~\ref{claim:simple_in_out} implies that $D'$ contains less simple arcs than $D$. By choice of $D$, either $\dic(D')\leq \Delta-1$ or $D'$ contains a copy of $\bid{K}_{\Delta}$ or $\vec{C_3} \boxplus \bid{K}_{\Delta-2}$.

        The former case does not hold, as every dicolouring of $D'$ is a dicolouring of $D$. Indeed, assume that $\phi$ is a dicolouring of $D'$, but $D$ coloured with $\phi$ contains a monochromatic directed cycle $\cSCR$. Clearly, $\cSCR$ contains $v$, for otherwise it is a cycle of $D'$. Let $u$ be the predecessor of $v$ in $\cSCR$, then $\dig{u,v}$ is a monochromatic digon of $D'$, a contradiction.

        Therefore $D'$ contains a copy $K$ of $\bid{K}_{\Delta}$ or $\vec{C_3} \boxplus \bid{K}_{\Delta-2}$. We can further assume that $K$ is an induced copy of $\bid{K}_{\Delta}$ or $\vec{C_3} \boxplus \bid{K}_{\Delta-2}$, as adding any arc to $\vec{C_3} \boxplus \bid{K}_{\Delta-2}$ yields a copy of $\bid{K}_{\Delta}$.

        Clearly, $v\in V(K)$, otherwise $K$ is a subdigraph of $D$. Since $d^-_{D'}(v) = d^+_{D'}(v) = \Delta-1$, and because $N^-_{D'}(v) = N^+_{D'}(v)$, $K$ is necessarily a copy of $\bid{K}_{\Delta}$, and $V(K)\setminus \{v\} = N^-(v)$ induces a biclique on $D$, as desired.
    \end{proofclaim}

    The following will allow us to find our specific vertex $v\in \mathcal{V}_{\Delta-1}$.

    \begin{claimlabeled}[thm:directed_bk_deltaplus]
        There exists $v\in \mathcal{V}_{\Delta-1}$ such that every simple in-neighbour $u$ of $v$ satisfies $d^-(u) \leq 2\Delta-4$.
    \end{claimlabeled}
    \begin{proofclaim}
        Assume that the statement does not hold, that is, every vertex $v\in \mathcal{V}_{\Delta-1}$ has a simple in-neighbour $u$ with in-degree at least $2\Delta-3$. We use a simple discharging argument to reach a contradiction.
        Let $X$ be the set of vertices with out-degree $\Delta-1$, and for every vertex $v$, let its initial charge $\ch(v)$ be exactly $d^-(v)$. Then the total charge in $D$ is:
        \begin{equation}
        \sum_{v\in V(D)} \ch(v) = |A(D)| = \sum_{v\in V(D)} d^+(v) = \Delta \cdot |V(D)| - |X|.
        \label{eq:sum_charge}
        \end{equation}
        We now apply the following rule:
        \begin{equation}
            \textit{For every $v\in \mathcal{V}_{\Delta-1}$ and $u\in N^-(v)\setminus N^+(v)$ with $d^-(u) \geq 2\Delta-3$, $v$ receives 2 from $u$.}
            \label{eq:discharging_rule}
            \tag{R}
        \end{equation}
        Let $\ch^\star(v)$ denote the charge of a vertex $v\in V(D)$ after applying~\eqref{eq:discharging_rule}. We have the following lower bounds on $\ch^\star(v)$.
        \begin{itemize}
            \item If $v\in \mathcal{V}_{\Delta-1}$, then $\ch^\star(v) \geq \Delta+1$.
            
            By assumption, $v$ has a simple neighbour $u$ with in-degree at least $2\Delta-3$. Therefore, $u$ sends $2$ to $v$. Moreover, $v$ does not send any charge as it has in-degree $\Delta-1$.

            \item If $\Delta \leq d^-(v) \leq 2\Delta-4$, then $\ch^\star(v) \geq \Delta$.

            By definition, the initial charge of $v$ is at least $\Delta$, and $v$ does not send any charge.

            \item If $d^-(v) \geq 2\Delta-3$, then $\ch^\star(v) \geq \Delta$.

            This is clear if $v$ does not send any charge. If it does send charge, then it has a simple out-neighbour $w\in \mathcal{V}_{\Delta-1}$.
            By Claim~\ref{claim:in-neighbourhood_biclique}, $v$ belongs to a biclique of size $\Delta-1$. Hence $v$ has at most two simple out-neighbours, so it sends charge at most twice. Its final charge is then at least $2\Delta-3 - 4 \geq \Delta$.

            \item If none of the cases above applies, then $d^-(v) = d^+(v)=\Delta-1$. In particular, $v\in X$ and $\ch^\star(v) = \Delta-1$.
        \end{itemize}
        Therefore, the total charge in $D$ is 
        \[ 
        \sum_{v\in V(D)} \ch^\star(v) \geq \Delta \cdot |V(D)| - |X| + |\mathcal{V}_{\Delta-1}| > \Delta \cdot |V(D)| - |X|,
        \]
        where in the last inequality we used that $\mathcal{V}_{\Delta-1}$ is non-empty by Claim~\ref{claim:V_deltam1_nonempty}. By~\eqref{eq:sum_charge}, we thus have $\sum_{v} \ch^\star(v) > \sum_{v} \ch(v)$, a contradiction.
    \end{proofclaim}
    From now on, we thus fix a vertex $v\in \mathcal{V}_{\Delta-1}$ such that each of its simple in-neighbours has in-degree at most $2\Delta-4$. We further fix such a vertex $u\in N^-(v) \setminus N^+(v)$, we denote $X=N^-[v]$, and we arbitrarily label $x_1,\dots,x_{\Delta-2}$ the vertices of $X\setminus \{u,v\}$. 
    
    In what remains, using arguments similar to the ones of Section~\ref{sec:directed_bk:structure}, especially Lemma~\ref{lemma:no_KD_minus_uv}, we justify that a $(\Delta-1)$-dicolouring of $D-X$ can be extended to $D$, thus yielding the contradiction. However, the fact that vertices may have arbitrarily large in-degree requires a more technical analysis. We start with the following analogue of Claim~\ref{claim:no_KD_minus_uv:B}.

    \begin{claimlabeled}[thm:directed_bk_deltaplus]
        \label{claim:wi_delplus}
        For every $i\in [\Delta-2]$, there exists a unique vertex $w_i\in N^+(x_i) \setminus X$. Moreover, for every $(\Delta-1)$-dicolouring $\phi$ of $D-X$, 
        \begin{enumerate}[label=$(\roman*)$]
            \item for every $i\in [\Delta-2]$, $D-X$ coloured with $\phi$ contains a monochromatic path from $w_i$ to $N^-(x_i)$; and
            \label{enum:claim_wi_delplus:1}
            \item for all $i,j\in [\Delta-2]$, $\phi(w_i) =\phi(w_j)$.
            \label{enum:claim_wi_delplus:2}
        \end{enumerate}
    \end{claimlabeled}
    \begin{proofclaim}
        Let $\phi$ be a $(\Delta-1)$-dicolouring of $D-X$. We first extend $\phi$ to $D-(X\setminus \{u,v\})$ by choosing for $u$ and $v$ a colour that is not appearing in $N^-(u) \cup N^-(v)$. This is possible because $uv$ is simple, no in-neighbour of $v$ is coloured, and $u$ has at most $(2\Delta-4) - (\Delta-2) = \Delta-2$ in-neighbours in $V(D)\setminus X$. 

        Now, one colour (namely $\phi(u)=\phi(v)$) is repeated in the out-neighbourhood of every vertex $x\in X\setminus \{u,v\}$. Then, for every $i$, $x_i$ has out-degree exactly $\Delta$, and admits an out-neighbour $w_i$ and an in-neighbour $y_i$ in $V(D)\setminus X$ such that $D-X$, coloured with $\phi$, contains a monochromatic path $P_i$ from $w_i$ to $y_i$ (we may have $w_i = y_i$). Indeed, if this is not the case, we can greedily colour the vertices of $X\setminus \{u,v\}$, keeping $x_i$ for the end, by choosing for every vertex a colour that is not appearing in its out-neighbourhood. The uniqueness of $w_i$ comes from the fact that $x_i$ has exactly $\Delta-1$ out-neighbours in $X$. This shows~\ref{enum:claim_wi_delplus:1}.

        To show~\ref{enum:claim_wi_delplus:2}, assume for a contradiction that, for some $i,j$, $\phi(w_i) \neq \phi(w_j)$. If $\phi(w_j) = \phi(u)$, then we can greedily extend $\phi$, keeping $x_j$ for the end, as one colour (namely $\phi(u)$) appears three times in $N^+(x_j)$.
        Hence assume that $\phi(w_j)\neq \phi(w_i)$ and $\phi(w_j) \neq \phi(u)$. We can thus extend $\phi$ by first colouring $w_i$ with $\phi(w_j)$, and then keeping $x_j$ for the end, as two colours are now repeated in $N^+(x_j)$.
    \end{proofclaim}
    
    From now on, we will consider two distinct partitions $(I^{\rm dig}, I^{\rm s})$ and $(I^{\pm},I^{\leq 6 },I^{>6})$ of $[\Delta-2]$. We define $(I^{\rm dig}, I^{\rm s})$ as follows:
    \[
        I^{\rm dig} = \Big\{i\in [\Delta-2] : \dig{x_i,w_i} \text{~is a digon of $D$}\Big\}
        \text{~~~and~~~}
        I^{\rm s} = \Big\{i\in [\Delta-2] : w_i \in N^+(x_i) \setminus N^-(x_i)\Big\},
    \]
    and $(I^{\pm},I^{\leq 6 },I^{>6})$ as follows:
    \begin{flalign*}
        I^{\pm} &= \{i\in [\Delta-2] : N^+(w_i) \cap X \neq \emptyset\},&\\
        I^{\leq 6 } &= \{i\in [\Delta-2] : w_i \text{ has at most six simple out-neighbours}\} \setminus I^{\pm}, \text{~~~and}&\\
        I^{> 6 } &= \{i\in [\Delta-2] : w_i \text{ has at least seven simple out-neighbours}\} \setminus I^{\pm}.&
    \end{flalign*}
    By definition, observe that $I^{\rm dig} \subseteq I^{\pm}$ and $(I^{\leq 6 } \cup I^{>6}) \subseteq I^{\rm s}$. Also, note that one set of $(I^{\pm},I^{\leq 6 },I^{>6})$ is large, as it has size at least $(\Delta-2)/3$. Using Claim~\ref{claim:wi_delplus}, our goal is now to show that each of the sets $(I^{\pm},I^{\leq 6 },I^{>6})$ has bounded size, thus yielding a contradiction.
    We first exclude the case of some vertex $w$ being adjacent with digons to all vertices of $X\setminus \{u,v\}$.
    
    \begin{claimlabeled}[thm:directed_bk_deltaplus]
        \label{claim:exclude_alldigons_delplus}
        For every $w\in V(D)\setminus X$, there exists $x\in X\setminus \{u,v\}$ such that $\dig{x,w}\nsubseteq A(D)$. 
    \end{claimlabeled}
    \begin{proofclaim}
        Assume for a contradiction that $(X\setminus \{u,v\}) \cup \{w\}$ is a biclique, for some $w\in V(D)\setminus X$. Recall that $X\setminus \{v\}$ is also a biclique. We now distinguish two cases, depending on whether $u$ is the unique simple in-neighbour of $v$ or not, see Figure~\ref{fig:exclude_alldigons_delplus} for an illustration.
        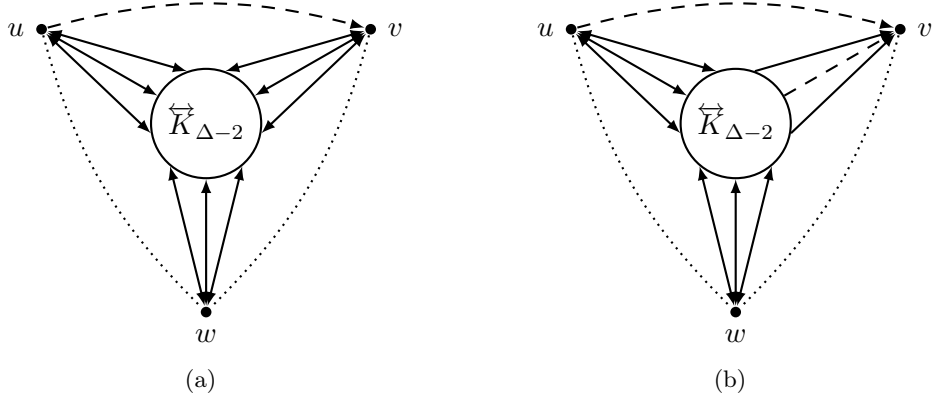
\begin{figure}[hbtp]
            \begin{center}
            \subfigure[]{
            \begin{tikzpicture}[thick,scale=1, every node/.style={transform shape}]
                \node[bigvertex] (K) at (0,0) {$\bid{K}_{\Delta-2}$};
                \node[vertex, label=left:$u$] (u) at (150:2.5) {};
                \node[vertex, label=right:$v$] (v) at (30:2.5) {};
                \node[vertex, label=below:$w$] (w) at (-90:2.5) {};
    
                \draw[edge] (u) to[out=15, in=165] (v);
                \draw[dotted] (u) to[out=-75, in=135] (w);
                \draw[dotted] (v) to[out=-105, in=45] (w);
                \draw[digon] (u) to (K.190);
                \draw[digon] (u) to (K.150);
                \draw[digon] (u) to (K.110);
                \draw[digon] (w) to (K.-130);
                \draw[digon] (w) to (K.-90);
                \draw[digon] (w) to (K.-50);
                \draw[digon] (v) to (K.30);
                \draw[digon] (v) to (K.70);
                \draw[digon] (v) to (K.-10);
              \end{tikzpicture}
            }
            \hspace*{1cm}
            \subfigure[]{
                \begin{tikzpicture}[thick,scale=1, every node/.style={transform shape}]
                \node[bigvertex] (K) at (0,0) {$\bid{K}_{\Delta-2}$};
                \node[vertex, label=left:$u$] (u) at (150:2.5) {};
                \node[vertex, label=right:$v$] (v) at (30:2.5) {};
                \node[vertex, label=below:$w$] (w) at (-90:2.5) {};
    
                \draw[edge] (u) to[out=15, in=165] (v);
                \draw[dotted] (u) to[out=-75, in=135] (w);
                \draw[dotted] (v) to[out=-105, in=45] (w);
                \draw[digon] (u) to (K.190);
                \draw[digon] (u) to (K.150);
                \draw[digon] (u) to (K.110);
                \draw[digon] (w) to (K.-130);
                \draw[digon] (w) to (K.-90);
                \draw[digon] (w) to (K.-50);
                \draw[arc] (K.70) to (v);
                \draw[edge] (K.30) to (v);
                \draw[arc] (K.-10) to (v);
            \end{tikzpicture}
            }
              \caption{The two distinguished cases in the proof of Claim~\ref{claim:exclude_alldigons_delplus}. When two vertices are linked with a dotted edge, it illustrates that there is at most one arc between them. A dashed arc is a simple arc.}
              \label{fig:exclude_alldigons_delplus}
            \end{center}
        \end{figure}
        
        Assume first that $u$ is the unique simple in-neighbour of $v$, see Figure~\ref{fig:exclude_alldigons_delplus}(a).
        In this case, $D[\{u,v,w\}]$ is acyclic, for otherwise $D$ contains either $\bid{K}_{\Delta}$ or $\vec{C_3}\boxplus \bid{K}_{\Delta-2}$.
        Let $\phi$ be a $(\Delta-1)$-dicolouring of $D-(X\cup \{w\})$. Since $u,v,w$ all have $\Delta-2$ uncoloured out-neighbours, we can extend $\phi$ by first giving to $u,v,w$ the same colour. We can then extend $\phi$ to $D$ by choosing for each vertex $x\in X\setminus \{u,v\}$ a colour that is not appearing in $N^+(x)$, which is possible because one colour (namely $\phi(v)$) is appearing three times, a contradiction.
    
        Assume then that $v$ has a simple in-neighbour $u'$ distinct from $u$, see Figure~\ref{fig:exclude_alldigons_delplus}(b).
        Again, let $\phi$ be a $(\Delta-1)$-dicolouring of $D-(X\cup \{w\})$. We now extend $\phi$ by choosing a common colour for both $u'$ and $v$ that is not appearing in $N^-(u')\cup N^-(v)$. This is possible because no in-neighbour of $v$ is coloured, and by choice of $v$, $d^-(u')\leq 2\Delta-4$ so $|N^-(u') \setminus (X\cup \{w\})|\leq \Delta-3$.
        We then extend $\phi$ by giving a common colour to $u$ and $w$. Again, this is possible as $\dig{u,w}$ is not a digon (otherwise $D$ contains $\bid{K}_\Delta$), and both $u$ and $w$ have $\Delta-3$ uncoloured out-neighbours.
        We can finally extend $\phi$ to $D$ by choosing for each vertex $x\in X\setminus \{u,u',v\}$ a colour that is not appearing in $N^+(x)$, which is possible because two colours are repeated in $N^+(x)$, a contradiction.
    \end{proofclaim}

    We now need a few extra definitions. 
    For every $i\in I^{\rm s}$, we let $N_i$ be $N^-(x_i) \setminus X$. We further let $D_i$ be the digraph obtained from $D-X$ by adding all arcs from $N_i$ to $w_i$ (note that $w_i\notin N_i$ by definition of $I^{\rm s}$).

    \begin{claimlabeled}[thm:directed_bk_deltaplus]
        \label{claim:delplus:Di_contains_K}
        For every $i\in I^{\rm s}$, $D_i$ contains a copy of $\bid{K}_{\Delta}$ or $\vec{C}_3 \boxplus \bid{K}_{\Delta-2}$.
    \end{claimlabeled}
    \begin{proofclaim}
        Assume that this is not the case. By construction, $\delplus(D_i)\leq \Delta$, so the minimality of $D$ implies the existence of a $(\Delta-1)$-dicolouring $\phi$ of $D_i$. 
        Then, by Claim~\ref{claim:wi_delplus}, since $\phi$ is a particular dicolouring of $D-X$, there exists a monochromatic path $P_i$ from $w_i$ to some $y_i \in N_i$. 
        The concatenation of $P_i$ and the arc $y_iw_i$ is a monochromatic directed cycle of $D_i$ coloured with $\phi$, a contradiction.
    \end{proofclaim}

    Claim~\ref{claim:delplus:Di_contains_K} allows us to do the following. For every $i\in I^{\rm s}$, we fix an induced copy $K_i$ of either $\bid{K}_{\Delta}$ or $\vec{C}_3 \boxplus \bid{K}_{\Delta-2}$ in $D_i$, and we fix $\Tilde{N}_i$ to:
    \[
        \Tilde{N}_i = (V(K_i)\cap N^-(x_i))\setminus N^-(w_i).
    \]
    The key point is that $\Tilde{N}_i$ can play the role of $N_i$, that is $\Tilde{N}_i \subseteq N^-(x_i)\setminus X$ and 
    in every $(\Delta-1)$-dicolouring of $D-X$, there is a monochromatic path from $w_i$ to $\Tilde{N}_i$, but we further have the extra property that $\Tilde{N}_i$ is almost the set of simple out-neighbours of $w_i$. We show these properties in the next two claims.
    
    \begin{claimlabeled}[thm:directed_bk_deltaplus]
        \label{claim:colouring_Ni_tilde}
        For every $(\Delta-1)$-dicolouring $\phi$ of $D-X$ and every $i\in I^{\rm s}$, there is a monochromatic path from $w_i$ to $\Tilde{N}_i$.
    \end{claimlabeled}
    \begin{proofclaim}
        Let $\phi_i$ be the restriction of $\phi$ to $V(K_i)$. Clearly, $\phi_i$ is not a dicolouring of $K_i$, as $\dic(K_i) = \Delta$. Let $\cSCR$ be a monochromatic directed cycle in $K_i$ coloured with $\phi_i$. Since the arcs of $A(K_i)\setminus A(D)$ all have the same head (namely $w_i$), and because $\cSCR$ is not a cycle of $D$, there exists a unique arc $xy$ of $\cSCR$ in $A(K_i)\setminus A(D)$. Since $xy\in A(D_i)\setminus A(D)$, we have $y=w_i$, and $x\in N^-(x_i)\setminus N^-(w_i)$. The claim follows.
    \end{proofclaim}

    \begin{claimlabeled}[thm:directed_bk_deltaplus]
        \label{claim:delplus:Ni_closeto_Si}
        Let $i\in I^{\rm s}$, and $S_i^+ = N^+(w_i)\setminus N^-(w_i)$ be the set of simple out-neighbours  of $w_i$. Then
        \[
            |S_i^+ \setminus \Tilde{N}_i|\leq 2 \text{~~~and~~~}  |\Tilde{N}_i\setminus S_i^+|\leq 1.
        \]
        Moreover, $\Tilde{N}_i\cap S_i^+ \neq \emptyset$.
    \end{claimlabeled}
    \begin{proofclaim}
        By definition, $K_i$ is obtained from $D[V(K_i)]$ by adding all arcs from $\Tilde{N}_i$ to $w_i$. Observe also that, in $K_i$, $w_i$ is linked with digons to at least $\Delta-2$ vertices. Let us call $Z_i = N^\pm_{K_i}(w_i)$ the set of these vertices. Suppose for a contradiction that we can find three simple out-neighbours of $w_i$ outside $\tilde N_i$. Note that all of them lie outside $Z_i$, regardless of whether they are in $V(K_i)$ or not. Therefore, $w_i$ has at least $3+(\Delta-2)=\Delta+1$ out-neighbours in $D_i$, a contradiction, because all of them are also out-neighbours in $D$.
        The first inequality follows.

        For the second inequality, recall that all vertices in $\Tilde{N}_i$ are in-neighbours of $w_i$ in $K_i$ because $K_i$ is an induced copy in $D_i$, and vertices in $\Tilde{N}_i$ belong to $N^-(x_i)$. However, they are not in-neighbours of $w_i$ in $D$ by definition of $\tilde{N}_i$. Since, in $K_i$, $w_i$ has at most one simple in-neighbour, it means that all vertices in $\Tilde{N}_i$, except possibly one, are simple out-neighbours of $w_i$ in $D$. The second inequality follows.

        Assume now, for a contradiction, that $\Tilde{N}_i$ contains no element of $S_i^+$. By the second inequality, $|\Tilde{N}_i| = 1$ ($\tilde{N}_i$ is not empty otherwise $K_i$ would be a subdigraph of $D$). Let us denote by $y$ the unique vertex in $\Tilde{N}_i$. Recall that $K_i\subseteq D\cup \{yw_i\}$, and $K_i$ is either $\bid{K}_\Delta$ or $\vec{C}_3 \boxplus \bid{K}_{\Delta-2}$. Since $y$ is not a simple out-neighbour of $w_i$ in $D$, then $yw_i$ is a simple arc in $K_i$, which implies that $K_i$ is exactly $\vec{C}_3 \boxplus \bid{K}_{\Delta-2}$, and that $D[V(K_i)]$ is exactly $\vec{P}_3\boxplus \bid{K}_{\Delta-2}$. This is a contradiction to Claim~\ref{claim:delplus:no_P3_Kdel}.
    \end{proofclaim}
    
    The last important property is that the sets $(\tilde{N}_i)_{i\in I^{\rm s}}$ have limited overlap, as we show now.
    
    \begin{claimlabeled}[thm:directed_bk_deltaplus]
    \label{claim:delplus:Ni_almost_disjoint}
        For all distinct $i,j,k\in I^{\rm s}$, we have $\tilde{N}_i \cap \tilde{N}_j \cap \tilde{N}_k = \emptyset$.
    \end{claimlabeled}
    \begin{proofclaim}
        Assume for a contradiction that $y \in \tilde{N}_i \cap \tilde{N}_j \cap \tilde{N}_k$.
        This implies in particular $|N^+(y) \cap X| \geq 3$ (as $\{x_i,x_j,x_k\} \subseteq N^+(y)$).
        
        We further know that $y$, in $K_i$, has at least $\Delta-1$ out-neighbours. 
        By definition of $K_i$, all of these except $w_i$ are out-neighbours of $y$ in $D$; we thus have $|N^+_D(y) \cap V(K_i)|\geq \Delta-2$. Since $V(K_i)\cap X=  \emptyset$, we get $d^+(y) \geq \Delta+1$, a contradiction.
    \end{proofclaim}

    We now have enough material to show that each of the sets $(I^{\pm},I^{\leq 6 },I^{>6})$ has bounded size. We start with the following, which intuitively implies that if $I\subseteq [\Delta-2]$ is large, then so is $W_I = \{w_i : i\in I\}$.

    \begin{claimlabeled}[thm:directed_bk_deltaplus]
        \label{claim:delplus:the_wis_are_almost_distinct}
        For every $w\in V(D)\setminus X$, $|N^-(w)\cap \{x_1,\ldots,x_{\Delta-2}\}|\leq 14$.
    \end{claimlabeled}
    \begin{proofclaim}
        Let $I_w$ be the set of integers $i\in [\Delta-2]$ such that $w_i = w$, and assume for a contradiction that $|I_w|\geq 15$.

        We first assume that $|I_w \cap I^{\rm dig}| \geq 3$, so $w$ is linked with digons to at least three vertices in $\{x_1,\dots,x_{\Delta-2}\}$. Let us fix $j\in I_w \cap I^{\rm dig}$. By Claim~\ref{claim:exclude_alldigons_delplus}, there exists $i\in [\Delta-2]$ such that $\dig{w,x_i}$ is not a digon of $D$. We claim that $w$ is not an out-neighbour of $x_i$. Indeed, if it is the case, then $w=w_i$ is a simple out-neighbour of $x_i$, and $w$ has at least $\Delta-1$ out-neighbours in $V(K_i) \subseteq V(D)\setminus X$. Since $w$ has three out-neighbours in $X$, this is a contradiction.

        Therefore, $w\neq w_i$. Let $\phi$ be a $(\Delta-1)$-dicolouring of $D-(X\cup \{w\})$. Since $w$ has three out-neighbours in $X$, we can extend $\phi$ to $D-X$ by choosing for $w$ a colour that is not appearing in $N^+(w)$ and that is further distinct from $\phi(w_i)$. We obtain a $(\Delta-1)$-dicolouring $\phi$ of $D-X$ with $\phi(w_i)\neq \phi(w_j)$, a contradiction to Claim~\ref{claim:wi_delplus}\ref{enum:claim_wi_delplus:2}.

        \medskip

        Henceforth, we assume that $I_w \cap I^{\rm s}$ has size at least $15-2 = 13$. Let $S_w^+$ be the set of simple out-neighbours of $w$. We conclude by distinguishing two cases, depending on the size of $S_w^+$. In both cases, we get a contradiction by showing that some $y\in S_w^+$ has out-degree at least $\Delta+1$.

        Assume first that $|S_w^+| \geq 7$, and consider  three distinct indices $i,j,k \in I_w \cap I^{\rm s}$. By the first inequality of Claim~\ref{claim:delplus:Ni_closeto_Si}, we have
        \[
            |S_w^+ \cap \Tilde{N}_i \cap \Tilde{N}_j \cap \Tilde{N}_k| \geq |S_w^+| - |(S_w^+ \setminus \Tilde{N}_i) \cup (S_w^+ \setminus \Tilde{N}_j) \cup (S_w^+ \setminus \Tilde{N}_k)| \geq 7-6=1.
        \]
        We thus have a simple out-neighbour $y$ of $w$ which belongs to $\Tilde{N}_i \cap \Tilde{N}_j \cap \Tilde{N}_k$. This is a contradiction to Claim~\ref{claim:delplus:Ni_almost_disjoint}.
        
        Assume then that $|S_w^+|\leq 6$. Recall that, for every $i\in I_w \cap I^{\rm s}$, $\Tilde{N}_i$ contains a simple out-neighbour of $w$ by the last part of Claim~\ref{claim:delplus:Ni_closeto_Si}. Hence, since $|I_w \cap I^{\rm s}| \geq 13$, and because $w$ has at most six simple out-neighbours, by the pigeonhole principle, some simple out-neighbour $y$ of $w$ belongs to $\Tilde{N}_i \cap \Tilde{N}_j \cap \Tilde{N}_k$ for three distinct $i,j,k\in I_w\cap I^{\rm s}$. Again, this is a contradiction to Claim~\ref{claim:delplus:Ni_almost_disjoint}.
    \end{proofclaim}

    In the next three claims, we finally show that each of the sets $(I^{\pm},I^{\leq 6 },I^{>6})$ has bounded size. We start with $I^{\pm}$.

    \begin{claimlabeled}[thm:directed_bk_deltaplus]
        \label{claim:delplus:Idig_bounded}
        $|I^{\pm}| \leq 70$.
    \end{claimlabeled}
    \begin{proofclaim}
        Assume for a contradiction that $|I^{\pm}| \geq 71$. By Claim~\ref{claim:delplus:the_wis_are_almost_distinct}, there exists $I \subseteq I^{\pm}$ such that $|I| = \lceil 71/14\rceil = 6$ and vertices $(w_i)_{i \in I}$ are pairwise distinct.
        For every distinct $i,j\in I$, we define $D_{ij}$ as the digraph obtained from $D-X$ by further adding the digon $\dig{w_i,w_j}$, see Figure~\ref{fig:delplus:construction_Dij}(a) for an illustration.
        
        Since $i,j$ belong to $I^{\pm}$, we have $\delplus(D_{ij})\leq \Delta$. 
        Furthermore, $\dic(D_{ij})\geq \Delta$, as any $(\Delta-1)$-dicolouring $\phi$ of $D_{ij}$ would yield a $(\Delta-1)$-dicolouring of $D-X$ satisfying $\phi(w_i) \neq \phi(w_j)$, a contradiction to Claim~\ref{claim:wi_delplus}.

        Therefore, the minimality of $D$ implies that each $D_{ij}$ contains an induced copy of either $\bid{K}_{\Delta}$ or $\vec{C_3}\boxplus \bid{K}_{\Delta-2}$. In both cases, we conclude on the existence of a biclique $\lSCR_{ij}$ in $D-X$ of size $\Delta-4$ such that 
        \[
        \lSCR_{ij} \subseteq N^{\pm}(w_i) \cap N^{\pm}(w_j).
        \]
        Let us fix two distinct $i,j\in I$.
        Observe that, for any $k \in I\setminus \{i,j\}$, $|\lSCR_{ij} \setminus \lSCR_{ik}| \leq 3$, for otherwise $w_i$ has at least $\Delta$ out-neighbours in $V(D)\setminus X$ and one out-neighbour in $X$, a contradiction. Therefore, 
        \[
        \left| \bigcap_{k\in I\setminus \{i\}} \lSCR_{ik} \right| \geq |\lSCR_{ij}| - \sum_{k\in I\setminus \{i,j\}} |\lSCR_{ij}\setminus \lSCR_{ik}| \geq (\Delta-4) - 3(|I|-2) > 0.
        \]
        Let us thus fix $y\in \bigcap_{k\in I\setminus \{i\}}\lSCR_{ik}$. Then $y$ is an in-neighbour of each $w_k$, $k\in I$. Furthermore, $y$ has $\Delta-5$ out-neighbours in $\lSCR_{ij}$. Note that $\{w_k : k\in I\}$ is disjoint from $\lSCR_{ij}$, as if $w_k\in \lSCR_{ij}$ then $\dig{w_i,w_k}$ is a digon of $D-X$, and $D_{ik}$ is a subdigraph of $D$, a contradiction.
        Hence, $d^+(y) \geq |I| + \Delta-5 = \Delta+1$, a contradiction.
    \end{proofclaim}

    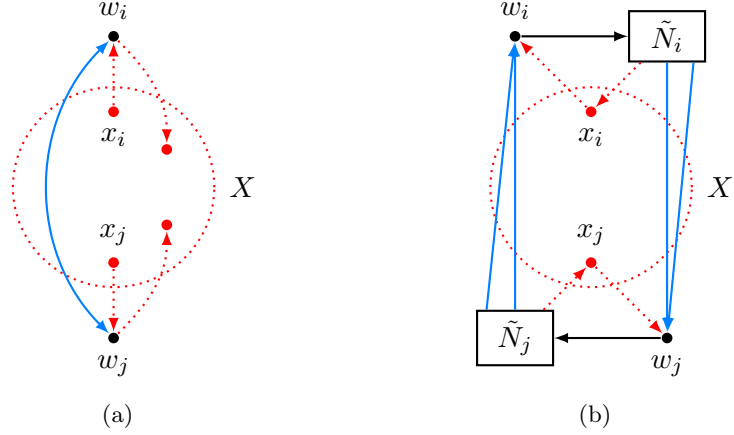
\begin{figure}
        \begin{center}
        \subfigure[]{
        \begin{tikzpicture}[thick,scale=1, every node/.style={transform shape}]
            \node[bigvertex, dotted, red, minimum size=16ex] (K) at (0,0) {};
            \node[vertex,red, label=below:${x_i}$] (xi) at (0,1) {};
            \node[vertex,red, label=above:${x_j}$] (xj) at (0,-1) {};
            \node[vertex, label=above:${w_i}$] (wi) at (0,2) {};
            \node[vertex, label=below:${w_j}$] (wj) at (0,-2) {};
            \node[] (X) at (1.7,0) {$X$};
            \node[vertex,white] (Y) at (-1.7,0) {};
            \draw[digon, g-blue] (wi) to[out=-135, in=135] (wj);
            \draw[edge, dotted, red] (xi) to (wi);
            \node[vertex, red] (xibis) at (0.7,0.5) {};
            \draw[edge, dotted, red] (wi) to[out=-45, in=90] (xibis);
            \draw[edge, dotted, red] (xj) to (wj);
            \node[vertex, red] (xjbis) at (0.7,-0.5) {};
            \draw[edge, dotted, red] (wj) to[out=45, in=-90] (xjbis);
          \end{tikzpicture}
        }
        \hspace*{2cm}
        \subfigure[]{
            \begin{tikzpicture}[thick,scale=1, every node/.style={transform shape}]
            \node[bigvertex, dotted, red, minimum size=16ex] (K) at (0,0) {};
            \node[vertex,red, label=below:${x_i}$] (xi) at (0,1) {};
            \node[vertex,red, label=above:${x_j}$] (xj) at (0,-1) {};
            \node[] (X) at (1.7,0) {$X$};
            \node[vertex,white] (Y) at (-1.7,0) {};
            
            \node[vertex, label=above:${w_i}$] (wi) at (-1,2) {};
            \node[vertex, label=below:${w_j}$] (wj) at (1,-2) {};
            \node[draw,minimum width=1cm,minimum height=0.6cm] (Ni) at (1,2)  {$\tilde{N}_i$};
            \node[draw,minimum width=1cm,minimum height=0.6cm] (Nj) at (-1,-2)  {$\tilde{N}_j$};

            \draw[edge, dotted, red] (xi) to (wi);
            \draw[edge, solid] (wi) to (Ni);
            \draw[edge, dotted, red] (Ni) to (xi);
            
            \draw[edge, dotted, red] (xj) to (wj);
            \draw[edge, solid] (wj) to (Nj);
            \draw[edge, dotted, red] (Nj) to (xj);
            
            \draw[edge, solid, g-blue] (Ni) to (wj);
            \draw[edge, solid, g-blue] (Ni.-45) to (wj);
    
            \draw[edge, solid, g-blue] (Nj) to (wi);
            \draw[edge, solid, g-blue] (Nj.135) to (wi);
        \end{tikzpicture}
        }
          \caption{The construction of $D_{ij}$ when $i,j\in I^{\pm}$ (Figure~(a)), and $i,j\in I^{\leq 6}\cup I^{> 6}$ (Figure~(b)). In both figures, $D_{ij}$ is built from $D$ by removing the dotted red part and adding the solid blue one.
          In Figure~(b), $w_i$ does not have any out-neighbour in $X$ nor any in-neighbour in $\tilde{N}_i$. Moreover, all but at most one vertex of $\tilde{N}_i$ are out-neighbours of $w_i$. The same properties hold for $w_j$ and $\tilde{N}_j$.
          }
          \label{fig:delplus:construction_Dij}
        \end{center}
    \end{figure}
    
    We now move to the cases of $I^{\leq 6}$ and $I^{>6}$. For two integers $i,j \in I^{\leq 6}\cup I^{>6}$ such that $w_i \neq w_j$, we define $D_{ij}$ as the digraph obtained from $D-X$ by adding all arcs from $\Tilde{N}_i$ to $w_j$ and all arcs from $\Tilde{N}_j$ to $w_i$. Observe that $w_j \notin \Tilde{N}_i$, for otherwise $w_j$ would have an out-neighbour in $X$ (namely $x_i$), which by definition would imply $j\in I^{\pm}$. See Figure~\ref{fig:delplus:construction_Dij}(b) for an illustration.

    For all $i,j\in I^{\leq 6}\cup I^{>6}$ with $w_i\neq w_j$, observe that $D_{ij}$ satisfies $\delplus(D_{ij})\leq \Delta$ by construction (this is because $\Tilde{N}_i \subseteq N^-(x_i)$ and $\Tilde{N}_j \subseteq N^-(x_j)$). Moreover, $\dic(D_{ij}) \geq \Delta$. To see this, assume for a contradiction that $D_{ij}$ admits a $(\Delta-1)$-dicolouring $\phi$. Then $\phi$ is a particular $(\Delta-1)$-dicolouring of $D-X$, so by Claim~\ref{claim:wi_delplus} we have $\phi(w_i) = \phi(w_j)$. Moreover, by Claim~\ref{claim:colouring_Ni_tilde} (recall that $I^{\leq 6}\cup I^{>6}\subseteq I^{\rm s}$), $D_{ij}$, coloured with $\phi$, contains a monochromatic directed path $P_i$ from $w_i$ to some $y_i \in \Tilde{N}_i$ and one $P_j$ from $w_j$ to $y_j\in \Tilde{N}_j$. Hence the concatenation of $P_i$, $y_iw_j$, $P_j$, and $y_jw_i$ is a monochromatic directed cycle in $D_{ij}$ coloured with $\phi$, a contradiction. 
    
    The choice of $D$ thus implies that, for all $i,j\in I^{\leq 6}\cup I^{>6}$ with $w_i\neq w_j$, $D_{ij}$ contains an induced copy $K_{ij}$ of either $\bid{K}_{\Delta}$ or $\vec{C_3}\boxplus \bid{K}_{\Delta-2}$.

    We now bound the size of $I^{\geq 6}$. The proof is analogous to that of Claim~\ref{claim:delplus:Idig_bounded}, but requires more technical arguments. Intuitively, for integers $i,j \in I^{\leq 6}$, the digraph $D_{ij}$ is obtained from $D-X$ by adding only a constant number of arcs. Using this, we will find a large biclique in $D-X$ that is linked with digons to many vertices in $\bigcup_{i\in I^{\geq 6}} \Tilde{N}_i$, hence implying that some vertex in this biclique has too large out-degree.
    
    \begin{claimlabeled}[thm:directed_bk_deltaplus]
        \label{claim:delplus:Ileq6_bounded}
        $|I^{\leq 6}| \leq 441$.
    \end{claimlabeled}
    \begin{proofclaim}
        Assume for a contradiction that $|I^{\leq 6}| \geq 442$.
        We claim that there exists $I\subseteq I^{\leq 6}$, $|I| = 22$, such that:
        \begin{enumerate}
            \item vertices $(w_i)_{i\in I}$ are pairwise distinct,
            \item sets $(\Tilde{N}_i)_{i\in I}$ are pairwise disjoint.
        \end{enumerate}

        To see the existence of $I$, we build an auxiliary undirected graph $G^{\leq 6}$ with vertex-set $V(G^{\leq 6}) = I^{\leq 6}$ and edge-set $E(G^{\leq 6}) = E_1 \cup E_2$, where
        \[
            E_1 = \{\{i,j\}  : w_i = w_j\}
            \text{~~~and~~~} E_2 = \{\{i,j\} : \Tilde{N}_i \cap \Tilde{N}_j \neq \emptyset \}.
        \]
        Our goal is thus to show that $G^{\leq 6}$ admits an independent set $I$ of size $22$. For this, fix an integer $i\in I^{\leq 6}$, which is also a vertex of $G^{\leq 6}$.
        By Claim~\ref{claim:delplus:the_wis_are_almost_distinct}, $i$ has at most 13 neighbours in $G^{\leq 6}$ via edges in $E_1$.
        Similarly, by Claim~\ref{claim:delplus:Ni_almost_disjoint}, each vertex in $\Tilde{N}_i$ belongs to $\Tilde{N}_j$ for at most one $j\neq i$. 
        Furthermore, note that $|\Tilde{N}_i|\leq 7$ by Claim~\ref{claim:delplus:Ni_closeto_Si} together with the fact that $w_i$ has at most $6$ simple out-neighbours. Hence, $i$ has at most $7$ neighbours in $G^{\leq 6}$ via edges in $E_2$. 
        This shows $\Delta(G^{\leq 6}) \leq 20$. Hence, $G^{\leq 6}$ has chromatic number at most $\chi(G^{\leq 6}) \leq 21$. Take an optimal proper colouring of $G^{\leq 6}$, its largest colour class is an independent set of $G^{\leq 6}$ of size at least $\lceil 442/21\rceil = 22$. This shows the existence of $I$.

        Recall that, for every distinct $i,j\in I$, $D_{ij}$ contains an induced copy $K_{ij}$ of either $\bid{K}_{\Delta}$ or $\vec{C_3}\boxplus \bid{K}_{\Delta-2}$. We fix such a $K_{ij}$ for every pair $i,j\in I$. As $K_{ij}$ is not a subgraph of $D$, it contains an arc from $\Tilde{N}_i$ to $w_j$ or an arc from $\Tilde{N}_j$ to $w_i$ that is not an arc of $D$ (it might contain both).

        We now build an auxiliary tournament $T^I$ as follows. We take $V(T^I) = I$, and we orient each edge $ij$ from $i$ to $j$ if $K_{ij}$ contains an arc from $\Tilde{N}_i$ to $w_j$ that is not an arc of $D$, and from $j$ to $i$ otherwise (if $K_{ij}$ contains both kind of arcs, we orient $ij$ arbitrarily).
        As $T^I$ has order $22$, a simple counting argument implies that it contains a vertex $i$ of in-degree at least $11$. We fix such a vertex $i$ and a set $J \subseteq I\setminus \{i\}$ of exactly $11$ in-neighbours of $i$ in $T^I$. Hence, by definition, for every $j\in J$, $K_{ij}$ contains an arc from some $y_j \in \Tilde{N}_j$ to $w_i$ that is not an arc of $D$. We fix such a $y_j$ for every $j\in J$. 

        For every $j\in J$, recall that $K_{ij}$ contains a biclique $\lSCR_{ij}'$ of size $\Delta-4$, whose vertices are all linked with digons to both $y_j$ and $w_i$.
        Note that the digons of $D_{ij}[\lSCR_{ij}' \cup \{y_j,w_i\}]$ do not necessarily all belong to $D$, as some of them may contain arcs from $\Tilde{N}_j$ to $w_i$ or arcs from $\Tilde{N}_i$ to $w_j$. In fact, $w_j \notin \lSCR_{ij}'$, for otherwise $w_j$ and $w_i$ are linked by a digon in $D_{ij}$, and so they are in $D$ since $w_j \notin \Tilde{N}_i$ and $w_i \notin \Tilde{N}_j$ ($w_j \in \Tilde{N}_i$ would imply $j\in I^{\pm}$ as $\Tilde{N}_i\subseteq N^-(x_i)$). Again, $w_i$ and $w_j$ cannot be linked with a digon in $D$ by Claim~\ref{claim:wi_delplus}.

        Therefore, in $D$, $\lSCR_{ij} = \lSCR_{ij}' \setminus \Tilde{N}_j$ is a biclique of $D$ of size at least $\Delta - 10$ (as $|\Tilde{N}_j|\leq 7$ and $y_j \notin  \lSCR_{ij}'$), and all of them are linked with digons to both $w_i$ and $y_j$.
        From now on, let us fix $j\in J$.
        Observe that, for any $k \in J\setminus \{j\}$, $|\lSCR_{ij} \setminus \lSCR_{ik}| \leq 10$, for otherwise $w_i$ has at least $\Delta+1$ out-neighbours in $\lSCR_{ij} \cup \lSCR_{ik}$, a contradiction. Therefore, 
        \[
        \left| \bigcap_{k\in J} \lSCR_{ik} \right| \geq |\lSCR_{ij}| - \sum_{k\in J\setminus \{j\}} |\lSCR_{ij}\setminus \lSCR_{ik}| \geq (\Delta-10) - 10(|J|-1) = \Delta-110 > 0.
        \]
        Let us thus fix $z\in \bigcap_{k\in J}\lSCR_{ik}$. Then $z$ is an in-neighbour of $w_i$, and an in-neighbour of $y_k$ for every $k\in J$.  Note that the $y_k$'s are pairwise distinct as we chose $I$ such that $\Tilde{N}_k$ and $\Tilde{N}_{k'}$ are disjoint for all distinct $k,k'\in I$.
        
        Finally, $z$ has at least $\Delta-11$ out-neighbours in $\lSCR_{ij}$. Observe that $\{y_k : k\in J\}$ is disjoint from $\lSCR_{ij}$, as for every $k\in J$, $y_kw_i$ is an arc of $A(D_{ik}) \setminus A(D)$. Altogether, we thus have
        \[
        d^+(z) \geq 1 + |J| + (\Delta-11) = \Delta+1,
        \]
        a contradiction.
    \end{proofclaim}

    We finally bound the size of $I^{> 6}$. This time, since the size of $\Tilde{N}_i$ might be arbitrarily large, we are not able to find a large biclique with many out-neighbours outside, like we did in the two previous claims. However, because each $\Tilde{N}_i$ is indeed large, we will manage to find distinct $i,j,k \in I^{> 6}$ for which $\Tilde{N}_i\cap \Tilde{N}_j \cap \Tilde{N}_k \neq \emptyset$, thus contradicting Claim~\ref{claim:delplus:Ni_almost_disjoint}.
    
    \begin{claimlabeled}[thm:directed_bk_deltaplus]
        \label{claim:delplus:Isup6_bounded}
        $|I^{> 6}| \leq 42$.
    \end{claimlabeled}
    \begin{proofclaim}
        Assume for a contradiction that $|I^{\leq 6}| \geq 43$. Again, by Claim~\ref{claim:delplus:the_wis_are_almost_distinct}, there exists $I\subseteq I^{> 6}$, $|I| = \lceil 43/14\rceil=4$, such that vertices $(w_i)_{i\in I}$ are pairwise distinct.
        
        Again, recall that for every distinct $i,j\in I$, $D_{ij}$ contains an induced copy $K_{ij}$ of either $\bid{K}_{\Delta}$ or $\vec{C_3}\boxplus \bid{K}_{\Delta-2}$. We fix such a $K_{ij}$ for every pair $i,j\in I$. As $K_{ij}$ is not a subgraph of $D$, it contains at least one vertex of $\{w_i,w_j\}$.
        Similarly to the previous case, we thus build an auxiliary tournament $T^I$ with $V(T^I) = I$, and we orient each edge $ij$ from $i$ to $j$ if $K_{ij}$ contains $w_j$, and from $j$ to $i$ otherwise (again, we orient $ij$ arbitrarily if it contains both).
        
        As $T^I$ has order $4$, it contains a vertex $i$ of in-degree at least $2$. We fix such a vertex $i$ and two distinct in-neighbours $j,k$ of $i$ in $T^I$. By definition, both $K_{ij}$ and $K_{ik}$ contain $w_i$.

        Let us denote by $S_i^+$ the set of simple out-neighbours of $w_i$. Observe that, in $K_{ij}$, $w_i$ is linked with digons to at least $\Delta-2$ vertices in $V(K_{ij})$, let us denote them by $Y_j$. By construction, every vertex $y\in Y_j$ is an out-neighbour of $w_i$ in $D$. Therefore, $Y_j$ contains all but at most two simple out-neighbours of $w_i$. These vertices from $Y_j\cap S_i^+$ necessarily belong to $\Tilde{N}_j$, as they are simple out-neighbours of $w_i$ in $D$, but linked with digons to $w_i$ in $D_{ij}$. Hence $|S_i^+\setminus \Tilde{N}_j| \leq 2$. Similarly, $|S_i^+\setminus \Tilde{N}_k| \leq 2$, and by Claim~\ref{claim:delplus:Ni_closeto_Si}, $|S_i^+\setminus \Tilde{N}_i| \leq 2$.
        Therefore, since $|S_i^+|\geq 7$ by definition of $I^{>6}$, we have
        \[ 
        |S_i^+ \cap \tilde{N}_i \cap \tilde{N}_j \cap \tilde{N}_k | \geq |S_i^+|-|(S_i^+ \setminus \Tilde{N}_i) \cup (S_i^+ \setminus \Tilde{N}_j) \cup (S_i^+ \setminus \Tilde{N}_k)| \geq 1,
        \]
        and in particular $\tilde{N}_i \cap \tilde{N}_j \cap \tilde{N}_k \neq \emptyset$, a contradiction to Claim~\ref{claim:delplus:Ni_almost_disjoint}.
    \end{proofclaim}
    
    By Claims~\ref{claim:delplus:Idig_bounded},~\ref{claim:delplus:Ileq6_bounded}~and~\ref{claim:delplus:Isup6_bounded}, and the fact that $(I^{\pm},I^{\leq 6 },I^{>6})$ is a partition of $[\Delta-2]$, we have $\Delta\leq 2 + 70 + 441 + 42 =555$, a contradiction.
\end{proof}

\section{A sufficient condition for \texorpdfstring{$\delmin$}{Deltamin}}
\label{sec:directed_bk:min}

Let $k$ be an integer. A {\bf $k$-obstruction} is a digraph $H$ of order $\ell\in \{k,k+1\}$ admitting a partition $(A,B)$ of its vertex-set, such that $H$ contains every arc from $A$ to $B$, and such that one of the following holds:
\begin{enumerate}[label=$(\roman*)$]
    \item $\ell = k$ and both $A$ and $B$ are bicliques of $H$,
    \label{enum:k_obstruction:1}
    \item $\ell = k+1$, $A$ is a biclique of $H$, and $H[B]$ is isomorphic to $\Vec{C_3} \boxplus \bid K_{|B|-3}$, or
    \item $\ell = k+1$, $H[A]$ is isomorphic to $\Vec{C_3} \boxplus \bid K_{|A|-3}$, and $B$ is a biclique of $H$.
\end{enumerate}

\begin{figure}
    \begin{center}
        \subfigure[]{
        \begin{tikzpicture}[thick,scale=0.8, every node/.style={transform shape}]
            \node[bigvertex] (K1) at (0,0) {$\bid{K}_{p}$};
            \node[bigvertex] (K2) at (2.4,0) {$\bid{K}_{k-p}$};
            \node[] (K) at (1.2,0) {\LARGE $\Rightarrow$};
          \end{tikzpicture}
        }
        \hspace*{1cm}
        \subfigure[]{
            \begin{tikzpicture}[thick,scale=0.8, every node/.style={transform shape}]
            \node[bigvertex] (K1) at (0,0) {$\bid{K}_{p-2}$};
            \node[bigvertex] (K2) at (-3,1.5) {$\bid{K}_{k-p}$};
            \begin{scope}[yshift=3cm]
                \node[vertex] (x) at (30:1) {};
                \node[vertex] (y) at (150:1) {};
                \node[vertex] (z) at (-90:1) {};
            \end{scope}
            \draw[digon] (x) to[out=-80, in=60] (K1);
            \draw[digon] (y) to[out=-100, in=120] (K1);
            \draw[digon] (z) -- (K1);
            \draw[edge,solid] (x) -- (y);
            \draw[edge,solid] (y) -- (z);
            \draw[edge,solid] (z) -- (x);
            \node[] (K) at (-1.7,1.5) {\LARGE $\Rightarrow$};
            \draw[draw=black, thin,dashed] (-1.2,4) rectangle ++(2.4,-5);
        \end{tikzpicture}
        }
        \hspace*{1cm}
        \subfigure[]{
            \begin{tikzpicture}[thick,scale=0.8, every node/.style={transform shape}]
            \node[bigvertex] (K1) at (0,0) {$\bid{K}_{p-2}$};
            \node[bigvertex] (K2) at (3,1.5) {$\bid{K}_{k-p}$};
            \begin{scope}[yshift=3cm]
                \node[vertex] (x) at (30:1) {};
                \node[vertex] (y) at (150:1) {};
                \node[vertex] (z) at (-90:1) {};
            \end{scope}
            \draw[digon] (x) to[out=-80, in=60] (K1);
            \draw[digon] (y) to[out=-100, in=120] (K1);
            \draw[digon] (z) -- (K1);
            \draw[edge,solid] (x) -- (y);
            \draw[edge,solid] (y) -- (z);
            \draw[edge,solid] (z) -- (x);
            \node[] (K) at (1.7,1.5) {\LARGE $\Rightarrow$};
            \draw[draw=black, thin,dashed] (-1.2,4) rectangle ++(2.4,-5);
        \end{tikzpicture}
        }
        \caption{The three possible types of a $k$-obstruction. In Figure~(a), $p$ is any value between $0$ and $k$. In Figures~(b) and~(c), $p$ is any value between $2$ and $k$. }
          \label{fig:delta_obstructions}
    \end{center}
\end{figure}
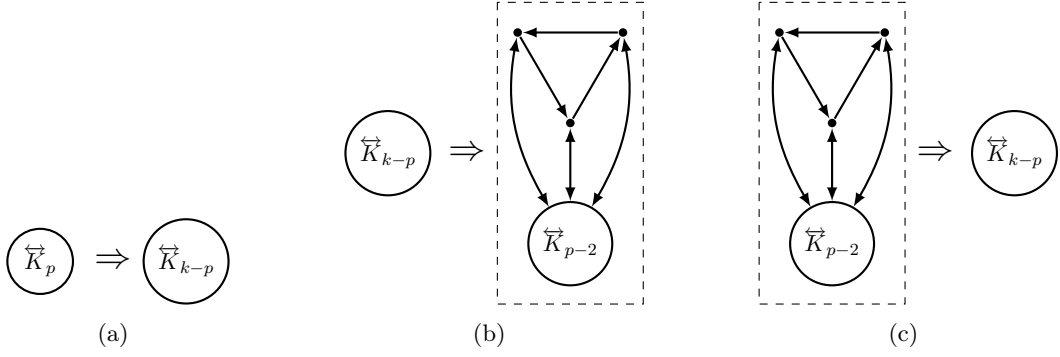

See Figure~\ref{fig:delta_obstructions} for an illustration. Note that $A$ or $B$ might be empty. Observe also that if a symmetric digraph $D$ contains a $k$-obstruction, then $\bic(D) \geq k$. 
Picasarri-Arrieta~\cite{picasarriJGT106} proved that every digraph $D$ satisfying $\dic(D) = \delmin(D)+1$ contains a $(\delmin(D)+1)$-obstruction satisfying~\ref{enum:k_obstruction:1}, thereby generalising Brooks' Theorem.
Similarly, the following is a third independent generalisation of Theorem~\ref{thm:reed_bk} to the directed setting. 

\begin{corollary}
    \label{thm:directed_bk_deltamin}
    The following holds for every integer $\Delta\geq\DeltaDirectedBKPlus$. Every digraph $D$ with $\delmin(D)\leq\Delta$ satisfying $\dic(D) \geq\Delta$ contains a $\Delta$-obstruction.
\end{corollary}
\begin{proof}
    Let $D$ be any digraph with $\delmin(D)\leq\Delta$ and $\dic(D) \geq \Delta$. Let $(A,B)$ be the bipartition of  $V(D)$ defined as:
    \[ 
    A=\{v\in V(D) : d^+(v) \leq \delmin(D)\} \text{~~~and~~~} B=\{v\in V(D) : d^+(v) > \delmin(D)\}.
    \]
    Let $\Tilde{D}$ be the digraph obtained from $D$ by removing all arcs from $B$ to $A$, and replacing all arcs from $A$ to $B$ by digons. Let then $\hat{D}$ be the digraph obtained from $\Tilde{D}$ by reversing every arc of $D[B]$.
    It was proved in~\cite[Proposition~7.3]{kawarabayashiArxiv24} that $\delplus(\hat{D})\leq \delmin(D) \leq \Delta$ and $\dic(\hat{D}) \geq \dic(D)\geq \Delta$, we omit the proof.

    By Theorem~\ref{thm:directed_bk_deltaplus}, $\hat{D}$ contains an induced copy $K$ of $\bid{K}_{\Delta}$ or $\vec{C_3}\boxplus \bid{K}_{\Delta-2}$. By construction, there is no simple arc between $A$ and $B$ in $\hat{D}$. Hence, if $K$ is isomorphic to $\vec{C_3}\boxplus \bid{K}_{\Delta-2}$, its vertices inducing $\vec{C}_3$ all belong to the same side of the bipartition $(A,B)$. 

    Since, by construction, all arcs of $\hat{D}$ going from $A$ to $B$ belong to $D$, we conclude that $D[V(K)]$ is a $\Delta$-obstruction, as desired.
\end{proof}

Noticing that every $k$-obstruction contains a biclique of size $\lceil \frac{1}{2}(k-1)\rceil$, we obtain the following.

\cordirectedbkmin*

Similarly, noticing that the underlying graph of any $k$-obstruction contains a clique of size $k$, we derive the following, which itself generalises Theorem~\ref{thm:reed_bk}.

\cordirectedbkminbis*

As mentioned in the introduction, the condition $\bic(D) < \frac{1}{2}(\Delta-1)$ above is best possible, as we have the following hardness result.

\propnpcompleteness*

\begin{proof} 
    The problem is clearly in NP, so let us show that it is NP-hard.
    The reduction is similar to that of~\cite[Theorem~9]{picasarriJGT106}.
    
    We reduce from the problem of deciding whether a given digraph $D$ has dichromatic number at most $k$, which has been shown to be NP-hard for every fixed integer $k\geq 2$ in~\cite{chenSIC37}.
    Let $D=(V,A)$ be an instance of this problem for a fixed $k\geq 2$.
    
    We build $D'=(V',A')$ from $D$ as follows. For each vertex $u\in V$, we associate a copy $H_u^-$ of $\vec{C_3} \boxplus \bid{K}_{\lceil (k-2)/2\rceil}$, a copy $H_u^+$ of $\bid{K}_{\lfloor k/2\rfloor}$, and we put all arcs from $V(H_u^-)$ to $V(H_u^+)$.
    For each arc $uv\in A$, we associate all possible arcs $xy$ in $A'$, such that $x \in V(H_u^+)$ and $y\in V(H_v^-)$.
    
    For every $u\in V$, vertices in $H_u^-$ have out-degree at most $k+1$, and vertices in $H_u^+$ have in-degree at most $k+1$.
    It is also straightforward that $\bic(D') = \lceil k/2\rceil$. Let us prove that $\dic(D) \leq k$ if and only if $\dic(D') \leq k$ to get the result.

    \medskip
    
    Assume first that  $\dic(D) \leq k$, and let $\phi$ be a $k$-dicolouring of $D$. Let $\phi'$ be the $k$-dicolouring of $D'$ defined as follows. For each vertex $u\in V$, choose arbitrarily $x^-_u \in V(H_u^-)$ and $x^+_u\in V(H_u^+)$, and set $\phi'(x^-_u) =\phi'(x^+_u) = \phi(u)$.
    Then choose arbitrarily two distinct vertices $y_u,z_u \in V(H_u^-)\setminus \{x_u^-\}$ that are not linked with a digon (which is possible as $H_u^-$ contains an induced directed triangle), and colour them with a common colour distinct from $\phi(u)$.
    Finally, choose a distinct colour for every other vertex of $V(H_u^-) \cup V(H_u^+)$.

    Observe that $\phi'$ is be a $k$-dicolouring of $D'$. Indeed, by construction, if $D'$, coloured with $\phi'$, contains a monochromatic directed cycle, then it is of the form $x_{u_1}^-,x_{u_1}^+,x_{u_2}^-,x_{u_2}^+,\ldots,x_{u_\ell}^-,x_{u_\ell}^+,x_{u_1}^-$. Hence, $u_1,u_2,\ldots,u_{\ell},u_1$ is a monochromatic directed cycle of $D$ coloured with $\phi$, a contradiction.

    \medskip
    
    Assume now that $\dic(D') \leq k$.
    Let $\phi' : V' \xrightarrow{} \{1,\dots,k\}$ be a $k$-dicolouring of $D'$. Let $\phi$ be the $k$-dicolouring of $D$ defined as follows. For each vertex $u\in V$, observe that there exist $x_u^- \in V(H_u^-)$ and $x_u^+ \in V(H_u^+)$ such that $\phi'(x_u^-) = \phi'(x_u^+)$. We set $\phi(u) = \phi'(x_u^-)$. Observe that $\phi$ is a $k$-dicolouring of $D$. Indeed, if $D$, coloured with $\phi$, contains a monochromatic directed cycle $u_1,u_2,\ldots,u_\ell,u_1$, then $x_{u_1}^-,x_{u_1}^+,x_{u_2}^-,x_{u_2}^+,\ldots,x_{u_\ell}^-,x_{u_\ell}^+,x_{u_1}^-$ is a monochromatic directed cycle in $D'$ coloured with $\phi'$, a contradiction.
\end{proof}

\bibliographystyle{abbrv}
\bibliography{refs}

\begin{thebibliography}{10}

\bibitem{aboulkerDM113193}
P.~Aboulker and G.~Aubian.
\newblock Four proofs of the {D}irected {B}rooks' {T}heorem.
\newblock {\em Discrete Mathematics}, page 113193, 2022.

\bibitem{aboulker2023digraph}
P.~Aboulker, G.~Aubian, and P.~Charbit.
\newblock Digraph colouring and arc-connectivity.
\newblock {\em preprint arXiv:2304.04690}, 2023.

\bibitem{aboulkerEJC31b}
P.~Aboulker, G.~Aubian, P.~Charbit, and S.~Thomassé.
\newblock ($\overrightarrow{P_6}$, triangle)-free digraphs have bounded dichromatic number.
\newblock {\em The Electronic Journal of Combinatorics}, 31(P4.60), 2024.

\bibitem{aboulkerSIAM36}
P.~Aboulker, G.~Aubian, and R.~Steiner.
\newblock Heroes in orientations of chordal graphs.
\newblock {\em SIAM Journal on Discrete Mathematics}, 36(4):2497--2505, 2022.

\bibitem{aboulkerEJC28}
P.~Aboulker, P.~Charbit, and R.~Naserasr.
\newblock Extension of {G}y{\'a}rf{\'a}s-{S}umner conjecture to digraphs.
\newblock {\em The Electronic Journal of Combinatorics}, 28(P2.27), 2021.

\bibitem{aboulkerEJC31a}
P.~Aboulker and Q.~Vermande.
\newblock Various bounds on the minimum number of arcs in a $k$-dicritical digraph.
\newblock {\em The Electronic Journal of Combinatorics}, 31(P1.22), 2024.

\bibitem{alon2008}
N.~Alon and J.~H. Spencer.
\newblock {\em The probabilistic method}.
\newblock Wiley Series in Discrete Mathematics and Optimization. Wiley-Blackwell, Hoboken, NJ, 3rd edition, 2008.

\bibitem{andresJGT79}
S.~D. Andres and W.~Hochst{\"a}ttler.
\newblock Perfect digraphs.
\newblock {\em Journal of Graph Theory}, 79(1):21--29, 2015.

\bibitem{axenovichCPC31}
M.~Axenovich, A.~Gir{\~a}o, R.~Snyder, and L.~Weber.
\newblock Strong complete minors in digraphs.
\newblock {\em Combinatorics, Probability and Computing}, 31(3):489--506, 2022.

\bibitem{azumaTMJ19}
K.~Azuma.
\newblock Weighted sums of certain dependent random variables.
\newblock {\em Tohoku Mathematical Journal, Second Series}, 19(3):357--367, 1967.

\bibitem{bang2009}
J.~Bang-Jensen and G.~Z. Gutin.
\newblock {\em {Digraphs: Theory, Algorithms and Applications}}.
\newblock Springer-Verlag, London, 2nd edition, 2009.

\bibitem{bokalJGT46}
D.~Bokal, G.~Fijavz, M.~Juvan, P.~M. Kayll, and B.~Mohar.
\newblock The circular chromatic number of a digraph.
\newblock {\em Journal of Graph Theory}, 46(3):227--240, 2004.

\bibitem{borodinJCTB23}
O.~V. Borodin and A.~V. Kostochka.
\newblock On an upper bound of a graph's chromatic number, depending on the graph's degree and density.
\newblock {\em Journal of Combinatorial Theory, Series B}, 23(2-3):247--250, 1977.

\bibitem{brooksMPCPS37}
R.~L. Brooks.
\newblock On colouring the nodes of a network.
\newblock {\em Mathematical Proceedings of the Cambridge Philosophical Society}, 37(2):194–197, 1941.

\bibitem{chenSIC37}
X.~Chen, X.~Hu, and W.~Zang.
\newblock A min-max theorem on tournaments.
\newblock {\em SIAM Journal on Computing}, 37(3):923–937, 2007.

\bibitem{cookEJC30}
L.~Cook, T.~Masařík, M.~Pilipczuk, A.~Reinald, and U.~S. Souza.
\newblock Proving a directed analogue of the {G}yárfás-{S}umner conjecture for orientations of {$P_4$}.
\newblock {\em The Electronic Journal of Combinatorics}, 30(P3.36), 2023.

\bibitem{erdosPNCN1979}
P.~Erd\H{o}s.
\newblock Problems and results in number theory and graph theory.
\newblock In {\em Proceedings of the ninth Manitoba Conference on Numerical Mathematics and Computing}, pages 3--21, 1979.

\bibitem{erdosIFS10}
P.~Erd\H{o}s and L.~Lov{\'a}sz.
\newblock Problems and results on 3-chromatic hypergraphs and some related questions.
\newblock {\em Infinite and Finite Sets}, 10(2):609--627, 1975.

\bibitem{farzadJCTB93}
B.~Farzad, M.~Molloy, and B.~Reed.
\newblock {($\Delta-k$)-critical graphs}.
\newblock {\em Journal of Combinatorial Theory, Series B}, 93(2):173--185, 2005.

\bibitem{golowichDM339}
N.~Golowich.
\newblock The m-degenerate chromatic number of a digraph.
\newblock {\em Discrete Mathematics}, 339(6):1734--1743, 2016.

\bibitem{gonccalves2024brooks}
D.~Gon{\c{c}}alves, L.~Picasarri-Arrieta, and A.~Reinald.
\newblock Brooks-type colourings of digraphs in linear time.
\newblock {\em preprint arXiv:2405.05222}, 2024.

\bibitem{harutyunyanSIDMA25}
A.~Harutyunyan and B.~Mohar.
\newblock Gallai's theorem for list coloring of digraphs.
\newblock {\em SIAM Journal on Discrete Mathematics}, 25(1):170--180, 2011.

\bibitem{harutyunyanEJC18}
A.~Harutyunyan and B.~Mohar.
\newblock Strengthened {B}rooks{\textquotesingle} theorem for digraphs of girth at least three.
\newblock {\em The Electronic Journal of Combinatorics}, 18(P195), 2011.

\bibitem{havetJGT107}
F.~Havet, L.~Picasarri-Arrieta, and C.~Rambaud.
\newblock On the minimum number of arcs in 4-dicritical oriented graphs.
\newblock {\em Journal of Graph Theory}, 107(4):778--809, 2024.

\bibitem{kawarabayashiArxiv24}
K.~Kawarabayashi and L.~Picasarri-Arrieta.
\newblock An analogue of {R}eed’s conjecture for digraphs.
\newblock {\em preprint arXiv:2407.05827}, 2024.

\bibitem{kawarabayashiSODA25}
K.~Kawarabayashi and L.~Picasarri-Arrieta.
\newblock An analogue of {R}eed’s conjecture for digraphs.
\newblock In {\em Proceedings of the 2025 Annual ACM-SIAM Symposium on Discrete Algorithms (SODA)}, pages 3310--3324. SIAM, 2025.

\bibitem{kostochkaGC36}
A.~V. Kostochka and M.~Stiebitz.
\newblock The minimum number of edges in 4-critical digraphs of given order.
\newblock {\em Graphs and Combinatorics}, 36(3):703--718, 2020.

\bibitem{meszarosJGT101}
T.~M{\'e}sz{\'a}ros and R.~Steiner.
\newblock Complete directed minors and chromatic number.
\newblock {\em Journal of Graph Theory}, 101(4):623--632, 2022.

\bibitem{moharJGT43}
B.~Mohar.
\newblock Circular colorings of edge-weighted graphs.
\newblock {\em Journal of Graph Theory}, 43(2):107--116, 2003.

\bibitem{moharLAA432}
B.~Mohar.
\newblock Eigenvalues and colorings of digraphs.
\newblock {\em Linear Algebra and its Applications}, 432(9):2273--2277, 2010.

\bibitem{molloyComb18}
M.~Molloy and B.~Reed.
\newblock A bound on the total chromatic number.
\newblock {\em Combinatorica}, 18(2):241--280, 1998.

\bibitem{molloyCPC7}
M.~Molloy and B.~Reed.
\newblock The size of the giant component of a random graph with a given degree sequence.
\newblock {\em Combinatorics, probability and computing}, 7(3):295--305, 1998.

\bibitem{molloyreed}
M.~Molloy and B.~Reed.
\newblock {\em Graph colouring and the probabilistic method}.
\newblock Algorithms and Combinatorics. Springer, Berlin, Germany, Nov. 2001.

\bibitem{molloyJCTB109}
M.~Molloy and B.~Reed.
\newblock Colouring graphs when the number of colours is almost the maximum degree.
\newblock {\em Journal of Combinatorial Theory, Series B}, 109:134--195, 2014.

\bibitem{neumannlaraJCT33}
V.~Neumann-Lara.
\newblock {The dichromatic number of a digraph}.
\newblock {\em Journal of Combinatorial Theory, Series B}, 33:265--270, 1982.

\bibitem{picasarriJGT106}
L.~Picasarri-Arrieta.
\newblock Strengthening the {D}irected {B}rooks' {T}heorem for oriented graphs and consequences on digraph redicolouring.
\newblock {\em Journal of Graph Theory}, 106(1):5--22, 2024.

\bibitem{picasarriDM347}
L.~Picasarri-Arrieta and M.~Stiebitz.
\newblock {Minimum number of arcs in $k$-critical digraphs with order at most $2k-1$}.
\newblock {\em Discrete Mathematics}, 347(9):114072, 2024.

\bibitem{reedJGT27}
B.~Reed.
\newblock {$\omega$, $\Delta$, and $\chi$}.
\newblock {\em Journal of Graph Theory}, 27(4):177--212, 1998.

\bibitem{reedJCTB76}
B.~Reed.
\newblock {A strengthening of Brooks' theorem}.
\newblock {\em Journal of Combinatorial Theory, Series B}, 76(2):136--149, 1999.

\bibitem{talagrand1995}
M.~Talagrand.
\newblock Concentration of measure and isoperimetric inequalities in product spaces.
\newblock {\em Publications Math{\'e}matiques de l'Institut des Hautes Etudes Scientifiques}, 81:73--205, 1995.

\end{thebibliography}

\clearpage

\appendix

\section{Proof of Lemma~\ref{lemma:talagrand}}
\label{sec:appendix}

As explained in Section~\ref{sec:notations}, we show in this appendix how Lemma~\ref{lemma:talagrand} is derived from the original result of Talagrand. The proof is identical to that of~\cite{molloyJCTB109}, applied to the particular case of $X$ being integer valued. Let us first formally state Talagrand's original result.

Let $T_1,\ldots,T_n$ be $n$ independent trials and $\mathcal{A}$ be the set of all sequences of possible outcomes of these. For any $A \subseteq \mathcal{A}$ and any $\ell\geq 0$, we let $A_\ell$ be the set of elements $b\in \mathcal{A}$ which are at distance at most $\ell$ from $A$ in the following sense. We have $b=(b_1,\dots,b_n)\in A_\ell$ if and only if,  for any sequence of real numbers $\lambda_1,\dots,\lambda_n$, there exists $a=(a_1,\dots,a_n)\in A$ such that
\[
    \sum_{\substack{i\in [n]\\a_i\neq b_i}} \lambda_i \leq \ell \cdot \Big(\sum_{i\in [n]} \lambda_i^2\Big)^{1/2}.
\]
\begin{theorem}[Talagrand~\cite{talagrand1995}]
    \label{thm:talagrand}
    For any $n$ independent trials $T_1,\dots,T_n$, any $A\subseteq \mathcal{A}$, and any real $\ell\geq 0$, we have 
    \[
    \PP(A) \cdot \PP(\overline{A_\ell}) \leq e^{-\ell^2/4}.
    \]
\end{theorem}

We now prove Lemma~\ref{lemma:talagrand}, which we first recall here for convenience.

\talagrand*

\begin{proof}
    We proceed in two steps. We first show that $X$ is concentrated around its median, and then show that its median is close to its expectation. Let $\mu$ be any integer such that $\PP(X\leq \mu) \geq \frac{1}{2}$ and $\PP(X\geq \mu) \geq \frac{1}{2}$.

    \begin{claimlabeled}[lemma:talagrand]
        \label{claim:talagrand:1}
        For every $t\geq 0$, $\PP(X< \mu- t) \leq 2\exp\left(\frac{-t^2}{4c^2r\mu}\right)$.
    \end{claimlabeled}
    \begin{proofclaim}
        Let $A = \{a\in \mathcal{A} : X(a) \geq \mu\}$, $B = \{ b\in \mathcal{A} : X(b) < \mu - t\}$, and $\ell= \frac{t}{c\sqrt{r\mu}}$. We first justify that $A\subseteq \overline{B_\ell}$. To see this, let us fix any $a=(a_1,\dots,a_n)\in A$. We will exhibit a specific sequence of real numbers $\lambda_1,\dots,\lambda_n$ such that, for every $b=(b_1,\dots,b_n)\in B$,
        \[
            \sum_{\substack{i\in [n]\\a_i\neq b_i}} \lambda_i > \ell ~\cdot \Big(\sum_{i\in [n]} \lambda_i^2\Big)^{1/2},
        \]
        hence implying $A\subseteq \overline{B_\ell}$.
        
        By definition of $A$, $X(a)\geq \mu$, and since $\mu$ is an integer, by assumption there exists a set $I$ of at most $r\mu$ indices such that the values of $(a_i)_{i\in I}$ certify $X(a)\geq \mu$. For each $i\in [n]$, we set $\lambda_i = c$ if $i\in I$, and $\lambda_i = 0$ otherwise. Let us fix any $b\in B$, and let $b'\in \mathcal{A}$ be the sequence equal to $a$ on $I$ and equal to $b$ on $[n]\setminus I$.
        By construction and by the first assumption of the statement, we have 
        \begin{align*}
            X(b) \geq X(b') - c \cdot |\{i\in I : a_i \neq b_i\}| \geq \mu - \sum_{\substack{i\in [n]\\a_i\neq b_i}} \lambda_i.
        \end{align*}
        In the second inequality, we used that the outcomes of $b'$ on $I$ certify $X(b')\geq \mu$. By definition of $B$, we thus have
        \[
        \sum_{\substack{i\in [n]\\a_i\neq b_i}} \lambda_i > t \geq \ell \cdot \Big(\sum_{i\in [n]}\lambda_i^2\Big)^{1/2},
        \]
        where the second inequality follows from
        \[
            \sum_{i\in [n]} \lambda_i^2 = c^2 \cdot |I| \leq c^2\cdot r\cdot \mu = \frac{t^2}{\ell^2}.
        \]
        This shows $A\subseteq \overline{B_\ell}$. By definition of $A$ and $\mu$, we thus have $\PP(\overline{B_\ell}) \geq \PP(A) \geq \frac{1}{2}$. By Theorem~\ref{thm:talagrand}, we thus have 
        \[
        \PP(X<\mu-t) = \PP(B) \leq 2e^{-\ell^2/4} = 2\exp\left(\frac{-t^2}{4c^2r\mu}\right),
        \]
        as desired.
    \end{proofclaim}

    \begin{claimlabeled}[lemma:talagrand]
        \label{claim:talagrand:2}
        For every $t\geq0$, $\PP(X> \mu + t) \leq 2\exp\left(\frac{-t^2}{4c^2r(\mu+t+1)}\right)$.
    \end{claimlabeled}
    \begin{proofclaim}
        The proof is analogous to that of Claim~\ref{claim:talagrand:1}, with some additional rounding arguments. Therefore, let $A=\{a\in \mathcal{A} : X(a) \leq \mu\}$, $B=\{b\in \mathcal{A}:X(b)>\mu +t\}$, $t'=\lfloor t\rfloor +1$, and $\ell = \frac{t}{c\sqrt{r(\mu + t + 1)}}$. Since $X$ is integer-valued, note that $B$ is also $\{b\in \mathcal{A}:X(b) \geq \mu + t'\}$.

        We now justify that $B \subseteq \overline{A_\ell}$. Let us fix $b=(b_1,\dots,b_n)\in B$. Since $X(b)\geq \mu+t'$, and because $\mu+t'$ is an integer, by assumption there exists a set $I\subseteq [n]$ of at most $r(\mu+t')$ indices such that $(b_i)_{i\in I}$ guarantees $X(b)\geq \mu+t'$. For each $i\in [n]$, we define $\lambda_i = c$ if $i\in I$ and $\lambda_i=0$ otherwise.

        Let $a=(a_1,\dots,a_n)$ be any element of $A$, and let $a'\in \mathcal{A}$ be the sequence equal to $b$ on $I$ and to $a$ on $[n]\setminus I$. We have
        \begin{align*}
            X(a) \geq X(a') - c\cdot |\{i\in I: a_i\neq b_i\}| \geq \mu+t'-\sum_{\substack{i\in [n]\\a_i\neq b_i}} \lambda_i.
        \end{align*}
        By definition of $A$, this implies
        \[ \sum_{\substack{i\in [n]\\a_i\neq b_i}} \lambda_i \geq t' > \ell \cdot \Big(\sum_{i\in [n]}\lambda_i^2\Big)^{1/2},
        \]
        where the second inequality follows from
        \[
            \ell^2\cdot \sum_{i\in [n]} \lambda_i^2 = \ell^2 \cdot c^2 \cdot |I| \leq \ell^2 \cdot c^2\cdot r\cdot (\mu+t') < {t'}^2.
        \]
        This shows that $B\subseteq \overline{A_\ell}$. Recall that $\PP(A) \geq \frac{1}{2}$ by definition of $\mu$, so by Theorem~\ref{thm:talagrand} we obtain
        \[
        \PP(X>\mu+t) = \PP(B) \leq \PP(\overline{A_\ell}) \leq 2e^{-\ell^2/4} \leq 2\exp\left(\frac{-t^2}{4c^2r(\mu+t+1)}\right),
        \]
        as desired.
    \end{proofclaim}

    By the union bound, it follows from Claims~\ref{claim:talagrand:1} and~\ref{claim:talagrand:2} that, for every $t\geq 0$,
    \begin{equation}
        \PP(|X-\mu| >t) \leq 4\exp\left( \frac{-t^2}{4c^2r(\mu+t+1)} \right).
        \label{eq:talagrand}
    \end{equation}
    It remains to justify that $\mu$ is close to $\EE(X)$.
    \begin{claimlabeled}[lemma:talagrand]
        \label{claim:talagrand:3}
        $|\EE(X) -\mu| \leq 63c\sqrt{r\EE(X)} + 172c^2r$.
    \end{claimlabeled}
    \begin{proofclaim}
        We first assume that $\mu \neq 0$.
        For every $i\in \mathbb{N}$, we let $I_i\subseteq \mathbb{R}$ be the interval 
        \[
            I_i = \Big(i \cdot c\sqrt{r\mu},\; (i+1)\cdot c\sqrt{r\mu}\Big].
        \]
        Observe that $|\EE(X)- \mu| = |\EE(X-\mu)|$, and that
        \begin{align*}
            |\EE(X-\mu)| &\leq \sum_{i\in \mathbb{N}} (i+1)\cdot c\sqrt{r\mu} \cdot \PP(|X-\mu| \in I_i)\\
            &= c\sqrt{r\mu} \sum_{i\in \mathbb{N}}\PP\left(|X-\mu| \in \bigcup_{j\geq i} I_j\right)\\
            &= c\sqrt{r\mu} \sum_{i\in \mathbb{N}}\PP\left(|X-\mu| > i\cdot c\sqrt{r\mu}\right)\\
            &\leq 4c\sqrt{r\mu} \sum_{i\in \mathbb{N}} \exp\left(\frac{-i^2 \mu}{4(\mu + ic\sqrt{r\mu}+1)}\right),
        \end{align*}
        where the last inequality follows from~\eqref{eq:talagrand}. This implies
        \begin{align*}
            \frac{1}{4c\sqrt{r\mu}}|\EE(X-\mu)| &\leq  \sum_{i\in \mathbb{N}} \exp\left(\frac{-i^2 \mu}{8\cdot \max(\mu , ic\sqrt{r\mu}+1)}\right)\\
            &\leq \sum_{i\in \mathbb{N}} \exp\left(\frac{-i^2}{8}\right)+\sum_{i\in \mathbb{N}}\exp\left(\frac{-i^2 \mu}{8( ic\sqrt{r\mu}+1)}\right)\\
            &\leq  \sum_{i\in \mathbb{N}} \exp\left(\frac{-i^2}{8}\right)+\sum_{i\in \mathbb{N}}\exp\left(\frac{-i \sqrt{\mu}}{16c\sqrt{r}}\right).
        \end{align*}

        Observe first that $\sum_{i\in \mathbb{N}} \exp\left(\frac{-i^2}{8}\right) \leq \sum_{i\in \mathbb{N}} \exp\left(\frac{-i}{8}\right) = \frac{1}{1-e^{-1/8}} \leq 9$.
        
        It remains to bound the second summation. Let $y=\frac{\sqrt{\mu}}{16c\sqrt{r}}$, and observe that 
        \[
        S = \sum_{i\in \mathbb{N}}\exp\left(\frac{-i \sqrt{\mu}}{16c\sqrt{r}}\right) = \frac{1}{1-e^{-y}}.
        \]
        Therefore, if $y\geq 1$ then $S\leq (1-e^{-1})^{-1} < 2$. Else if $y\leq 1$ then $e^{-y} \leq 1-\frac{1}{2}y$, and $S \leq \frac{2}{y}$. We conclude that 
        \begin{equation}
            \label{eq:talagrand:2}
            |\EE(X-\mu)| \leq  44c\sqrt{r\mu} + 128c^2r.    
        \end{equation}
        Since $X$ is non-negative, $\mu \neq 0$, and $\PP(X\geq \mu)\geq \frac{1}{2}$, by Markov's Inequality we have $\EE(X) \geq \frac{1}{2}\mu$, which implies
        \begin{equation*}
            |\EE(X) - \mu| = |\EE(X-\mu)| \leq  63c\sqrt{r\EE(X)} + 128c^2r.            
        \end{equation*}

        This proves the claim when $\mu\neq 0$. Henceforth, we assume that $\mu=0$, that is $\PP(X=0) \geq \frac{1}{2}$. Consider the random variable $Y=X+1$. Observe that all assumptions on $X$ hold for $Y$, and that $\PP(Y\leq 1) \geq \frac{1}{2}$, $\PP(Y\geq 1)=1$. Hence inequality~\eqref{eq:talagrand:2} holds for $Y$ and gives
        \[
            |\EE(Y-1)| \leq  44c\sqrt{r} + 128c^2r \leq 172c^2r,  
        \]
        and the result follows as $\EE(X) = \EE(Y-1)$.
    \end{proofclaim}

    We are now ready to prove the statement. Let us fix any real number $t>126c\sqrt{r\EE(X)} + 344 c^2 r$. We have
    \begin{align*}
        \PP(|X-\EE(X)| > t) &\leq \PP\left(|X-\EE(X)| > \frac{1}{2}t+63c\sqrt{r\EE(X)}+172c^2r\right) &\text{by choice of $t$}\\
        &\leq \PP(|X-\mu|>\frac{1}{2}t)&\text{by Claim~\ref{claim:talagrand:3}}\\
        &\leq 4\exp\left( \frac{-t^2}{16c^2r(\mu+2t)} \right) &\text{by~\eqref{eq:talagrand} and $t\geq1$}\\
        &\leq 4\exp\left( \frac{-t^2}{32c^2r(\EE(X)+t)} \right)
        &\text{\!by Markov's Inequality.}
    \end{align*}
    The result follows.
\end{proof}

\end{document}